\documentclass{gtart}
\usepackage{amsmath, amssymb, graphicx, epsfig, verbatim, graphs, psfrag}
\usepackage{placeins}
\usepackage{epic}
\usepackage{amscd
}\usepackage{amssymb}
\newcommand{\C}{{\mathbb C}}
\usepackage[curve]{xypic}
\renewcommand{\comment}[1]{}
\newbox\mybox
\def\overtag#1#2#3{\setbox\mybox\hbox{$#1$}\hbox to
  0pt{\vbox to 0pt{\vglue-#3\vglue-\ht\mybox\hbox to \wd\mybox
      {\hss$\scs#2$\hss}\vss}\hss}\box\mybox}
\def\undertag#1#2#3{\setbox\mybox\hbox{$#1$}\hbox to 0pt{\vbox to
    0pt{\vglue#3\vglue\ht\mybox\hbox to \wd\mybox
      {\hss$\scs#2$\hss}\vss}\hss}\box\mybox}
\def\lefttag#1#2#3{\hbox to 0pt{\vbox to 0pt{\vss\hbox to
      0pt{\hss$\scs#2$\hskip#3}\vss}}#1}
\def\righttag#1#2#3{\hbox to 0pt{\vbox to 0pt{\vss\hbox to
      0pt{\hskip#3$\scs#2$\hss}\vss}}#1}
\let\scs\scriptstyle

\def\Dot{\lower.2pc\hbox to 2.5pt{\hss$\bullet$\hss}}
\def\Circ{\lower.2pc\hbox to 2.5pt{\hss$\circ$\hss}}
\def\Vdots{\raise5pt\hbox{$\vdots$}}
\def\splicediag#1#2{\xymatrix@R=#1pt@C=#2pt@M=0pt@W=0pt@H=0pt}

\renewcommand\frame[2][3pt]{\hbox{$\vcenter{\hbox{\vrule\vbox
{\hrule\kern#1\hbox{\kern#1$#2$\kern#1}\kern#1\hrule}\vrule}}$}}
\newcommand\lineto{\ar@{-}}
\newcommand\dashto{\ar@{--}}
\newcommand\dotto{\ar@{.}}

\newtheorem{thm}{Theorem}[section]

\newtheorem{cor}[thm]{Corollary}
\newtheorem{lem}[thm]{Lemma}
\newtheorem{prop}[thm]{Proposition}
\newtheorem{defn}[thm]{Definition}

\newtheorem{exa}[thm]{Example}
\newtheorem{rem}[thm]{Remark}

\numberwithin{equation}{section}

\newcommand{\bt}{\bullet}

\newcommand{\bfz}{{\mathbb {Z}}}
\newcommand{\bfn}{{\mathbb {N}}}
\newcommand{\bfq}{{\mathbb {Q}}}

\newcommand{\rank}{\mbox{rk}}

\newcommand{\s}{\mathbf s}

\newcommand{\cpkk}{{\overline {{\mathbb C}{\mathbb P}^2}}}
\newcommand{\cpk}{{\mathbb {CP}}^2}
\newcommand{\cpot}{{\mathbb {CP}}^2\# 5{\overline {{\mathbb C}{\mathbb P}^2}}}
\newcommand{\cphat}{{\mathbb {CP}}^2\# 6{\overline {{\mathbb C}{\mathbb P}^2}}}

\newcommand{\mg}{M_{\Gamma}}
\newcommand{\xg}{X_{\Gamma}}
\newcommand{\bg}{B_{\Gamma}}
\newcommand{\frg}{{\mathcal {G}}}

\newcommand{\fra}{{\mathcal {A}}}
\newcommand{\frb}{{\mathcal {B}}}
\newcommand{\frc}{{\mathcal {C}}}
\newcommand{\frw}{{\mathcal {W}}}
\newcommand{\frs}{{\mathcal {S}}}
\newcommand{\frn}{{\mathcal {N}}}
\newcommand{\farm}{{\mathcal {M}}}

\begin{document}

\title[Rational blow--downs and smoothings]{Rational blow--downs and
  smoothings of surface singularities}

\author{Andr\'as I. Stipsicz}
\address{R{\'e}nyi Institute of Mathematics\\
Hungarian Academy of Sciences\\
H-1053 Budapest\\ 
Re{\'a}ltanoda utca 13--15, Hungary}
\email{stipsicz@math-inst.hu}

\author{Zolt\'an Szab\'o}
\address{Department of Mathematics,\\
Princeton University,\\
Princeton, NJ, 08544}
\email{szabo@math.princeton.edu}

\author{Jonathan Wahl}
\address{Department of Mathematics,\\
The University of North Carolina,\\ 
Chapel Hill, NC 27599--3250}
\email{jmwahl@email.unc.edu}


\begin{abstract}
  In this paper we give a necessary combinatorial condition for a
  negative--definite plumbing tree to be suitable for rational
  blow--down, or to be the graph of a complex surface singularity
  which admits a rational homology disk smoothing. New examples of
  surface singularities with rational homology disk smoothings are
  also presented; these include singularities with resolution graph
  having valency four nodes.
\end{abstract}
\maketitle

\section{Introduction}
One of the basic operations in algebraic geometry is the blow--up
process, which in complex dimension 2 (and for a point) is simply the
replacement of a point with a rational curve of homological square
$-1$.  From the differential topological point of view this operation
corresponds to replacing a closed tubular neighbourhood of a point
(which is simply the closed 4--disk) by the tubular neighbourhood of a
sphere with self--intersection $-1$.  The inverse operation, i.e., the
replacement of the tubular neighbourhood of a $(-1)$--sphere with the
4--disk, is called the \emph{blow--down process}. It is well--known
how these operations affect topological and complex analytic
invariants; and the effect on Seiberg--Witten invariants is the
blow--up formula of \cite{FS2}.

This
operation has been generalized by Fintushel and Stern \cite{FS1}: in the
\emph{rational blow--down operation} the (closed) tubular neighbourhood of
certain chains of 2--spheres (whose self--intersections are the negatives of
the continued fraction coefficients of $p^2/(p-1)$, $p\geq 2$ an integer) have
been replaced by 4--manifolds with boundary, with rational homology equal to
the rational homology of the 4--disk. This idea was extended by J. Park
\cite{Pratb}, using the linear chains of spheres arising from the continued
fraction coefficients of $p^2/(pq-1)$, for $p>q>0$ relatively prime.  (We will
denote such a plumbing chain by $\Gamma _{p/q}$.)  The usefulness of this
operation stems from the fact that the Seiberg--Witten invariants of the
4--manifold we get by rationally blowing down these chains of
2--spheres can be computed using a fairly simple formula from the
Seiberg--Witten invariants of the original manifold \cite{FS1}, \cite{Pratb}. 
This scheme admits many applications in finding 4--manifolds with
various properties, for example in searching for exotic smooth
structures on 4--manifolds with small Euler characteristic
\cite{FSexo}, \cite{P}, \cite{PSS}, \cite{SS}.

The boundary of a tubular neighbourhood of such a chain of spheres is a
lens space, and it had already been known that a lens space $L(p^2,
pq-1)$ bounds a rational homology disk, by constructions of Casson and
Harer \cite{CH}.  Even more, it had been known that the
two--dimensional cyclic quotient singularity determined by $p^2/(pq-1)$
(whose resolution graph is $\Gamma_{p/q}$) admitted a smoothing with
vanishing Milnor number, hence its link bounds a rational homology
disk (the "Milnor fiber") admitting a Stein structure (see
\cite[Example~(5.9.1)]{Wahl} and Sections~\ref{s:kirby} and
\ref{s:smoothing} below).  Of course, one cannot expect to be able to
do complex--analytic surgery, and replace one compact complex surface
by another, cf. \cite{LPark}, \cite{upa} for related discussion.

The simplicity of the formula relating the Seiberg--Witten invariants
of the 4--manifolds before and after rational blow--down follows
from three facts:
\begin{enumerate}
\item if the given configuration of spheres is tautly embedded for a given
  spin$^c$ structure $\s $ (see Subsection~\ref{ss:taut}), then the restriction
  of the spin$^c$ structure extends from the complement of the configuration
  of spheres to the rational homology disk;
\item this extended spin$^c$ structure gives rise to a Seiberg--Witten moduli
  space with the same expected dimension as $\s$; and finally
\item the 3--manifold along which we cut and glue admits the simplest
possible Seiberg--Witten Floer homology (since it is a lens space, possessing
a metric of positive scalar curvature).
\end{enumerate}

It has been shown by Symington \cite{Sym0}, \cite{Sym} that a configuration
$\Gamma _{p/q}$ of \emph{symplectic} spheres in a symplectic
4--manifold can be blown down symplectically; that is, the symplectic
form, when restricted to the complement of the spheres, extends to the
rational homology disk.  This observation explains the appearance of
properties (1) and (2) above. The existence of the appropriate
symplectic structure on the rational homology disk is a consequence of
the fact that it is the Milnor fiber of a smoothing of the
corresponding quotient singularity, hence admits a Stein structure. It
is known that a surface singularity admitting a rational homology disk
smoothing is necessarily rational (see also
Subsection~\ref{ss:singularities}), and the link of a rational
singularity is an $L$--space \cite{nemethi} (i.e., has the simplest
possible Heegaard Floer homology, which is conjectured to be
equivalent to the link admitting simple Seiberg--Witten Floer
homology). This observation explains (3) above.

After the success of the rational blow--down process in constructing
new and interesting smooth and symplectic 4--manifolds (and more
recently complex surfaces \cite{LPark}, \cite{upa}), it was natural to
ask which other plumbing trees possess similar properties. Casson and
Harer~\cite{CH} provided many examples of 3--manifolds which bound
rational homology disks; but when performing the rational blow--down
process along them, the lack of (some of) the properties (1--3) listed
above implied that the resulting 4--manifolds were usually
uninteresting (i.e., had trivial Seiberg--Witten invariants).
Correspondingly, it has been known for some time that the only cyclic
quotient singularities admitting smoothings with Milnor number 0
(i.e., rational homology disk smoothings) are those with graph
$\Gamma_{p/q}$, see \cite[Remark~(5.10)]{LW}.  On the other hand, a
triply--infinite family of singularities admitting Milnor number 0
smoothings appeared already in 1980 \cite[(5.9.2)]{Wahl}.

We will see that minimal negative--definite plumbing graphs $\Gamma$
of interest in these problems satisfy very strong combinatorial
restrictions.  Properties (1) and (2) above lead to the following

\begin{defn}\label{d:symptree}
{\rm The plumbing tree $\Gamma$ on $n$ vertices is a \emph{symplectic
plumbing tree} if $\Gamma$ admits an embedding $\varphi$ into the
negative--definite diagonal lattice $(\bfz ^n, Q_n)$ with $Q_n=n\langle -1
\rangle$ such that
\begin{itemize}
\item for vertices $v_1\neq v_2\in \Gamma$ we have $Q_n (\varphi (v_1),
\varphi (v_2))=1$ or $0$ depending on whether $v_1$ and $v_2$ are 
adjacent in $\Gamma$, 
\item $Q_n (\varphi (v), \varphi (v))$ is equal to the 
decoration of $v$ for all $v\in \Gamma$, and
\item with the basis $\{ E_1, \ldots , E_n\}$ of $Q_n$ satisfying $Q_n
(E_i, E_j)=-\delta _{ij}$ and with the notation  
$K=\sum _{i=1}^n E_i$ 
we have the adjunction equality
\begin{equation}\label{e:adjunction}
Q_n (v, K)+Q_n (v,v)=-2
\end{equation}
for every vertex $v$ of $\Gamma$.  
\end{itemize}
A graph $\Gamma$ is called \emph{minimal} if there is no vertex in $\Gamma$
with decoration $(-1)$. We denote the set of minimal, connected symplectic
plumbing trees by $\frs$.}
\end{defn}
\begin{rem}
{\rm The name for the set $\frs$ comes from the fact that the
adjunction equality is naturally satisfied by a tree of
\emph{symplectically} embedded spheres in a symplectic 4--manifold.
Notice that if a resolution of a complex surface singularity has dual
graph given by $\Gamma$, then $K$, the restriction of the first Chern
class of the complex structure, satisfies \eqref{e:adjunction} --- as
a result of the adjunction formula.}
\end{rem}
\begin{rem}\label{lattice}
{\rm From the point of view of \cite{LW}, consider the
negative--definite lattice
$$L=\bigoplus_{v\in \Gamma}\mathbb Z \cdot v$$ generated by the
vertices of $\Gamma$, using intersections $v\cdot v'$ as given by
the graph; one has as well an element $K\in L^*$.  Then
Definition~\ref{d:symptree} requires that $L$ be a sublattice of a
unimodular lattice $L'$ of the same rank, which contains the
characteristic element $K$. As in \cite{LW}, the existence of such
an $L'$ and $K$ is equivalent to the existence of a
self--isotropic subspace of the discriminant quadratic group
$(L^*/L,q)$ (see also Section~\ref{ss:singularities}); however, in this case examples indicate that one would not necessarily  have that $L'$ is \emph{diagonal}, as is further required by 
Definition~\ref{d:symptree} .}
\end{rem}

Our main goal is to determine a necessary combinatorial property
of a plumbing tree $\Gamma$ (and its associated plumbing
4--manifold $M_{\Gamma}$) such that there exists a rational
homology disk $B_{\Gamma}$ (with $\partial M_{\Gamma}\cong
\partial B_{\Gamma}$), and for which one of the following holds: (a) for $M_{\Gamma }\subset (X, \s )$ tautly embedded
into a spin$^c$ 4--manifold with $SW_X(\s )\neq 0$, the new manifold
$(X-M_{\Gamma })\cup B _{\Gamma }$ has nontrivial Seiberg--Witten
invariants; (b) a symplectically embedded $M_{\Gamma }$ can be blown
down in the symplectic category; (c) there is a singularity with
resolution graph $\Gamma$ admitting a rational homology disk
smoothing.  As will be explained in the next section, the graphs in
$\frs$ are the only candidates for these three problems. (We will also
explain the relation among the three problems mentioned above.) In
summary, we can compile the results of the next section into

\begin{prop} \mbox{\bf {(cf. 
Corollaries~\ref{cswbd}, \ref{c2-3}, \ref{Milnor})}} Let $\Gamma$ be a
negative--definite minimal plumbing tree of spheres which can be
tautly embedded and rationally blown down, or can be symplectically
blown down, or gives rise to a complex surface singularity admitting a
rational homology disk smoothing.  Then $\Gamma \in \frs$.
\end{prop}

The main theorem of the present paper (Theorem~\ref{t:main}) gives a
combinatorial description of the set $\frs$.  Aside from the
set $\Gamma_{p/q}$, the graphs consist of 3 triply--infinite families
$\frw, \frn, \farm$, built on the 3 basic examples shown by
Figure~\ref{f:n=4}.
\begin{figure}[ht]
\begin{center}
\setlength{\unitlength}{1mm}
\unitlength=0.7cm
\begin{graph}(16,7.5)(0,-7)
\graphnodesize{0.2}

 \roundnode{m1}(0,0)
 \roundnode{m2}(2,0)
 \roundnode{m3}(4,0)
 \roundnode{m4}(2,-2)

 \roundnode{m5}(12,0)
 \roundnode{m6}(14,0)
 \roundnode{m7}(16,0)
 \roundnode{m8}(14,-2)

 \roundnode{m9}(6,-3)
 \roundnode{m10}(8,-3)
 \roundnode{m11}(10,-3)
 \roundnode{m12}(8,-5)

\edge{m1}{m2}
\edge{m2}{m3}
\edge{m2}{m4}

\edge{m5}{m6}
\edge{m6}{m7}
\edge{m6}{m8}

\edge{m9}{m10}
\edge{m10}{m11}
\edge{m10}{m12}

\freetext(2,-3){(a)}
\freetext(14,-3){(b)}
\freetext(8,-6){(c)}

  \autonodetext{m1}[n]{{\small $-3$}}
  \autonodetext{m2}[n]{{\small $-4$}}
  \autonodetext{m3}[n]{{\small $-3$}}
  \autonodetext{m4}[e]{{\small $-3$}}

  \autonodetext{m5}[n]{{\small $-4$}}
  \autonodetext{m6}[n]{{\small $-3$}}
  \autonodetext{m7}[n]{{\small $-4$}}
  \autonodetext{m8}[e]{{\small $-2$}}

  \autonodetext{m9}[n]{{\small $-2$}}
  \autonodetext{m10}[n]{{\small $-2$}}
  \autonodetext{m11}[n]{{\small $-6$}}
  \autonodetext{m12}[e]{{\small $-3$}}
\end{graph}
\end{center}
\caption{\quad Graphs of the basic examples, giving the families
$\frw, \frn, \farm$} 
\label{f:n=4}
\end{figure}
There are 3 more classes $\fra, \frb, \frc$ built out of variations of
these basic examples, cf. Definition~\ref{d:pelda} below.
     
The other key results of the paper (e.g., Theorems~\ref{t:kirby},
\ref{t:smoothing}) show that many of the graphs in $\frs$ actually do
occur in one of our situations (cf. also \cite{GaS}).  But, not all
$\Gamma \in \frs$ occur; further constraints and a better
understanding of the geometric picture will be addressed in a
subsequent project, see also Subsection~\ref{ss:spec}. In order to
state the main result of the paper, we need a definition.

\begin{defn}\label{d:pelda}
{\rm 
\begin{itemize}
\item For integers $p>q>0$ relatively prime, consider the 
continued fraction
expansion of $p^2/(pq-1)$, i.e.,
\[
p^2/(pq-1)=a_1-\cfrac1{a_2-\cfrac1{\ddots - \cfrac1{a_{k}}}}\ \ ,
\] 
where $a_i\in \bfn , \ a_i\geq 2$.  Figure~\ref{f:gpq} shows the
corresponding negative--definite plumbing tree, denoted by $\Gamma
_{p/q}$. We denote the set of all such graphs by $\frg$.

\begin{figure}[ht]
\begin{center}
\setlength{\unitlength}{1mm}
\unitlength=0.7cm
\begin{graph}(10,4)(0,-1)
\graphnodesize{0.2}

 \roundnode{m1}(3,0)
 \roundnode{m2}(4.5,0)

 \roundnode{m6}(6.5,0)
 \roundnode{m7}(8,0)

\edge{m1}{m2}
\edge{m6}{m7}

  \autonodetext{m1}[n]{{\small $-a_1$}}
  \autonodetext{m2}[n]{{\small $-a_2$}}
  \autonodetext{m6}[n]{{\small $-a_{k-1}$}}
  \autonodetext{m7}[n]{{\small $-a_k$}}

\freetext(5.5,0){{\Large $\cdots$}}

\freetext(1,0)
{$\Gamma _{p/q}=$}

\end{graph}
\end{center}
\caption{\quad The graph $\Gamma _{p/q}$ in $\frg $ }
\label{f:gpq}
\end{figure}

\item 
The plumbing tree given by Figure~\ref{f:wahltype} will be denoted by
$\Gamma _{p,q,r}$ ($p,q,r\geq 0$).
\begin{figure}[ht]
\begin{center}
\setlength{\unitlength}{1mm}
\unitlength=0.7cm
\begin{graph}(10,6)(0,-5)
\graphnodesize{0.2}

 \roundnode{m1}(1.5,0)
 \roundnode{m2}(2.5,0)
 \roundnode{m3}(4,0)
 \roundnode{m4}(5,0)
 \roundnode{m5}(6,0)
 \roundnode{m6}(7.5,0)
 \roundnode{m7}(8.5,0)
 \roundnode{m8}(5,-1)  
 \roundnode{m9}(5,-2.5)
 \roundnode{m10}(5,-3.5)

\edge{m1}{m2}
\edge{m3}{m4}
\edge{m4}{m5}
\edge{m6}{m7}
\edge{m4}{m8}
\edge{m9}{m10}

  \autonodetext{m1}[sw]{{\small $-(p+3)$}}
  \autonodetext{m2}[n]{{\small $-2$}}
  \autonodetext{m3}[n]{{\small $-2$}}
  \autonodetext{m4}[n]{{\small $-4$}}
  \autonodetext{m5}[n]{{\small $-2$}}
  \autonodetext{m6}[n]{{\small $-2$}}
  \autonodetext{m7}[se]{{\small $-(q+3)$}}
  \autonodetext{m8}[w]{{\small $-2$}}
  \autonodetext{m9}[w]{{\small $-2$}}
  \autonodetext{m10}[s]{{\small $-(r+3)$}}
  \autonodetext{m3}[w]{{\Large $\cdots$}}
  \autonodetext{m5}[e]{{\Large $\cdots$}}

\freetext(3.25,-0.8)
{$\underset{{\textstyle q}}{\underbrace{\hspace{1.2cm}}}$}

\freetext(6.75,-0.8)
{$\underset{{\textstyle r}}{\underbrace{\hspace{1.2cm}}}$}

\freetext(5.7,-1.75)
{{\Huge $\rbrace$}}

\freetext(5,-1.4){\Large $.$}
\freetext(5,-1.7){\Large $.$}
\freetext(5,-2){\Large $.$}

\freetext(6.2,-1.75){$p$}

\end{graph}
\end{center}
\caption{\quad The graph $\Gamma _{p,q,r}$ in $\frw$ }
\label{f:wahltype}
\end{figure}
We denote the set of these plumbing trees by $\frw$.

\item The plumbing tree of Figure~\ref{f:masik}
will be denoted by $\Delta _{p,q, r}$ ($p\geq 1$ and $q,r\geq 0$).
\begin{figure}[ht]
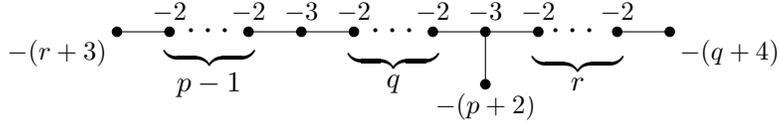

\begin{center}
\setlength{\unitlength}{1mm}
\unitlength=0.7cm
\begin{graph}(6,2)(-1,-1)
\graphnodesize{0.2}

 \roundnode{m12}(-2,0)
 \roundnode{m11}(-1,0)
 \roundnode{m10}(0.5,0)
 \roundnode{m1}(1.5,0)
 \roundnode{m2}(2.5,0)
 \roundnode{m3}(4,0)
 \roundnode{m4}(5,0)
 \roundnode{m5}(6,0)
 \roundnode{m6}(7.5,0)
 \roundnode{m7}(8.5,0)
 \roundnode{m8}(5,-1)  

\edge{m1}{m2}
\edge{m3}{m4}
\edge{m4}{m5}
\edge{m6}{m7}
\edge{m4}{m8}
\edge{m12}{m11}
\edge{m10}{m1}

  \autonodetext{m1}[n]{{\small $-3$}}
  \autonodetext{m2}[n]{{\small $-2$}}
  \autonodetext{m3}[n]{{\small $-2$}}
  \autonodetext{m4}[n]{{\small $-3$}}
  \autonodetext{m5}[n]{{\small $-2$}}
  \autonodetext{m6}[n]{{\small $-2$}}
  \autonodetext{m7}[se]{{\small $-(q+4)$}}
  \autonodetext{m8}[s]{{\small $-(p+2)$}}
 \autonodetext{m10}[n]{{\small $-2$}}
 \autonodetext{m11}[n]{{\small $-2$}}
 \autonodetext{m12}[sw]{{\small $-(r+3)$}} 
  \autonodetext{m3}[w]{{\Large $\cdots$}}
  \autonodetext{m5}[e]{{\Large $\cdots$}}
 \autonodetext{m10}[w]{{\Large $\cdots$}}

\freetext(3.25,-0.8)
{$\underset{{\textstyle q}}{\underbrace{\hspace{1.2cm}}}$}
\freetext(-0.25,-0.8)
{$\underset{{\textstyle p-1}}{\underbrace{\hspace{1.2cm}}}$}

\freetext(6.75,-0.8)
{$\underset{{\textstyle r}}{\underbrace{\hspace{1.2cm}}}$}
\end{graph}
\end{center}
\caption{\quad The graph $\Delta _{p,q,r}$ in for $p\geq 1$ and $q,r\geq 0$}
\label{f:masik}
\end{figure}
The slight modification of the graph $\Delta _{p,q,r}$ when $p=0$ is shown
in Figure~\ref{f:specmasik}. The set of graphs $\Delta _{p,q,r}$ with
$p,q,r\geq 0$ will be denoted by $\frn$.

\begin{figure}[ht]
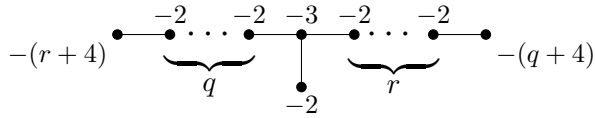

\begin{center}
\setlength{\unitlength}{1mm}
\unitlength=0.7cm
\begin{graph}(10,4)(0,-1)
\graphnodesize{0.2}

 \roundnode{m1}(1.5,0)
 \roundnode{m2}(2.5,0)
 \roundnode{m3}(4,0)
 \roundnode{m4}(5,0)
 \roundnode{m5}(6,0)
 \roundnode{m6}(7.5,0)
 \roundnode{m7}(8.5,0)
 \roundnode{m8}(5,-1)  

\edge{m1}{m2}
\edge{m3}{m4}
\edge{m4}{m5}
\edge{m6}{m7}
\edge{m4}{m8}

  \autonodetext{m1}[sw]{{\small $-(r+4)$}}
  \autonodetext{m2}[n]{{\small $-2$}}
  \autonodetext{m3}[n]{{\small $-2$}}
  \autonodetext{m4}[n]{{\small $-3$}}
  \autonodetext{m5}[n]{{\small $-2$}}
  \autonodetext{m6}[n]{{\small $-2$}}
  \autonodetext{m7}[se]{{\small $-(q+4)$}}
  \autonodetext{m8}[s]{{\small $-2$}}
  \autonodetext{m3}[w]{{\Large $\cdots$}}
  \autonodetext{m5}[e]{{\Large $\cdots$}}

\freetext(3.25,-0.8)
{$\underset{{\textstyle q}}{\underbrace{\hspace{1.2cm}}}$}

\freetext(6.75,-0.8)
{$\underset{{\textstyle r}}{\underbrace{\hspace{1.2cm}}}$}
\end{graph}
\end{center}
\caption{\quad The graph $\Delta _{0,q,r}$}
\label{f:specmasik}
\end{figure}

\item The plumbing graph of Figure~\ref{f:lambda} is $\Lambda _{p,q,r}$
with $p, r\geq 1$ and $q\geq 0$.
\begin{figure}[ht]
\begin{center}
\setlength{\unitlength}{1mm}
\unitlength=0.7cm
\begin{graph}(-3,2)(0,-1)
\graphnodesize{0.2}

 \roundnode{m15}(-5.5,0)
 \roundnode{m14}(-4.5,0)
 \roundnode{m13}(-3,0)
 \roundnode{m12}(-2,0)
 \roundnode{m11}(-1,0)
 \roundnode{m10}(0.5,0)
 \roundnode{m1}(1.5,0)
 \roundnode{m2}(2.5,0)
 \roundnode{m3}(4,0)
 \roundnode{m4}(5,0)
 \roundnode{m5}(6,0)
  \roundnode{m8}(5,-1)  

\edge{m1}{m2}
\edge{m3}{m4}
\edge{m4}{m5}
\edge{m4}{m8}
\edge{m12}{m11}
\edge{m10}{m1}
\edge{m15}{m14}
\edge{m13}{m12}

  \autonodetext{m1}[n]{{\small $-3$}}
  \autonodetext{m2}[n]{{\small $-2$}}
  \autonodetext{m3}[n]{{\small $-2$}}
  \autonodetext{m4}[n]{{\small $-2$}}
  \autonodetext{m5}[se]{{\small $-(p+2)$}}
  \autonodetext{m8}[s]{{\small $-(r+3)$}}
 \autonodetext{m10}[n]{{\small $-2$}}
 \autonodetext{m11}[n]{{\small $-2$}}
 \autonodetext{m12}[n]{{\small $-3$}} 
  \autonodetext{m13}[n]{{\small $-2$}}
 \autonodetext{m14}[n]{{\small $-2$}}  
 \autonodetext{m15}[sw]{{\small $-(q+4)$}} 
\autonodetext{m3}[w]{{\Large $\cdots$}}
  \autonodetext{m10}[w]{{\Large $\cdots$}}
\autonodetext{m13}[w]{{\Large $\cdots$}}

\freetext(3.25,-0.8)
{$\underset{{\textstyle q}}{\underbrace{\hspace{1.2cm}}}$}
\freetext(-0.25,-0.8)
{$\underset{{\textstyle r-1}}{\underbrace{\hspace{1.2cm}}}$}

\freetext(-3.75,-0.8)
{$\underset{{\textstyle p-1}}{\underbrace{\hspace{1.2cm}}}$}
\end{graph}
\end{center}
\caption{\quad The graph $\Lambda _{p,q,r}$ in for $p,r\geq 1$ and $q\geq 0$}
\label{f:lambda}
\end{figure}
Modifications of these graphs for $p=0, r\geq 1$ and $p\geq 1, r=0$
are shown by Figures~\ref{f:lambdaspec1} and \ref{f:lambdaspec2}. The
last degeneration, when $p=r=0$ (as shown by Figure~\ref{f:cs}) will
appear in the  family $\frc$ defined below. The set of graphs  
$\{ \Lambda _{p,q,r}\mid p,q,r\geq 0, \ (p,r)\neq (0,0)\}$ will be
denoted by $\farm$.

\begin{figure}[ht]
\begin{center}
\setlength{\unitlength}{1mm}
\unitlength=0.7cm
\begin{graph}(1,2)(0,-1)
\graphnodesize{0.2}

 \roundnode{m12}(-2,0)
 \roundnode{m11}(-1,0)
 \roundnode{m10}(0.5,0)
 \roundnode{m1}(1.5,0)
 \roundnode{m2}(2.5,0)
 \roundnode{m3}(4,0)
 \roundnode{m4}(5,0)
 \roundnode{m5}(6,0)
  \roundnode{m8}(5,-1)  

\edge{m1}{m2}
\edge{m3}{m4}
\edge{m4}{m5}
\edge{m4}{m8}
\edge{m12}{m11}
\edge{m10}{m1}

  \autonodetext{m1}[n]{{\small $-3$}}
  \autonodetext{m2}[n]{{\small $-2$}}
  \autonodetext{m3}[n]{{\small $-2$}}
  \autonodetext{m4}[n]{{\small $-2$}}
  \autonodetext{m5}[s]{{\small $-2$}}
  \autonodetext{m8}[s]{{\small $-(r+3)$}}
 \autonodetext{m10}[n]{{\small $-2$}}
 \autonodetext{m11}[n]{{\small $-2$}}
 \autonodetext{m12}[sw]{{\small $-(q+5)$}} 
\autonodetext{m3}[w]{{\Large $\cdots$}}
 \autonodetext{m10}[w]{{\Large $\cdots$}}

\freetext(3.25,-0.8)
{$\underset{{\textstyle q}}{\underbrace{\hspace{1.2cm}}}$}
\freetext(-0.25,-0.8)
{$\underset{{\textstyle r-1}}{\underbrace{\hspace{1.2cm}}}$}

\end{graph}
\end{center}
\caption{\quad The graph $\Lambda _{0,q,r}$ for $r\geq 1$ and $q\geq 0$}
\label{f:lambdaspec1}
\end{figure}

\begin{figure}[ht]
\begin{center}
\setlength{\unitlength}{1mm}
\unitlength=0.7cm
\begin{graph}(1,2)(0,-1)
\graphnodesize{0.2}

 \roundnode{m12}(-2,0)
 \roundnode{m11}(-1,0)
 \roundnode{m10}(0.5,0)
 \roundnode{m1}(1.5,0)
 \roundnode{m2}(2.5,0)
 \roundnode{m3}(4,0)
 \roundnode{m4}(5,0)
 \roundnode{m5}(6,0)
  \roundnode{m8}(5,-1)  

\edge{m1}{m2}
\edge{m3}{m4}
\edge{m4}{m5}
\edge{m4}{m8}
\edge{m12}{m11}
\edge{m10}{m1}

  \autonodetext{m1}[n]{{\small $-4$}}
  \autonodetext{m2}[n]{{\small $-2$}}
  \autonodetext{m3}[n]{{\small $-2$}}
  \autonodetext{m4}[n]{{\small $-2$}}
  \autonodetext{m5}[se]{{\small $-(p+2)$}}
  \autonodetext{m8}[s]{{\small $-3$}}
 \autonodetext{m10}[n]{{\small $-2$}}
 \autonodetext{m11}[n]{{\small $-2$}}
 \autonodetext{m12}[sw]{{\small $-(q+4)$}} 
\autonodetext{m3}[w]{{\Large $\cdots$}}
 \autonodetext{m10}[w]{{\Large $\cdots$}}

\freetext(3.25,-0.8)
{$\underset{{\textstyle q}}{\underbrace{\hspace{1.2cm}}}$}
\freetext(-0.25,-0.8)
{$\underset{{\textstyle p-1}}{\underbrace{\hspace{1.2cm}}}$}

\end{graph}
\end{center}
\caption{\quad The graph $\Lambda _{p,q,0}$ for $p\geq 1$ and $q\geq 0$}
\label{f:lambdaspec2}
\end{figure}

\begin{figure}[ht]
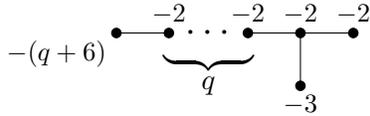

\begin{center}
\setlength{\unitlength}{1mm}
\unitlength=0.7cm
\begin{graph}(1,2)(3,-1)
\graphnodesize{0.2}

 \roundnode{m1}(1.5,0)
 \roundnode{m2}(2.5,0)
 \roundnode{m3}(4,0)
 \roundnode{m4}(5,0)
 \roundnode{m5}(6,0)
  \roundnode{m8}(5,-1)  

\edge{m1}{m2}
\edge{m3}{m4}
\edge{m4}{m5}
\edge{m4}{m8}

  \autonodetext{m1}[sw]{{\small $-(q+6)$}}
  \autonodetext{m2}[n]{{\small $-2$}}
  \autonodetext{m3}[n]{{\small $-2$}}
  \autonodetext{m4}[n]{{\small $-2$}}
  \autonodetext{m5}[n]{{\small $-2$}}
  \autonodetext{m8}[s]{{\small $-3$}}
  \autonodetext{m3}[w]{{\Large $\cdots$}}

\freetext(3.25,-0.8)
{$\underset{{\textstyle q}}{\underbrace{\hspace{1.2cm}}}$}

\end{graph}
\end{center}
\caption{\quad The graph $\Lambda _{0,q,0}$ for $q\geq 0$}
\label{f:cs}
\end{figure}

\item Let us define $\fra$ as the family of graphs we get in the
  following way: start with the graph of Figure~\ref{f:alap}(a), blow
  up its $(-1)$--vertex or any edge emanating from the $(-1)$--vertex
  and repeat this procedure of blowing up (either the new
  $(-1)$--vertex or an edge emanating from it) finitely many times,
  and finally modify the single $(-1)$--decoration to
  $(-4)$. Depending on the number and configuration of the chosen
  blow--ups, this procedure defines an infinite family of
  graphs. Define $\frb$ similarly, when starting with
  Figure~\ref{f:alap}(b) and substituting $(-1)$ in the last step with
  $(-3)$, and finally define $\frc$ in the same vein by starting with
  Figure~\ref{f:alap}(c) and putting $(-2)$ instead of $(-1)$ in the
  final step.
\begin{figure}[ht]
\begin{center}
\setlength{\unitlength}{1mm}
\unitlength=0.7cm
\begin{graph}(14,8.5)(0,-7)
\graphnodesize{0.2}

 \roundnode{m1}(0,0)
 \roundnode{m2}(2,0)
 \roundnode{m3}(4,0)
 \roundnode{m4}(2,-2)

 \roundnode{m5}(12,0)
 \roundnode{m6}(14,0)
 \roundnode{m7}(16,0)
 \roundnode{m8}(14,-2)

 \roundnode{m9}(6,-3)
 \roundnode{m10}(8,-3)
 \roundnode{m11}(10,-3)
 \roundnode{m12}(8,-5)

\edge{m1}{m2}
\edge{m2}{m3}
\edge{m2}{m4}

\edge{m5}{m6}
\edge{m6}{m7}
\edge{m6}{m8}

\edge{m9}{m10}
\edge{m10}{m11}
\edge{m10}{m12}

\freetext(2,-3){(a)}
\freetext(14,-3){(b)}
\freetext(8,-6){(c)}

  \autonodetext{m1}[n]{{\small $-3$}}
  \autonodetext{m2}[n]{{\small $-1$}}
  \autonodetext{m3}[n]{{\small $-3$}}
  \autonodetext{m4}[e]{{\small $-3$}}

  \autonodetext{m5}[n]{{\small $-4$}}
  \autonodetext{m6}[n]{{\small $-1$}}
  \autonodetext{m7}[n]{{\small $-4$}}
  \autonodetext{m8}[e]{{\small $-2$}}

  \autonodetext{m9}[n]{{\small $-2$}}
  \autonodetext{m10}[n]{{\small $-1$}}
  \autonodetext{m11}[n]{{\small $-6$}}
  \autonodetext{m12}[e]{{\small $-3$}}
\end{graph}
\end{center}
\caption{\quad Nonminimal plumbing trees giving rise to the families $\fra,
  \frb$ and $\frc$}
\label{f:alap}
\end{figure}
\end{itemize}
}
\end{defn}

\begin{rem}
{\rm 
Figure~\ref{f:blowup}
\begin{figure}[ht]
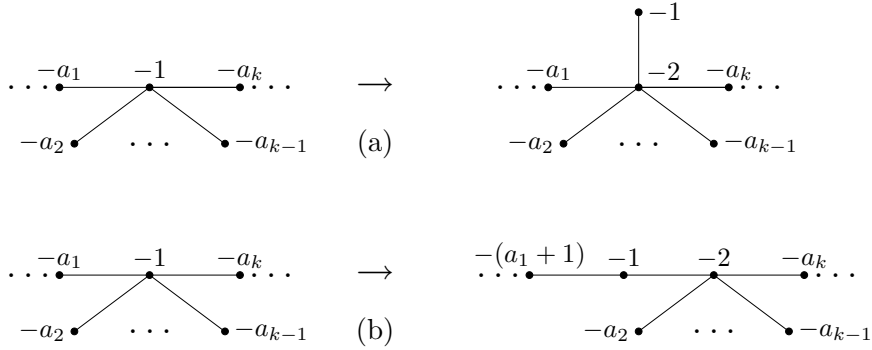

\begin{center}
\setlength{\unitlength}{1mm}
\unitlength=0.5cm
\begin{graph}(14,10)(0,-8)
\graphnodesize{0.2}

 \roundnode{m1}(-1,0)
 \roundnode{m2}(1.4,0)
 \roundnode{m3}(-3.4,0)
\roundnode{m4}(-3,-1.5)
\roundnode{m5}(1,-1.5)

\edge{m1}{m2}
\edge{m2}{m3}
\edge{m1}{m4}
\edge{m1}{m5}

\autonodetext{m1}[n]{{\small $-1$}}
\autonodetext{m2}[n]{{\small $-a_k$}}
\autonodetext{m3}[n]{{\small $-a_1$}}
\autonodetext{m3}[w]{{\Large $\ldots$}} 
\autonodetext{m2}[e]{{\Large $\ldots$}}

\autonodetext{m4}[w]{\small $-a_2$}
\autonodetext{m5}[e]{\small $-a_{k-1}$}

\freetext(5,0){\Large $\to$}
\freetext(5,-1.5){(a)}
\freetext(-1,-1.5){{\Large $\ldots$}}

\freetext(12,-1.5){{\Large $\ldots$}}

 \roundnode{n1}(12,0)
 \roundnode{n2}(14.4,0)
 \roundnode{n3}(9.6,0)
 \roundnode{n4}(10,-1.5)
 \roundnode{n5}(14,-1.5)

   \roundnode{n7}(12,2)

\edge{n1}{n2}
\edge{n2}{n3}
\edge{n1}{n7}
\edge{n1}{n4}
\edge{n1}{n5}

\autonodetext{n1}[ne]{{\small $-2$}}
\autonodetext{n2}[n]{{\small $-a_k$}}
\autonodetext{n3}[n]{{\small $-a_1$}}
\autonodetext{n7}[e]{{\small $-1$}}
\autonodetext{n2}[e]{{\Large $\ldots$}}
\autonodetext{n3}[w]{{\Large $\ldots$}}
\autonodetext{n4}[w]{\small $-a_2$}
\autonodetext{n5}[e]{\small $-a_{k-1}$}

 \roundnode{v1}(-3.4,-5)
 \roundnode{v2}(-1,-5)
 \roundnode{v3}(-3,-6.5)
 \roundnode{v4}(1,-6.5)
 \roundnode{v5}(1.4,-5)

\edge{v1}{v2}
\edge{v2}{v3}
\edge{v2}{v4}
\edge{v2}{v5}

\autonodetext{v1}[n]{{\small $-a_1$}}
\autonodetext{v2}[n]{{\small $-1$}}
\autonodetext{v1}[w]{{\Large $\ldots$}}
\autonodetext{v5}[e]{{\Large $\ldots$}}
\autonodetext{v3}[w]{\small $-a_2$}
\autonodetext{v4}[e]{\small $-a_{k-1}$}
\autonodetext{v5}[n]{\small $-a_{k}$}

\freetext(5,-5){\Large $\to $}
\freetext(-1,-6.5){{\Large $\ldots $}}

 \roundnode{u1}( 9.1,-5)
 \roundnode{u2}( 11.6,-5)
 \roundnode{u3}( 14,-5)
 \roundnode{u4}( 12,-6.5)
 \roundnode{u5}( 16,-6.5)
 \roundnode{u6}( 16.4,-5)

\edge{u1}{u2}
\edge{u2}{u3}
\edge{u3}{u4}
\edge{u3}{u5}
\edge{u3}{u6}

\autonodetext{u1}[n]{{\small $-(a_1+1)$}}
\autonodetext{u2}[n]{{\small $-1$}}
\autonodetext{u3}[n]{{\small $-2$}}
\autonodetext{u6}[e]{{\Large $\ldots$}}
\autonodetext{u1}[w]{{\Large $\ldots$}}
\autonodetext{u4}[w]{{\small $-a_2$}}
\autonodetext{u5}[e]{{\small $-a_{k-1}$}}
\autonodetext{u6}[n]{{\small $-a_k$}}

\freetext(14,-6.5){{\Large $\ldots$}}
\freetext(5,-6.5){(b)}

\end{graph}
\end{center}
\caption{\quad The blow--up of (a) a $(-1)$--vertex and (b) an edge
emanating from a $(-1)$--vertex}
\label{f:blowup}
\end{figure}
gives a pictorial description of what we mean by blowing up a $(-1)$--vertex
(Figure~\ref{f:blowup}(a)) and an edge emanating from a $(-1)$--vertex
(Figure~\ref{f:blowup}(b)). Notice that in the plumbing 4--manifold both
operations correspond to blowing up the $(-1)$--sphere defined by the vertex,
either in a generic point or in an intersection with another sphere of the
plumbing configuration. The reverse of the blow--up operation will be called
\emph{blowing down}; notice that only a vertex of valency 1 or 2 can be blown
down in the category of plumbing trees.  }
\end{rem}

\begin{rem}{\rm In classes $\frn$ and $\farm$, the need to make special graphs when some of $p,q,r$ are 0 disappears if we look instead at the "dual graphs", which arise in smoothings of negative weight (see Subsection~\ref{ss:negwt}).}
\end{rem}

After these preparations we are ready to state the main
theorems of the paper:

\begin{thm}\label{t:main}
The set $\frs $ is equal to the union 
$\frg \cup \frw \cup \frn \cup \farm \cup \fra \cup \frb \cup \frc$.
\end{thm}

We will show that these 7 types of graphs are in $\frs$ by exhibiting
explicit embeddings in $n\langle -1 \rangle$: for $\frg$, use
Proposition~\ref{p:g-s}; $\frw, \frn, \farm$ are done in the Appendix;
for $\fra, \frb, \frc$, proceed by induction as in
Proposition~\ref{p:c-s}.  (One can give an abstract proof for the
first four classes by combining the Examples in
Section~\ref{s:smoothing} with Corollary~\ref{Milnor}).  That $\frs$
consists only of these types is the heart of this paper.

The next question is which $\Gamma\in\frs$ are suitable for rational
blow--down, or which can be the graph of a normal surface singularity
admitting a rational homology disk smoothing. On the positive side,
in many cases there are several constructions for the appropriate
rational homology disks.  Using Kirby calculus, one can prove
(Section~\ref{s:kirby}):

\begin{thm}\label{t:kirby}
A plumbing graph $\Gamma \in \frg \cup \frw \cup \frn \cup
\farm$ defines the plumbing 4--manifold $\mg$ with the property that
$\partial \mg$ bounds a rational homology disk.
\end{thm}

In fact, a stronger result is true;  two methods for constructing smoothings of normal surface
singularities (Section~\ref{s:smoothing}) lead to:

\begin{thm}\label{t:smoothing}
Any graph $\Gamma \in \frg \cup \frw \cup \frn \cup \farm$ gives rise
to a normal surface singularity which admits a smoothing with
vanishing Milnor number. Some infinite families of graphs in each of
$\fra$, $\frb$, and $\frc$, including some with a node of valency 4,
admit the same property.
\end{thm}
(Concerning the valency 4 Examples~\ref{e:newfourv}, \ref{e:intex} and
\ref{e:valency4}, see \cite{kollar} for a related discussion.)  While
Theorem~\ref{t:main} excludes graphs with a node of valency $\geq 5$
or two nodes of valency $4$, it allows graphs with two (or more) nodes
of valency $\geq 3$; but at present none of these examples are known
to admit a smoothing with vanishing Milnor number.  Finally, it is
worth noting that any $\Gamma \in \frs$ which is star--shaped with a
node of valency 3 is \emph{taut} in the sense of Laufer
\cite{laufert}: there is a \emph{unique} singularity (necessarily
rational and weighted homogeneous) with graph $\Gamma$.  
On the other hand, \cite[Theorem 4.1]{laufert}
shows that any singularity with star--shaped $\Gamma\in\frs$ with
a node of valency 4 is weighted homogeneous, and its analytic type
is determined by the cross--ratio of the 4 special points on the
central curve in the minimal resolution.  But in
Example~\ref{e:intex}, for instance, there exists a homology disk
smoothing for \emph{only one} of these analytic types;
in fact, the deformation space of that particular singularity contains a
"smoothing component" of dimension 1, cf.
\cite[Theorem~3.13(c)]{Wahl}. 

On the other hand, certain classes in $\fra, \frb$ and $\frc$ give
$3$--manifolds which do not even bound rational homology disks (because
of their $\mu$--invariant).  For this and other obstructions, see
Subsection~\ref{ss:spec}.

The paper is organized as follows: In Section~\ref{s:motiv} we
motivate our definition of $\frs$ through a more detailed discussion
of the rational blow--down operation in the smooth and symplectic
category, and briefly review some basic facts about normal surface
singularities. In Section~\ref{s:plumbtree}, we consider a graph
$\Gamma\in\frs$, and write each vertex $v\in\Gamma$ as an integral
combination of the basis vectors $E_i$, deducing quickly strong
conditions on the coefficients, hence on the graphs.  The hard
technical work in
Sections~\ref{s:casea}, \ref{s:caseb} and \ref{s:casec} provides the proof of
Theorem~\ref{t:main} (with the addendum of describing the actual embeddings
of graphs in $\frw\cup \farm\cup \frn$, which is deferred to the Appendix
given in Section~\ref{s:append}). Section~\ref{s:kirby} is devoted to the
proof of the existence of rational homology disks for plumbing graphs in 
$\frg \cup \frw\cup \farm \cup \frn$; here we use techniques of smooth
4--manifold topology. Finally, in Section~\ref{s:smoothing} two 
algebro--geometric methods are described and applied to verify
the existence of rational homology disk smoothings for large families of 
singularities given by plumbing graphs of $\frs$.

{\bf {Acknowledgements:}} The first author would like to thank the
Institute for Advanced Study in Princeton for their hospitality while
some part of this work has been carried out.  AS was partially
supported by OTKA T49449; ZSz was supported by NSF grant number DMS
0406155; and JW was partially supported by NSA grant number
FA9550-06-1-0063.  The first and second authors were also supported by
EU Marie Curie TOK project BudAlgGeo. We would like to thank the
referee for his/her many valuable comments and suggestions.

\section{Plumbing trees}
\label{s:motiv}
The main purpose of this section is to justify
Definition~\ref{d:symptree}.  We show that if the plumbing 4--manifold
$\mg$ can be (Seiberg--Witten or symplectically) blown down then
$\Gamma \in \frs$; similarly, if $\Gamma$ is the resolution graph of a
normal surface singularity with a smoothing with vanishing Milnor
number, then $\Gamma \in \frs$.

\subsection{Tautly embedded plumbing trees}
\label{ss:taut}
Suppose that $\Gamma $ is a given negative--definite plumbing tree of
spheres and $\mg$ is the corresponding plumbed 4--manifold with
boundary. Let $\Sigma _v\subset \mg $ denote the 2--sphere
corresponding to the vertex $v\in \Gamma$. We say that $\mg$ is
\emph{tautly embedded} in a closed, oriented spin$^c$ 4--manifold $(X,
\s )$ if
\begin{equation}\label{e:taut}
\langle c_1( \s), [\Sigma _v]\rangle +[\Sigma _v ]^2=-2
\end{equation}
for all $v\in \Gamma$. With a slight abuse of terminology, sometimes we say
that the graph $\Gamma $ is tautly embedded. 

\begin{defn} \label{d:swblow}
{\rm The 4--manifold $\mg$ (or the plumbing graph $\Gamma$) can be
\emph{Seiberg--Witten blown down} if there is a 4--manifold $\bg$,
with the same rational homology as the disk $D^4$ and with $\partial
\bg \cong \partial \mg$, with the following property: for a taut
embedding $\mg\subset (X, \s) $ for which $SW_X(\s)\neq 0$,
 \begin{itemize}
 \item the spin$^c$ structure $\s \vert _{X-\mg}$
extends to $\s ' \in Spin ^c(X')$, where 
 $X'=(X-\mg )\cup \bg $, 
 \item  the Seiberg--Witten moduli spaces
${\mathcal {M}}_{X'}(\s ')$ and
${\mathcal {M}}_{X}(\s )$ have the same expected dimensions, and
\item $SW_{X}(\s ) =\pm  SW_{X'}(\s ')$.
\end{itemize}}
\end{defn}

Since the expected dimension of the Seiberg--Witten moduli space
is given by the formula 
\[
\dim {\mathcal {M}}_{X}(\s ) = \frac{1}{4}(c_1^2(\s ) -3 \sigma (X) -2e(X))
\]
(where, as usual,  $\sigma (X) $ and $e(X)$ denote the signature and
Euler characteristic of the 4--manifold $X$), the assumption 
$\dim {\mathcal {M}}_{X}(\s ) = \dim {\mathcal {M}}_{X'}(\s ')$
together with the facts
\[
\sigma (X')=\sigma (X) +\vert \Gamma \vert \mbox{  and  }
e(X')=e(X)-\vert \Gamma \vert
\]
(following from the negative--definiteness of $\Gamma$)
imply 
\[
c_1^2(\s )-c_1^2 (\s ')=-\vert \Gamma \vert .
\]
The symbol $\vert \Gamma \vert $ denotes the rank of $\Gamma$, that
is, the cardinality of its vertex set.

Let $\xg = \mg \cup {\overline {\bg}}$, where ${\overline {\bg}}$
denotes the 4--manifold $\bg $ with the opposite orientation.  Since
$\Gamma $ is a negative--definite graph and $\bg$ is a rational
homology disk, it follows that $\xg$ is a negative--definite closed,
oriented 4--manifold, and $\rank \ H_2(\xg ; \bfz ) = \vert \Gamma
\vert$.  According to Donaldson's famous diagonalizability theorem
\cite{D1}, \cite{Donor} we have that the intersection form $Q$ on $H_2 (\xg ;
\bfz )/Torsion$ is diagonalizable, providing an embedding of $\Gamma$
on $n$ vertices into the diagonal form $(\bfz ^n , Q_n )$ with $Q_n=n
\langle -1\rangle$.  (For a Heegaard Floer theoretic proof of
Donaldson's result, see \cite[Theorem~9.1]{Oszabs}.)  By gluing the
spin$^c$ structure $\s \vert _{\mg}$ to the extension of $\s \vert
_{\partial M_{\Gamma}}$ over ${\overline {\bg }}$ we get a spin$^c$
structure $\s _{\Gamma}\in Spin ^c (X _{\Gamma })$.  The assumption on
the expected dimensions of the Seiberg--Witten moduli spaces and the
fact that $H_2(\bg ; \bfq )=0$ readily implies that
\[
c_1^2(\s _{\Gamma})=c_1^2(\s ) - c_ 1^2 (\s ')=-\vert \Gamma \vert =-n
.
\]
Since in the negative--definite diagonal lattice of rank $n$ there is
essentially a unique characteristic element of square $-n $ (which is
$\sum _{i=1}^nE_i$ for an orthonormal basis $\{ E_1, \ldots , E_n
\}$), the assumption on $\mg$ being tautly embedded implies that
Equation~\eqref{e:adjunction} holds for every $v\in \Gamma$, providing
\begin{cor}\label{cswbd}
  If the minimal, negative--definite plumbing tree $\Gamma$ gives rise to $\mg$
  which can be Seiberg--Witten blown down, then
  $\Gamma \in \frs$. \hfill $\Box$
\end{cor}

\subsection{Trees which can be symplectically blown down}
Next we deal with a necessary condition for a negative--definite
plumbing tree $\Gamma$ to be blown down symplectically.  Suppose that
$\mg \subset (X, \omega )$ is a symplectic embedding into the closed
4--manifold $X$ with symplectic form $\omega$, that is, the spheres
$\Sigma _v$ are symplectic submanifolds in $(X, \omega )$ for all
$v\in \Gamma$ and they intersect each other
$\omega$--perpendicularly. Suppose that the 3--manifold $\partial
M_{\Gamma}$ bounds a rational homology disk $B_{\Gamma}$.  We say that
$\Gamma$ can be \emph{blown down symplectically} if $\omega \vert
_{X-M_{\Gamma}}$ can be extended over the glued--up rational homology
disk $B_{\Gamma}$ to produce a symplectic structure $\omega '$ on
$X'=(X-M_{\Gamma })\cup B_{\Gamma}$, cf.  \cite[Definition~1.1]{Sym0}.
If $\s _{\omega }\in Spin ^c (X)$ denotes the spin$^c$ structure
induced by $\omega$ then this property implies that
\begin{itemize}
\item $\mg $ is tautly embedded in $(X, \s _{\omega })$ (since the spheres $\Sigma _v$
are symplectic submanifolds),
\item $SW_X(\s _{\omega})=\pm 1$ and the expected dimension of the Seiberg--Witten
moduli space ${\mathcal {M}}_X(\s _{\omega} )$ is equal to zero, and finally
\item if $\omega '$ induces the spin$^c$ structure $\s '_{\omega
'}\in Spin ^c (X')$ then $SW_{X'}(\s ' _{\omega '}) = \pm 1$ and the
expected dimension of the Seiberg--Witten moduli space ${\mathcal
{M}}_{X'}(\s '_{\omega '})$ is zero.
\end{itemize}
Therefore $\Gamma $ fits in the category of the previous subsection,
hence if $\Gamma$ can be symplectically blown down then (at least for
the spin$^c$ structure induced by the symplectic structure) it can be
Seiberg--Witten blown down. This observation immediately yields

\begin{cor}\label{c2-3}
  If the minimal, negative--definite plumbing tree $\Gamma $ can be
  symplectically blown down, then $\Gamma \in \frs$. \hfill $\Box$
\end{cor}

\subsection{Normal surface singularities}
\label{ss:singularities}

For more details about the following, see e.g. \cite{LW}.  Let
$\rho\colon (\tilde{Y},E) \to (Y,o)$ be the minimal good resolution of
a germ of a normal surface singularity with rational homology sphere
link $\Sigma$ (i.e., $H_{1}(\Sigma;\mathbb Q)=0$).  Then the
exceptional divisor $E=\cup_{v}E_{v}=\rho ^{-1}(0)$ decomposes into
smooth rational curves, giving a dual resolution graph $\Gamma$, whose
vertices are indexed by the $v$'s; $\Gamma$ is a negative--definite
tree which (except for lens spaces) determines and is determined by
the homotopy type of the 3--manifold $\Sigma$.  Consider the lattice
$\mathbb E= \oplus _{v\in \Gamma}\mathbb Z[E_{v}] \subset
H_{2}(\tilde{Y}; \bfz )$ (called $L$ in Remark \ref{lattice}).  The
intersection pairing gives the discriminant group $D(\Gamma)\cong
\mathbb E^{*}/\mathbb E$ and
$$D(\Gamma)\cong H_{1}(\Sigma; \bfz ).$$ The complex structure on $\tilde{Y}$
gives rise to a relative cohomology class $K_{\tilde{Y}}$, and a rational
invariant $K_{\tilde{Y}}\cdot K_{\tilde{Y}}\in \mathbb Q$ computable from
$\Gamma$.  The analytic invariant $p_{g}(Y)= \text{dim}\ R^{1}\rho_{*}\mathcal
O_{\tilde{Y}}$ (the \emph{geometric genus}) is generally not computable from
$\Gamma$.  We say that $Y$ has a \emph{rational} singularity if $p_{g}(Y)=0$,
and this is determined by the graph condition $Z\cdot (Z+K_{\tilde{Y}})=-2,$
where $Z$ (the \emph{fundamental cycle}) is the smallest effective cycle for
which $Z\cdot E_{v}\leq 0$, all $v$.  For example, any negative--definite tree
$\Gamma$ for which the valency of every vertex is at most the negative of the
self--intersection is automatically rational.

Consider a smoothing $f\colon (\mathcal Y,o)\to (\Delta,o)$ of
$(Y,o)$, where $(\Delta,o)$ is the germ of an open disk in $\mathbb C$
(in particular,  $(f^{-1}(o),o)=(Y,o)$).  With appropriate attention to the boundaries
of representatives of these germs, one can define a ``general fiber''
of the smoothing, the Milnor fiber $M$.  $M$ is a compact oriented
4--manifold with boundary $\Sigma$, and its first Betti number is 0;
it also has a rational invariant $K_{M}\cdot K_{M}\in \mathbb Q$.  The
Milnor number $\mu$ of the smoothing is the rank of $H_{2}(M; \bfq)$,
and in general depends on which smoothing of $Y$ one considers.
Denoting the Sylvester invariants of the intersection pairing on
$H_{2}(M; \bfq)$ by $(\mu_{0},\mu_{+},\mu_{-})$, one has the general
formulae
$$\mu_{0}+\mu_{+}=2p_{g}(Y),$$
$$K_{M}\cdot K_{M} +\chi (M)=K_{\tilde{Y}}\cdot K_{\tilde{Y}}+\chi
(\tilde{Y})+12p_{g}(Y).$$ 
Now suppose that one has a singularity and a smoothing whose Milnor fibre is a rational homology disk
(i.e., $\mu=\mu _0+\mu _+ +\mu _-=0$); then for topological reasons $K_M=0$.  But the above formulas and discussion imply further that $Y$ has
a rational singularity, and hence (if $n$ denotes the number of exceptional
curves)
$$n+K_{\tilde{Y}}\cdot K_{\tilde{Y}}=0.$$

For cyclic quotient singularities $\mathbb C^{2}/G$, $\Gamma$ is a
string, and the condition above is valid only for the type
$p^{2}/(pq-1)$ (see \cite[Theorem (5.10)]{LW}).

The finite group $H_{1}(\Sigma ; \bfz )$ has a natural non--degenerate
linking pairing into $\mathbb Q/\mathbb Z$, which is the same as the
one on $D(\Gamma)$; but this is induced via a finer object, a
\emph{quadratic function} $q\colon H_{1}(\Sigma ; \bfz )\to \mathbb
Q/\mathbb Z,$ coming from the function on $\mathbb E$ given by $e
\mapsto (1/2)(e\cdot e +e\cdot K_{\tilde{Y}})$.  Theorem~4.5 of
\cite{LW} implies that if a Milnor fiber $M$ has $H_{2}(M;\mathbb
Q)=0$, then the kernel $I$ of the natural surjection $H_{1}(\Sigma ;
\bfz)\to H_{1}(M; \bfz )$ is $q$--isotropic; in particular,
$I=I^{\perp}$, so the order of $D(\Gamma)$ is a square, and the class
of $K_{\tilde Y}$ is in $I$.  In addition, by gluing the Milnor fiber
$M$ with $\mu =0$ (with opposite orientation) to $\tilde{Y}$ along
$\Sigma$ we get a negative--definite 4--manifold; hence, by
Donaldson's Theorem we conclude that $\Gamma$ embeds into the
negative--definite diagonal lattice of rank $\vert \Gamma \vert =
n$. We collect these results in

\begin{prop}\label{rat} Suppose that the  normal surface singularity $(Y,o)$ 
with resolution graph $\Gamma$ admits a smoothing whose Milnor fiber
  $M$ is a rational homology disk.  Then
    \begin{enumerate}
    \item $Y$ has a rational singularity,
    \item on the minimal good resolution
    $\vert \Gamma \vert +K_{\tilde{Y}}\cdot K_{\tilde{Y}}=0$,
    \item there is a self--isotropic subgroup $I\subset H_{1}(\Sigma ;
    \bfz )$, and $H_{1}(M; \bfz )\cong H_{1}(\Sigma ; \bfz )/I$
    (hence is$ \cong I)$ and
    \item $\Gamma $ admits an embedding into the diagonal lattice
     $(\bfz ^n , Q_n)$ of the same rank, and the class of $K_{\tilde{Y}}$ is a characteristic element of the lattice. \hfill $\Box$
    \end{enumerate}
    \end{prop}
In conclusion
\begin{cor}\label{Milnor}
If the minimal plumbing tree $\Gamma$ gives rise to a complex surface 
singularity which admits a rational homology disk smoothing, then
$\Gamma \in \frs$. \hfill $\Box$
\end{cor}

Since the only rational singularities with trivial $D(\Gamma)$ are a
smooth point and rational double point of type $E_{8}$, it follows
that every Milnor fiber $M$ with $\mu=0$ has non--trivial $H_{1}(M;
\bfz )$.  Examples with $\pi_{1}(M)$ finite and non--abelian will be
shown in Section~\ref{s:smoothing}, e.g. Example~\ref{e:intex}.  On
the other hand, recall that Milnor fibers of hypersurface and complete
intersection singularities are always simply connected.

\begin{rem}
{\rm As it was pointed out earlier, if $\Gamma$ can be blown down
symplectically, then it can be Seiberg--Witten blown down. Correspondingly,
if a singularity with resolution graph $\Gamma$ admits a rational 
homology disk smoothing, then $\Gamma$ can be blown down symplectically.
This last observation is a direct consequence of the main result of 
\cite{GSuj}.}
\end{rem}

\subsection{Further considerations}\label{ss:spec}
Notice that graphs in $\frs$ only capture some aspects of the
combinatorics for a configuration to admit a rational homology disk
smoothing or to be one which can be rationally blown down. We make
some short remarks concerning the question of which elements of $\frs$
actually arise in one of our three situations.  Further discussion is
postponed to a future project.

 If $\Sigma$ is a compact rational homology $3$--sphere with no
$2$--torsion in its first homology, then there is a well--known
topological obstruction for $\Sigma$ to bound (smoothly) a rational
homology disk.  The \emph{$\mu$-invariant} of $\Sigma$ is an integer
mod $16$ computed from the signature of a spin 4--manifold bounded by
$\Sigma$ (e.g., \cite[p.~46]{hirz}); it must be 0 if $\Sigma$ bounds a
rational homology disk.  Neumann and Raymond \cite{raym} show how to
compute $\mu$ for a plumbed 3--manifold.  One easily finds that the
class $\frc$ example given by Figure~\ref{f:classcex}
has $\mu = 8$, hence its boundary cannot bound a rational homology
disk.  Similar examples exist for types $\fra$ and $\frb$.

\begin{figure}[ht]
\begin{center}
\setlength{\unitlength}{1mm}
\unitlength=0.7cm
\begin{graph}(10,4)(0,-1)
\graphnodesize{0.2}

 \roundnode{m1}(1.5,0)
 \roundnode{m2}(2.5,0)
 \roundnode{m3}(3.5,0)
 \roundnode{m4}(4.5,0)
 \roundnode{m5}(5.5,0)
 \roundnode{m8}(4.5,-1)  

\edge{m1}{m2}
\edge{m3}{m2}
\edge{m4}{m3}
\edge{m5}{m4}
\edge{m4}{m8}

  \autonodetext{m1}[n]{{\small $-7$}}
  \autonodetext{m2}[n]{{\small $-2$}}
  \autonodetext{m3}[n]{{\small $-2$}}
  \autonodetext{m4}[n]{{\small $-3$}}
  \autonodetext{m5}[n]{{\small $-2$}}
  \autonodetext{m8}[s]{{\small $-3$}}
\end{graph}
\end{center}
\caption{\quad A graph in $\frc$ with $\mu =8$}
\label{f:classcex}
\end{figure}

There are further exclusions when considering smoothability of surface
singularities; as already mentioned, any such must have the graph of a
rational surface singularity.  For example the graph of
Figure~\ref{f:notgood} is of class $\frb$, but represents a minimally
elliptic singularity (of geometric genus 1 \cite{laufermin}), so could
not give a rational homology disk smoothing.  Futher obstruction for a
surface singularity to admit rational homology disk smoothing is
provided by the integer valued invariant ${\overline {\mu}}$ of
Neumann, cf. \cite{SPAMS}.

\begin{figure}[ht]
\begin{center}
\setlength{\unitlength}{1mm}
\unitlength=0.7cm
\begin{graph}(12,6)(3,-8)
\graphnodesize{0.2}

 \roundnode{m9}(6,-4)
 \roundnode{m10}(8,-4)
 \roundnode{m11}(10,-4)
 \roundnode{m12}(8,-6)
 \roundnode{m13}(8,-2)

\edge{m9}{m10}
\edge{m10}{m11}
\edge{m10}{m12}
\edge{m10}{m13}

  \autonodetext{m9}[n]{{\small $-2$}}
  \autonodetext{m10}[ne]{{\small $-2$}}
  \autonodetext{m11}[n]{{\small $-4$}}
  \autonodetext{m12}[e]{{\small $-3$}}
  \autonodetext{m13}[e]{{\small $-4$}}

\end{graph}
\end{center}
\caption{\quad Graph in $\frb\subset \frs$ with no rational homology disk 
smoothing}
\label{f:notgood}
\end{figure}

\section{Symplectic plumbing trees}
\label{s:plumbtree}
The next four sections are devoted to the proof of
Theorem~\ref{t:main}. The arguments in these sections are purely
combinatorial.  We start with some generalities regarding plumbing
trees in $\frs$. From now on, we will identify the vertex $v$ of
$\Gamma$ with its image $\varphi (v)$ (a linear combination of the
$E_i$'s) in the diagonal lattice $\langle E_1, \ldots , E_n \rangle =
n \langle -1 \rangle $. We will say that $E_i$ is in $v$ (or $E_i$
contributes to $v$) if $Q(E_i,v)\neq 0$. In this case the
\emph{multiplicity} of $E_i$ in $v$ is equal to $-Q(E_i,v)$.  With a
slight abuse of notation, the product $Q(v,w)$ sometimes will be
written as $v\cdot w$ or even as $vw$.  

\begin{lem}\label{l:forma}
Suppose that $\Gamma \in \frs$. For a vertex $v\in \Gamma $ we have
either
\[ 
v =E_{i_v} - \sum _{j\in J_v} E_j
\]
or
\[ 
v =-2E_{i_v} - \sum _{j\in J_v} E_j
\]
with $i_v$ not being in the index set $J_v\subset \{ 1, \ldots ,n\}$.
\end{lem}
\begin{proof}
Suppose that $v=\sum _k \alpha _k E_k$ with $\alpha _k \in \bfz
$. From the fact that
\begin{equation}\label{e:adj}
(\sum _{i=1}^n E_i)\cdot v +v \cdot v=-2,
\end{equation}
we conclude that 
\[
-\sum _k \alpha _k (\alpha _k +1)=-2,
\]
implying that either exactly one $\alpha _k$ is equal to 1 and all the
others are 0 or $-1$, or exactly one is equal to $-2$ and all the
others are 0 or $-1$. This observation clearly implies the result.
\end{proof}
It is easy to see that (in some sense) the converse of this statement also
holds: if the minimal plumbing tree $\Gamma$ is embedded into the diagonal
lattice $n\langle -1\rangle$ with $n = \vert \Gamma \vert$ and for all $v\in
\Gamma$ we have that $v$ (in the orthogonal basis of the diagonal lattice) has
one of the forms listed in the statement of the lemma above, then $v$
satisfies \eqref{e:adj} and so $\Gamma $ is, in fact, in $\frs$.

\begin{cor}\label{c:subtree}
Suppose that $\{ v_{i_1}, \ldots , v_{i_k}\}$ are vertices of
$\Gamma\in \frs$ determining a connected subtree. Then
\[
\sum _{j=1}^k v_{i_j}
\]
is of the form given by Lemma~\ref{l:forma}.
\end{cor}
\begin{proof}
Notice that if $v_1$ and $v_2$ both satisfy \eqref{e:adj} and
$v_1,v_2$ are adjacent vertices in $\Gamma$ (that is, $v_1 \cdot v_2
=1$) then
\[
(\sum _{i=1}^n E_i)\cdot  (v_1+v_2) +(v_1+v_2)^2=-2,
\] 
hence the sum of the classes $v_1$ and $v_2$ also satisfies
\eqref{e:adj}.  Then Lemma~\ref{l:forma} together with a simple
inductive argument implies the corollary.
\end{proof}
Applying the above idea for $\Gamma $ itself we get the following
\begin{cor}\label{c:sum}
If $\Gamma \in \frs$ is defined on $n$ vertices $\{ v_1, \ldots ,
v_n\}$ then
\[
\vert \sum _{i=1}^n v_i^2 \vert \leq 3n+1.
\]
\end{cor}
\begin{proof}
  Notice that $(\sum _{i=1}^n v_i )^2=2\sum _{i<j} v_iv_j+\sum _{i=1}^n
  v_i^2$, and since $\Gamma $ is a tree, we have $\sum _{i<j}v_iv_j=n-1$. Now
  the fact $\sum _{i=1}^n v_i=-2E_i-\sum _{j\in J } E_j$ or $\sum _{i=1}^n
  v_i=E_i-\sum _{j\in J } E_j$ (with $J \subset \{ 1, \ldots , n \} -\{ i\}$)
  implies that $(\sum _{i=1}^n v_i)^2 $ is either $-4-\vert J \vert \geq -3-n$
  or $-1-\vert J \vert\geq -n$. This observation implies the statement.
\end{proof}

\begin{prop}\label{p:atmost1}
Suppose that for $\Gamma$ as above there is a vertex $v$ of the form
$-2E_{i_v} -\sum _{j\in J_v} E_j$. Then for all other vertices $w\in
\Gamma $ we have
\[
w=E_{i_w}-\sum _{j\in J_w}E_j.
\]
\end{prop}
\begin{proof}
Suppose that there are two vertices $v_1, v_2$ 
in $\Gamma$ such that 
\[
v_k=-2E_{i_{v_k}}-\sum _{j\in J_{v_k}}E_j \qquad (k=1,2).
\]
Consider the path $\{ w_l \mid l=1, \ldots , m\}$ of vertices
connecting $v_1$ and $v_2$ in $\Gamma$, that is, $w_1=v_1$, $w_n=v_2$
and $w_i$ is adjacent to $w_{i+1}$ for $i=1,\ldots , m-1$. Assume
furthermore that the elements $w_l$ for $l$ different from $1$ and $m$
are of the form $E_{i_{w_l}}-\sum _{j\in J_{w_l}}E_j$.  (Otherwise we
could choose a shorter chain with endpoints admitting the same
properties as $v_1,v_2$ and satisfying this additional requirement.)
It is easy to see that if $m>2$ then $w_1+w_2$ (and so by induction
$w_1 +\ldots +w_{m-1}$ is of the form $-2E_i-\sum _{j\in J} E_j$ for
some $i$ and index set $J$): if $w_1=-2E_{i_{w_1}}-\sum _{j\in
J_{w_1}}E_j$ and $w_2=E_{i_{w_2}}-\sum _{j\in J_{w_2}}E_j$ then
$Q(w_1,w_2)=1$ implies that either $i_{w_2}\in J_{w_1}$ and
$J_{w_1}\cap J_{w_2}=\emptyset$ (giving $w_1+w_2=-2E_{i_{w_1}}-\sum
_{j \in J_{w_1}\cup J_{w_2}-{i_{w_2}}}$), or $i_{w_1}=i_{w_2}$ and
$J_{w_1}\cap J_{w_2}=\{ p\}$. In this latter case $w_1+w_2=-2E_p -\sum
_{j \in J_{w_1}\cup J_{w_2}\cup \{ i_{w_1} \} -\{ p\}}E_j$.  Then the
sum $w_1+\ldots +w_{m}$ would violate \eqref{e:adj}, since it contains
at least two $E_k$'s with multiplicity $-2$, contradicting
Corollary~\ref{c:subtree}.
\end{proof}

In short, by Proposition~\ref{p:atmost1} in $\Gamma\in \frs$ there is at most
one vertex $v$ which is of the form $-2E_{i_v}-\sum _{j\in J_v} E_j$ and all
other vertices are of the shape $E_k-\sum _{j\in J} E_j$.  A slight
generalization of the above proposition follows from the same principle:
\begin{cor}\label{c:alttree}
Suppose that $\Gamma \in \frs $  and $\Gamma _1, \Gamma _2 \subset \Gamma$
are disjoint connected subtrees. Then at least one of the vectors
$w_k=\sum _{v\in \Gamma _k}v$ ($k=1,2$) is of the form
$E_{i_{w_k}}-\sum _{j\in J_{w_k}}E_j$. \hfill $\Box$
\end{cor}

\begin{lem}\label{l:dist}
Suppose that for two vertices $v_1\neq v_2$  we have
\[
v_k=E_{i_{v_k}}-\sum _{j\in J_{v_k}}E_j \qquad (k=1,2).
\]
Then $i_{v_1}\neq i_{v_2}$.
\end{lem}
\begin{proof}
Recall that the pairing $Q (v_1,v_2)$ is either 1 or 0 (depending on whether
$v_1$ and $v_2$ are adjacent in $\Gamma$ or not). If $i_{v_1}=i_{v_2}$ then
\[
)\leq Q (v_1,v_2)=(E_{i_{v_1}}-\sum _{j\in J_{v_1}}E_j)\cdot
(E_{i_{v_2}}-\sum _{j\in J_{v_2}}E_j)=-1+\sum _{j\in J_{v_1}\cap
J_{v_2}}(-1)\leq -1,
\]
providing the desired contradiction.
\end{proof}
  
Consequently, each $E_i$ appears in at most one $v_j\in \Gamma$ with
positive coefficient.  By Corollary~\ref{c:subtree}, the sum
$w=\sum _{v\in \Gamma} v$ is either of the form $E_i-\sum _{j\in
J}E_j$ or of the form $-2E_k-\sum _{j\in J}E_j$; so, summing all the coefficients of $E_i$ gives $1$, $0$, $-1$, or $-2$.  Hence, for an index $i$ there are
ten possibilities for $E_i$ to be contained by vectors of
$\Gamma$:

\begin{lem}\label{l:esetek}
If $i$ is an index between $1$ and $n$ then for the basis vector $E_i$
of the diagonal lattice $Q_n$ one of the following ten possibilities
can occur:
\end{lem}
\begin{enumerate}
{\it \item $E_i$ appears in a single vector of $\Gamma$ with multiplicity $1$;
\item $E_i$ appears in a single vector of $\Gamma$ with multiplicity $-1$;
\item $E_i$ appears in a single vector of $\Gamma$ with multiplicity $-2$;
\item $E_i$ appears in two vectors of $\Gamma$, once with multiplicity $1$ and in another
with multiplicity $-1$;
\item $E_i$ is in two vectors of $\Gamma$, with multiplicities $1$ and $-2$;
\item $E_i$ is in two vectors of $\Gamma$, with both multiplicities $-1$;
\item $E_i$ is in three vectors of $\Gamma$, with multiplicities $1$, $-1$, $-1$;
\item $E_i$ is in three vectors of $\Gamma$, with multiplicities $1$, $-1$ and $-2$;
\item $E_i$ is in four vectors of $\Gamma$, once with multiplicity 1, and in three
with multiplicity $-1$;
\item $E_i$ does not appear in any vector of $\Gamma$ at all.} \hfill $\Box$
\end{enumerate}

We start with a few observations regarding the types of indices
occuring in an embedding $\Gamma \subset n \langle -1 \rangle$.
First, it is easy to see that (8) cannot occur, since if $v_1$
contains $-E_i$ and $v_2$ contains $-2E_i$ then $Q(v_1, v_2)\leq -1$,
contradicting our assumptions.  If $i$ is an index of type (10), then
adding $-E_i$ to a single (but otherwise arbitrary) $v\in \Gamma$ we
change $\Gamma$ and $\varphi$ (through the decoration of the vertex to
which we added $-E_i$) in a way that the resulting $\Gamma '$ is still
in $\frs$ but the index of type (10) is turned into an index of type
(2).  (Later we will see that if $i$ is an index of type (2) then
$-E_i$ can appear only in a vector of square $-2$, hence we cannot
change an index of type (2) to one of type (10) within the category of
minimal plumbing trees.)  In the following we will always assume that
$\Gamma $ admits no index of type (10).

If $\Gamma \subset n \langle -1 \rangle$ involves an index of type (1)
then by changing both $\Gamma$ and the embedding $\varphi$ we get a
new embedding $\varphi '\colon \Gamma ' \to n \langle -1 \rangle $
where no index of type (1) exists: if $v=E_{i_v}-\sum _{j \in J_v}
E_j$ and $E_{i_v}$ does not appear in any other vector (that is, $i_v$
is of type (1)) then by changing $v$ to $v'=-2E_{i_v}-\sum _{j \in
J_v}E_j$ we get another $\Gamma '\in \frs$ which only differs from
$\Gamma$ by the decoration of $v$ (which has been reduced by 3).
Therefore a type (1) index can be changed to be of type (3).
Therefore we can assume that $\Gamma$ contains no index of type (1).
The next theorem shows that, in fact, a type (3) index can never
appear in a graph of $\frs$ defined on more than one vertex:

\begin{thm}\label{t:notype3}
There is no minimal symplectic plumbing tree $\Gamma \in \frs$ on 
$n>1$ vertices involving an index of type (3).
\end{thm}
\begin{proof}
A simple check shows that the theorem holds for $n=2$. Suppose that
$\Gamma$ is a plumbing tree on $n$ vertices ($n\geq 3$) and it
involves an embedding into $n\langle -1 \rangle $ with an index $i$ of
type (3). Consider such a $\Gamma$ with the smallest possible $n$. Let
$v=E_{i_v}-\sum _{j \in J_v}E_j$ be a \emph{leaf} of the graph $\Gamma
$ (i.e., a vertex of valency one). The existence of such $v$ follows
from our assumption $n>1$. It follows from its form that $E_i$
does not appear in $v$.  If $i_v$ is of type (5), (7) or (9) then the
sum
\[
u=\sum _{v_j\neq v}v_j
\]
of all the other vertices (which form a connected subgraph, since $v$
is a leaf) will contain $-E_{i_v}$ with multiplicity at least 2, but
$u$ also contains $-2E_i$ (recall that $i$ is the index of type (3))
contradicting Corollary~\ref{c:subtree}.  If $i_v$ is of type (4) (the
only remaining possibility), then consider the unique vertex $v'\in
\Gamma$ with $Q(v,v')=1$, and distinguish two cases according to
whether $-E_{i_v}$ does or does not appear in $v'$.  Suppose first
that the unique $-E_{i_v}$ is in $v'$. By changing $\Gamma $ to
$\Gamma '$ via replacing the edge ${\overline {vv'}}$ (and its two
vertices) by one vertex $v+v'$, we eliminate the index $i_v$ from the
embedding; but $\Gamma '\in \frs$ is now a minimal symplectic plumbing
tree on $(n-1)$ vertices (with an embedding into $(n-1) \langle -1
\rangle$) which still has $i$ as an index of type (3), contradicting
our choice of $\Gamma$.  (Notice that since $(v+v')^2=v^2
+(v')^2+2\leq v^2<-1$, minimality of $\Gamma '$ obviously follows.)
Finally, if $-E_{i_v}$ appears in a vertex $v''$ different from $v'$
then change $\Gamma $ by deleting $v$ and replacing the vector $v''$
by $v''+E_{i_v}$. Once again, the resulting $\Gamma ''$ will be
contained by $(n-1)\langle -1 \rangle $ (since the index $i_v$ has
been eliminated) and $i$ is still an index of type (3) in $\Gamma ''$,
hence if it is minimal, then it contradicts our choice for
$\Gamma$. The plumbing graph $\Gamma ''$ can contain only a unique
vector with square $-1$: this is $v''+E_{i_v}$, implying that
$v''=E_{i_{v''}}-E_{i_v}$.  In this case, however
\[
0=Q(v,v'')=1+Q(E_{i_v}-\sum _{j\in J_v}E_j, E_{i_{v''}})\geq 1 ,
\]
which is a contradiction.
\end{proof}
\begin{rem}\label{r:ifandonlyif} {\rm Notice that for $n=1$ the plumbing graph
    $(-4)$ on one vertex (with the embedding $-2E_1$) does admit an
    index of type (3), hence our assumption $n>1$ is essential. In
    fact, the above theorem can be rephrased as follows: $\Gamma \in
    \frs $ contains an index of type (3) if and only if $\vert \Gamma
    \vert =1$. This form of the previous result will be very helpful
    in our later arguments, cf. Lemma~\ref{l:red3caseb}.}
\end{rem}
\begin{rem}
{\rm We also note that the line of argument above
(considering a leaf $v=E_{i_v}-\sum _{j\in J_v}E_j$, distinguishing two cases
according to whether $i_v$ is of types (5),(7), (9), or of type (4), and in 
this latter case examining where $-E_{i_v}$ is) will be repeatedly applied
throughout the proof of Theorem~\ref{t:main}.}
\end{rem}

From the result above we conclude that an index of type (1) in $\Gamma
\in \frs$ is also impossible once $\vert \Gamma \vert >1$. (For a
graph $\Gamma $ with $\vert \Gamma \vert =1$ the presence of an index
of type (1) contradicts the minimality of $\Gamma$.)  It follows then
that the sum $w=\sum _{v\in \Gamma } v$ is of the form
\[
-2E_{w}-\sum _{j\in J_w}E_j .
\]

In order to analyze more systematically which types of indices can
actually occur, let $\{ v_1, \ldots , v_n \}$ denote the vertices of
$\Gamma$ and notice that the sum
\begin{equation}\label{e:osszeg}
\sum _{i<j} Q(v_i, v_j)
\end{equation}
is equal to the number of edges in $\Gamma$, and since $\Gamma$ is
assumed to be a connected tree, it is equal to $n-1$. The sum
\eqref{e:osszeg} can be expanded according to the expansions of the
vectors $v_i$ in the basis $\{ E_1, \ldots , E_n\}$, and since $Q(E_i,
E_j)=0$ once $i\neq j$, the sum $n-1= \sum _{i<j} Q(v_i, v_j)$
decomposes as the sum of contributions of individual indices through
the multiplicity of $Q(E_i, E_i)$ appearing in the expansion of
\eqref{e:osszeg}.  It is easy to see that for indices of types (4) and
(7) this contribution is 1, for type (5) indices it is 2, if $i$ is of
type (6) then this contribution is $-1$ and finally for type (2) and
(9) indices it is 0.

Now we can have a better picture about the possible types of indices
$\Gamma$ can have. Since $w=\sum _{i=1}^n v_i$ is of the form
$-2E_i-\sum _{j\in J_i}E_j$, there is exactly one index $i$ of type
(6) or (9).
A simple analysis shows that 
\begin{itemize}
\item {\bf {(A)}} If this index $i$ is of type (6), then the remaining
$n-1$ indices contribute $n$ to the sum \eqref{e:osszeg}, therefore
one of them must be of type (5). By Proposition~\ref{p:atmost1} at
most one index of type (5) can exist, so the rest of the indices are
either of type (4) or of type (7).

 \item {\bf {(B)}} If $i$ is of type (9) and there is an index of type
(5) among the remaining ones (and then it is necessarily unique),
then there  exists  a unique index of type (2), and the
rest of the indices are of types (4) or (7).
  
\item {\bf {(C)}} If $i$ is of type (9) and there is no index of type
(5)  then all others are of types (4) or (7).
\end{itemize} 

In the next three sections we will analyize these three possibilities
separately.  More precisely, we will verify the following three
statements, the combination of which (together with
Propositions~\ref{p:g-s}, \ref{p:c-s} and \ref{p:sok-s}) prove
Theorem~\ref{t:main}.
\begin{thm}\label{t:casea}
If $\Gamma \in \frs$ with rank $\vert \Gamma \vert =n$ embeds into
$n\langle -1\rangle$ with indices as described in Case {\bf {(A)}} above
then $\Gamma \in \frg $.
\end{thm}
\begin{thm}\label{t:caseb}
If $\Gamma \in \frs$ with rank $\vert \Gamma \vert =n$ embeds into
$n\langle -1\rangle$ with indices as described in Case {\bf {(B)}} above
then $\Gamma \in \frc$.
\end{thm}

\begin{thm}\label{t:casec}
If $\Gamma \in \frs$ with rank $\vert \Gamma \vert =n$ embeds into
$n\langle -1\rangle$ with indices as described in Case {\bf {(C)}} above
then $\Gamma \in \frw\cup \frn \cup \farm \cup \fra\cup \frb$. 
\end{thm}

Before starting the proofs, however, we make some preparatory definitions and
observations.

\begin{defn}
{\rm A vector $v$ is \emph{full} if it is of the form $E_{i_v}-\sum
_{j\in J_v } E_j$ and $i_v$ is of type (5), (7) or (9). The plumbing
graph $\Gamma$ is \emph{full} if it has an embedding with no index of
type (4).}
\end{defn}

\begin{lem}
A symplectic plumbing tree $\Gamma \in \frs$ with $n$ vertices is full
if and only if (after possibly reordering the indices)
\[
\sum _{i=1}^n v_i=-E_1-E_2-\ldots - E_{n-1}-2E_n. 
\] 
\end{lem}
\begin{proof}
  By the exclusion of indices of types (1), (3), (4), (8) and
  (10), in the sum $w=\sum _{v\in \Gamma }v$ each $E_i$ comes with
  multiplicity $-1$ or $-2$, and there is a single one with
  multiplicity $-2$, verifying the lemma.
\end{proof}
\begin{lem}
A symplectic plumbing tree $\Gamma \in \frs$ is full if and only if
\[
\sum _{i=1}^n v_i^2=-3n-1.
\]
\end{lem}
\begin{proof}
Consider the vector $w=\sum _{i=1}^n v_i$. It is easy to 
see that $w^2=\sum _{i=1}^n v_i^2+2\sum _{i<j}v_iv_j$.
Since $\sum _{i<j}v_iv_j$ counts the number of edges in
$\Gamma$, and it is a connected tree, we get that 
\[
w^2=\sum _{i=1}^n v_i^2+2(n-1).
\]
Since by the previous observation $\Gamma$ is full if and only
if $w^2=-n-3$, the result follows.
\end{proof}
\begin{rem}
{\rm Notice that according to Corollary~\ref{c:sum} the graph is full
 if and only if $\vert \sum _{i=1}^n v_i^2\vert $ is maximal.  }
\end{rem}

\begin{prop}
The plumbing trees  $\Gamma _{p,q,r}\in \frw $
are full for all $p,q,r\geq 0$. 
\end{prop}
\begin{proof}
A  simple direct check provides the
proof: 
\[
\sum v_i^2=-2(p+q+r)-p-3-q-3-r-3-4=-3(p+q+r+4)-1.
\]
\end{proof}
Similarly, a somewhat more tedious, but straightforward computation shows that
\begin{prop}
The plumbing trees $\Delta_{p,q,r}, \Lambda _{p,q,r}\in \farm \cup
\frn $ are full for all $p,q,r\geq 0$. \hfill $\Box$
\end{prop}

It is not hard to see that the graphs of Figure~\ref{f:alap} become
full after adding $-1,-2,-3$ appropriately to the $(-1)$--vertex.
However, the graphs in the families $\fra, \frb, \frc$ are full only
in the cases when during the blow--up process we only blow up edges of
the configuration.  

\section{The classification of graphs of Case {\bf {(A)}}}
\label{s:casea}
Before starting the proof of Theorem~\ref{t:casea} we need a little
preparation, showing in particular that the graphs in $\frg $ are all
full.

\subsection{Continued fraction computations}
\label{ss:contfrac}
In this subsection we flesh out Remark (2.8.2) of \cite{wahl1}.
Define the class of plumbing trees $\frg _r$ as the minimal set of
plumbing trees which {\bf {(a)}} contains $\Gamma = (-4)$ and {\bf
{(b)}} if the linear chain $\Gamma =(-a_1, \ldots , -a_k) $ is in
$\frg _r$ ($a_i\in \bfn, a_i \geq 2$) then so are $\Gamma _1 = (-2,
-a_1, \ldots , -a_n-1)$ and $\Gamma _2= (-a_1-1, -a_2, \ldots , -a_n,
-2)$.  Notice that for all $\Gamma \in \frg _r$ we have that the tree
is, in fact, a chain. Let $[a_1, \ldots , a_k]$ (with $a_i\geq 2$
integers) denote the value of the continued fraction expansion
\[
a_1-\cfrac1{a_2-\cfrac1{\ddots - \cfrac1{a_{k}}}}
\] 
\begin{prop}
The plumbing tree $\Gamma = (-a_1, \ldots , -a_k)$ 
(with $a_i\geq 2$ integers) is in $\frg _r$
if and only if there are relatively prime  $p>q>0$ such that 
\[
[a_1, \ldots , a_k]=\frac{p^2}{pq-1}.
\]
\end{prop}
\begin{proof}
We start the proof with a simple auxiliary result concerning continued
fractions.  Suppose therefore that $[a_1, \ldots ,
a_k]=\frac{p^2}{pq-1}$; we  determine the value of $[a_1 +1,
\ldots , a_k, 2]$.  In the computation we will use the fact
\cite[Appendix]{OW} that if $[b_1, \ldots , b_n]=\frac{u}{v}$ then
$[b_n , \ldots , b_1]=\frac{u}{v'}$, where $vv'\equiv 1 \ \ (\bmod
\ u)$.  According to this principle and the definition, we get that
\[
[a_1+1, \ldots , a_k] = \frac{p^2}{pq-1}+1=\frac{p^2+pq-1}{pq-1},
\]
hence
\[
[a_k,\ldots ,a_1+1]=\frac{p^2+pq-1}{p^2-q^2-2},
\]
consequently
\[
[2,a_k,\ldots ,a_1 +1]=\frac{(p+q)^2}{p^2+pq-1},
\]
implying that 
\[
[a_1+1, \ldots , a_k, 2]=\frac{(p+q)^2}{q(p+q)-1}.
\]
In a similar manner we get that 
\[
[2,a_1, \ldots , a_k +1]=\frac{(p+q)^2}{p(p+q)-1}.
\]

Consider now a tree $\Gamma \in \frg _r$, given by the array $(-a_1,
\ldots , -a_k)$.  Using induction on $k$ it is fairly easy to see that
the statement of the Proposition holds: for $\Gamma = (-4)$ choose
$p=2$ and $q=1$; otherwise $\Gamma$ is constructed from $\Gamma '$ of
shorter length, hence induction, the above computation and the fact
that $p+q$ and $q$ are relatively prime if and only if $p$ and $q$
are, conclude one direction of the proof.

Suppose now that $\Gamma$ is given by an array $(-a_1, \ldots , -a_k)$
such that $[a_1, \ldots , a_k ]=\frac{p^2}{pq-1}$ for some relatively
prime integers $p$ and $q$. Induction on the length of the Euclidean
algorithm for determining $(p,q)=1$, together with the above formulae
and the uniqueness of the continued fraction expansion readily imply
the converse direction.
\end{proof}
Notice that, as a consequence of the above result, elements of $\frg$
are exactly the plumbing chains considered in \cite{Pratb} as
\emph{generalized rational blow--downs}. In particular, $\frg _r= \frg
$.

\begin{prop}
For any $p>q>0$ relative prime, the symplectic plumbing tree
$\Gamma _{p/q}$ is full.
\end{prop}
\begin{proof}
Since for $\vert \Gamma \vert =1$ the graph $(-4)$ is full, and the
inductive step constructing elements of $\frg _r$ decreases $\sum v_i^2$
by 3 when increasing the number of vertices by 1, it is obvious that
elements of $\frg _r$ are full. Now the above identification of $\frg$
with $\frg _r$ concludes the proof.
\end{proof}

\subsection{The proof of Theorem~\ref{t:casea}}

After these preliminary results we are ready to prove
Theorem~\ref{t:casea}.  It turns out to be more convenient to prove a
slightly stronger statement than Theorem~\ref{t:casea}, which
obviously implies the result.

\begin{prop}
Suppose that $\Gamma \in \frs$ involves an index $t$ of type (6),
that is, $\Gamma $ is listed in Case {\bf {(A)}}.  Then $\Gamma \in
\frg$ with the additional property that the two vectors containing
$-E_t$ are the two endpoints of the plumbing chain.
\end{prop}
\begin{proof}
The proof will proceed by induction on $n=\vert \Gamma \vert$.  The
statement is easy to check for $n=2$. Consider a graph $\Gamma$ on $n$
vertices; by our inductive hypothesis we can assume that for all $k<n$
the symplectic plumbing trees of Case {\bf {(A)}} defined on $k$
vertices satisfy the conclusion of the proposition. We can also assume
that $n\geq 3$.

First we show that if the two vectors $v_1,v_2$ containing $-E_t$ are
leaves of $\Gamma $ then $\Gamma \in \frg$. Add $E_t$ to $v_1,v_2$ and
delete the one (call it $v$) which has smaller decoration in absolute
value.  The resulting plumbing tree $\Gamma '\in \frs$, defined on
$(n-1)$ vertices, is still connected (since we deleted a leaf), and it
is embedded in $(n-1)\langle -1\rangle =\oplus _{i\neq t}\langle E_i
\rangle$. It follows from the above construction that $\Gamma '$ still
has no index of type (9), therefore it is still of Case {\bf {(A)}}.
We claim that $\Gamma '$ is also minimal. If it contains a vertex of
the form $E_{i}$, then the two leaves $v_1, v_2$ of $\Gamma$
containing $-E_t$ had the shape $v_1=E_{v_1}-E_t$ and
$v_2=E_{v_2}-E_t$.  This implies $Q(v_1, v_2)=-1$, a contradiction.
Therefore $\Gamma '$ is minimal, and so by our inductive hypothesis we
have that $\Gamma ' \in \frg$ and the index $t'$ of type (6) appears
in the two leaves of $\Gamma '$.  Since a graph in $\frg$ is full, for
the deleted vector $v$ we have $\vert v^2 \vert \leq 2$, therefore it
can only have the shape $v=E_{t'}-E_t$. This observation readily
implies that $\Gamma$ can be built from $\Gamma '$ by adding $-E_t$ to
the vector on one of its ends, and concatenating $\Gamma '$ with a new
vector $v$ having $v^2=-2$ on the other end.  Since $\Gamma '\in
\frg$, this shows that $\Gamma \in \frg$.  We also see that in
$\Gamma$ the vectors involving type (6) indices are on the two ends of
the plumbing chain.

Suppose now that (at least) one of the vectors $v_1,v_2$ containing
$-E_t$ is not a leaf, say it is $v_1$.  Our goal is to derive a
contradiction from this assumption.  Since $\vert \Gamma \vert \geq
3$, it has at least two leaves, therefore we have a leaf $v$ with
$Q(v, E_t)=0$.  The vector $v$ cannot be full, since then $\sum
_{u\neq v}u$ would contain both $E_{i_v}$ and $E_t$ with multiplicity
$-2$. Similarly, by Corollary~\ref{c:alttree} the vector $v$ cannot be
of the form $-2E_{i_v}-\sum _{j\in J_v}E_j$ either.  So we can assume
that $v=E_{i_v}-\sum _{j\in J_v} E_j$ and $i_v$ is of type (4), that
is, there is a single further $-E_{i_v}$ in one of the vectors of
$\Gamma$.  Let $v'$ denote the single vector adjacent to $v$ in
$\Gamma$.  If $-E_{i_v}$ is not in $v'$, then we can erase $v$ from
$\Gamma$ and $-E_{i_v}$ from the vector $w$ containing it, and the
resulting graph $\Gamma '$ (now on $(n-1)$ vertices, with an embedding
into $(n-1)\langle -1 \rangle=\oplus _{i\neq i_v}\langle E_i \rangle$)
is still minimal. Indeed, we changed only $w$, hence $\Gamma '$ is
nonminimal only if $w$ was of the form $E_{i_w}-E_{i_v}$, in which
case $Q(v,w)\geq 1$, giving a contradiction, since $w\neq v'$ and $v'$
is the only vector having nontivial pairing with $v$.  Noting that the
adjunction equality \eqref{e:adjunction} still holds at $w$, $\Gamma
'\in \frs$. The graph $\Gamma '$ still does not have any index of type
(9), hence by induction $\Gamma '\in \frg $ and so $\Gamma '$ is full.
As before, this fact implies $\vert v ^2 \vert \leq 2$, hence either
$v=E_{i_v}$ (contradicting the minimality of $\Gamma$) or
$v=E_{i_v}-E_j$. Then $v'=E_j -\sum _{k\in J_{v'}}E_k$, hence after
deleting $v$ the vector $v'$ in $\Gamma '$ will not be full, providing
the desired contradiction.

In the last case to consider, when $-E_{i_v}$ is in $v'$, we need to
refine the above argument in order to prove minimality for $\Gamma
'$. Modify $\Gamma$ first by adding $v$ to $v'$; notice that with this
move we eliminated the index $i_v$, hence the resulting $\Gamma ''$
embeds into $(n-1)\langle -1 \rangle$. Since $(v+v')^2 \leq v^2\leq
-2$ and the adjunction equality \eqref{e:adjunction} holds for $v+v'$,
it follows that $\Gamma ''$ is in $\frs$ involving no index of type
(9), hence by induction $\Gamma ''\in \frg $ and $v_1,v_2$ (containing
$-E_t$) are both leaves in $\Gamma ''$, although $v_1$ was not a leaf
in $\Gamma$. This shows that $v'=v_1$, hence by deleting $v$ from
$\Gamma$ and $-E_{i_v}$ from $v'=v_1$ we get a minimal graph $\Gamma
'$, since $v_1=v'$ still contains $-E_t$. Therefore the argument given
for the case above applies and provides the desired contradiction,
completing the proof.
\end{proof}
For the proof of Theorem~\ref{t:main} we also need 
\begin{prop}\label{p:g-s}
If $\Gamma \in \frg$ then $\Gamma \in \frs$.
\end{prop}
\begin{proof}
  For the graph $\Gamma \in \frg$ with $\vert \Gamma \vert =1$ the statement
  is obvious: since $\Gamma $ is full, it is $(-4)$ and so the embedding
  $v\mapsto -2E$ is suitable. If $\Gamma \subset (n-1)\langle -1 \rangle =
  \oplus _{i=1}^n \langle E_i \rangle$ with $-E_t$ in $v_1, v_2$ on the two
  ends of the plumbing chain (recall that $t$ is of type (6)), then the new
  vector $E_t-E_n$ and the modification $v_2-E_n$ provides an appropriate
  embedding of the next element $\Gamma _{next}\in \frg$ into $n \langle -1
  \rangle$. (See also the arguments in Sections~\ref{s:kirby} and
  \ref{s:smoothing}.)
\end{proof}

\section{Classification of graphs of Case {\bf {(B)}}}
\label{s:caseb}

Before examining the graphs listed under Case {\bf {(B)}}, we need a little
preparation.  

\subsection{Generalities for embeddings of Cases {\bf {(B)}} and {\bf {(C)}}}
Notice first that a simple case--by--case check shows that on 4
vertices there are exactly three graphs in $\frs$ with embeddings into
$4\langle -1 \rangle$ of Cases {\bf {(B)}} and {\bf {(C)}}; these
configurations are shown by Figure~\ref{f:n=4}. Notice that these graphs are
all full and the graph of Figure~\ref{f:n=4}(a) is both in $\frw$ and in
$\fra$; similarly Figure~\ref{f:n=4}(b) belongs to both the families $\frn$
and $\frb$, and Figure~\ref{f:n=4}(c) to both $\farm$ and $\frc$.

\begin{defn}
{\rm A vector $v=E_{i_v}-\sum _{j\in J_v}E_j$ in the minimal plumbing
tree $\Gamma$ is called \emph{reducible} if the further $-E_{i_v}$'s
are contained exactly by the vectors adjacent to $v$.  That is,
$Q(v,w)=1$ if and only if $Q(E_{i_v},w)=1$.}
\end{defn}
Note that for a reducible vector $v=E_{i_v}-\sum _{j\in J_v}E_j$ we
have that $Q(-\sum _{j\in J_v}E_j, w)=0$ for all $w\neq v$ in
$\Gamma$: If $Q(v,w)=0$ then $w$ contains no $-E_{i_v}$ by the
reducibility of $v$, hence $Q(-\sum _{j\in J_v}E_j, w)=Q(v,w)=0$.  If
$Q(v,w)=1$ then
\[
1=Q(v,w)=Q(E_{i_v},w)+Q(-\sum _{j\in J_v}E_j, w)=1+ Q(-\sum _{j\in J_v}E_j, w),
\]
which implies the statement. For example the central (that is, of
valency three) vertices of the graphs of Figure~\ref{f:n=4} under the
embeddings given by Propositions~\ref{p:c-s} and \ref{p:sok-s} are
reducible. Notice that the valency of a reducible vector $v$ is 1, 2
or 3, and in the last case $v=E_1-\sum _{j\in J_1}E_j$ where 1 is the
index of type (9).  The presence of a reducible vector will be
convenient when applying inductive proofs.

\begin{prop}\label{p:full}
Consider the graph $\Gamma\in \frs $ with indices either of Case {\bf {(B)}}
or of Case {\bf {(C)}}.  Suppose that in $\Gamma$ all leaves are full. Then
$\Gamma $ is full.
\end{prop}
\begin{proof}
  We start with the observation that such a graph admits a unique
  vertex of valency three and all others are of valency $\leq 2$. This
  follows from the fact that if $1$ denotes the index of type (9) then
  every leaf $v\in \Gamma $ must contain $-E_1$; otherwise, the sum
  $\sum _{u\neq v}u$ would have both $-2E_{i_v}$ and $-2E_1$.  (The
  case $i_v=1$ is easily excluded.)  Therefore there are at most three
  leaves and if there are only two, then the sums of the two sides of
  $v_1=E_1-\sum _{j\in J_1}E_j$ are $-2E_1-\sum _{j\in J'}E_j$ and
  $E_i-E_1-\sum _{j\in J''}E_j$, and these two vectors pair
  negatively.  In addition, if $v=E_{i_v}-\sum _{j\in J_v}E_j$ is a
  full leaf, then one $-E_{i_v}$ is in another leaf: by considering
  the sum $u$ of all vectors which are not leaves we see that $u$
  contains $E_1$ with multiplicity $+1$ and so it cannot contain
  $-2E_{i_v}$.

Turning to the proof of the proposition, suppose first that there is a
reducible vector $v=E_{i_v}-\sum _{j\in J_v}E_j$ in $\Gamma$.  Since
all leaves are full, no reducible vector of valency 1 exists in
$\Gamma$.  If $v$ is of valency two, then change $v$ to $E_{i_v}$,
blow down $E_{i_v}$ and add $v-E_{i_v}=-\sum _{j\in J_v}E_j$ to the
potential $(-1)$--vertex.  (If there is no $(-1)$--vertex in the blown
down graph then add it to any vector of the graph.)  In this way we
get a minimal tree in $\frs$, still of Case {\bf {(B)}} or {\bf
{(C)}}, with full leaves, therefore induction applies and so the blown
down graph is full. Since we changed the sum of the squares by 3, it
shows that $\Gamma$ was full.  If the reducible vector is of valency
three, then it is equal to $v_1=E_1-\sum _{j\in J_1}E_j $ (where $1$
is the index of type (9)), and since the leaves are full, we have
$n=4$: any leaf must contain $-E_1$ but since $v_1$ is reducible, all
vectors containing $-E_1$ are adjacent to $v_1$.  The graphs in $\frs$
on 4 vertices (given by Figure~\ref{f:n=4}) are known to be full.

  In the case when there is no reducible vector we need a lemma:
\begin{lem}\label{l:noredv}
Let $w\in \Gamma $ be a vector different from the central vector $v_1$
(of valency three) and let the vector next to $w$ in the path
connecting $w$ to the central vector $v_1$ be denoted by $w'$.
Suppose that there is no reducible vector on the path between $w$ and
$v_1$.  Then $-E_{i_w}$ is in $w'$.
\end{lem}
\begin{proof}
  If $w_1=w, w_2=w', w_3, \ldots , w_k=v_1$ denote the vectors on the path
  connecting $w$ with $v_1$ then $w_{k-1}$ does not admit $E_{i_{v_1}}=E_1$
  (since besides $v_1$ this basis vector is in leaf only), so $w_k=v_1$ must
  contain $-E_{w_{k-1}}$ in order $Q(w_{k-1}, w_k)=1$ to hold. Since $w_{k-1}$
  is not reducible, the potential other $-E_{w_{k-1}}$ cannot be in $w_{k-2}$,
  implying that $w_{k-1}$ must contain $-E_{w_{k-2}}$. Applying the same
  principle we arrive to the fact that $w'=w_2$ contains $-E_{i_w}$.
\end{proof}
Returning to the proof of Proposition~\ref{p:full}, suppose now that
$\Gamma$ contains no reducible vector.  Take a leaf $u$ on a leg of
length $>1$ and denote the vertex adjacent to $u$ by $u'$ (which, by
the length assumption, is different from $v_1$).  Since there is no
reducible vector on the leg of $u$, according to Lemma~\ref{l:noredv}
we get that $-E_{i_u}$ is in $u'$.  Therefore by replacing the edge
${\overline {uu'}}$ with $u+u'$ and deleting $-E_{i_u}$ from another
leaf we get a graph $\Gamma '$, and since $\Gamma '$ embeds into
$(n-1)\langle -1 \rangle = \oplus _{i \neq i_u} \langle E_i \rangle$
and it is obviously minimal (since leaves contain $-E_1$), we get that
$\Gamma '\in \frs$.  Since it still has an index of type (9), it is of
Case {\bf {(B)}} or {\bf {(C)}}. In order to apply induction, we only
need to show that $\Gamma '$ has full leaves, and since we changed
only one leaf, it boils down to showing that $u+u'=E_{i_{u'}}-\sum _j
E_j$ is full.  In other words, we have to find two $-E_{i_{u'}}$ in
$\Gamma$. The first one is easy to find, since the lack of reducible
vectors (by Lemma~\ref{l:noredv}) implies that $-E_{i_{u'}}$ is in the
vector $u_3\neq u$ adjacent to $u'$.  To find the second one, we need
a longer argument. Let $v_l\in \Gamma$ denote the leaf containing
$-E_{i_u}$.  If the leg of $v_l$ is of length $>1$, then (again by the
above principle) one $-E_{i_{v_l}}$ must be in the (unique) vector
adjacent to $v_l$, and the other one (which exists, since $v_l$ is
full) in a leaf, hence there is none in $u'$, although $Q(u',v_l)=0$.
Since both vectors contain $-E_{i_u}$, this implies that $-E_{i_{u'}}$
is in $v_l$, finishing the argument. We still have to deal with the
case when the leg of $v_l$ is of length 1, that it, $v_l$ is adjacent
to the central vertex $v_1$. If the third leg of $\Gamma$ (i.e. the
one with leaf $w$ distinct from $u$ and $v_l$) is of length $>1$, then
using this $w$ as the starting vector instead of $u$ the inductive
argument proceeds.  Therefore we only need to consider the last case
when two legs of $\Gamma$ are of length one. Suppose that the leaves
of these short legs are $v_l$ and $w$ as above. Then, as above, we can
assume that $-E_{i_u}$ is in $v_l$, and if $-E_{i_{u'}}$ is not in
$v_l$ then $-E_{i_{v_l}}$ must be in $u'$.  This implies that
$-E_{i_{v_l}}$ cannot be in $u$, hence it must be in $w$, and so
either $-E_{i_{u'}}$ is in $w$ (finishing our search) or $-E_{i_w}$ is
in $u'$, which contradicts the fact that $Q(u,u')=1$ since $-E_{i_w}$
also contributes to $u$. This final (long) observation concludes the
proof that $u+u'$ is full, hence by induction we can assume that
$\Gamma '$ is full, and since when forming $\Gamma '$ from $\Gamma$ we
increased the sum of squares by 3, this implies that $\Gamma$ is full.
\end{proof}

A useful consequence of the above result provides an easy
characterization of full graphs in $\frs$:
\begin{cor}\label{c:newfull}
The plumbing graph $\Gamma \in \frs$ of Cases {\bf {(B)}} or {\bf {(C)}}
is full if and only if $\Gamma$ has three leaves.
\end{cor}
\begin{proof}
The first paragraph of the proof of Proposition~\ref{p:full} shows
that if $\Gamma$ is full (hence admits only full leaves) then $\Gamma$
has three leaves.  Conversely, if $\Gamma$ has three leaves and all
leaves are full, then by Proposition~\ref{p:full} the plumbing graph
$\Gamma$ is full.  Consider now the case when $\Gamma$ has three
leaves and $v$ is a nonfull leaf. By denoting the unique vector
adjacent to $v$ by $v'$ we distinguish two cases: If $-E_{i_v}$ is
not in $v'$ then by erasing $v$ and adding $E_{i_v}$ to the vector $w$
$(\neq v')$ containing $-E_{i_v}$ we get a minimal graph with three
leaves on less vertices, hence induction and the fact $v^2\leq -2$
implies that $\Gamma$ is full. If $-E_{i_v}$ is in $v'$, then the
above argument breaks down only if $v'=E_{i_{v'}}-E_{i_v}$, i.e.,
after deleting $v$ and $-E_{i_v}$ the resulting plumbing graph is not
minimal. Now adding $v$ to $v'$ (and applying induction) we get that
the resulting graph $\Gamma '$ is full. This, however, means that the
other two leaves $w, u$ of $\Gamma$ are full and (by denoting the type
(9) index by 1) the basis element $-E_1$ is in $v$ (since it is in
$v+v'$ but $v'=E_{i_{v'}}-E_{i_v}$). Since $u$ and $w$ are full, a
copy of $-E_1$ is in both $u$ and $w$. Now $Q(v,u)=Q(v,w)=0$ implies
that $-E_{i_u}$ and $-E_{i_w}$ are in $v$, and so are in $v+v'$. In
this case $Q(v',u)=Q(v',w)=0$ shows that both $u$ and $w$ contain
$-E_{i_{v'}}$. In this point we reached the desired contradiction,
since a vector other than $v$ can connect to $v'$ only through
$-E_{i_{v'}}$. This final contradiction shows that a graph with three
leaves is full.
\end{proof}

\subsection{Embeddings with indices of Case {\bf {(B)}}}
We are now ready to consider graphs $\Gamma \in \frs$ with indices of
Case {\bf {(B)}}. As before, we denoted the type (9) index by 1, while
the type (2) index will be denoted by 2. Let $v$ denote the vector
containing $-E_2$.  In the classification we follow the strategy of
first proving that a graph $\Gamma $ of Case {\bf {(B)}} always admits
a reducible vector $v$; this vector always contains $-E_2$; and,
either $v$ is of valency three (in which case we show that $\vert
\Gamma \vert =4$ and the graph is given by Figure~\ref{f:n=4}(c)), or
we can contract it and apply induction. 

\begin{prop}\label{p:red2}
Suppose that $\Gamma $ is of Case {\bf {(B)}} and $v$ contains $-E_2$
(where 2 is the index of type (2)). Then for $v$ we have $v^2=-2$; in
particular, $v=E_{i_v}-E_2$ is reducible.
\end{prop}
\begin{proof}
We proceed by induction on $\vert \Gamma \vert $. For $n=4$ a
case--by--case check verifies the result. Suppose that the statement is
verified for all graphs on less than $n$ vertices, and $\Gamma $ is of
Case {\bf {(B)}} defined on $n\geq 5$ vertices.

Suppose first that $\Gamma$ has a nonfull leaf $u=E_{i_u}-\sum _{j\in
J_u}E_j $, which is adjacent to $u'$.  Assume first that $-E_{i_u}$ is
not in $u'$.  If $u$ contains $-E_2$ then either $u^2=-2$, verifying
the result, or we can move $-E_2$ to any other vector and keep the
resulting graph minimal.  Therefore in the following we can assume
that $-E_2$ is not in $u$.  Delete $u$ from $\Gamma$ and the single
$-E_{i_u}$ from the vector $w\in\Gamma$ containing it (hence $w\neq
u'$), and get $\Gamma '$ embedded in $(n-1)\langle -1 \rangle $.  The
usual argument shows that $\Gamma '$ is minimal: if
$w=E_{i_w}-E_{i_u}$ then $Q(u,w)\neq 0$.  Also, since $-E_2$ is in
$\Gamma ' \in \frs$, it is obviously of Case {\bf {(B)}}. Therefore by
induction $-E_2$ is in a vector of square $-2$ in $\Gamma '$.  The
basis element $-E_2$ cannot be in $w\in \Gamma $, since then by
induction $(w+E_{i_u})^2=-2$ and so $Q(u,
E_{i_w}-E_{i_u}-E_2)=Q(u,w)\neq 0$.  This shows that the vector
containing $-E_2$ was not changed when creating $\Gamma '$ from
$\Gamma$ hence the statement of the proposition follows by induction
on $\Gamma$.  We have to examine the case when $u'$ contains
$-E_{i_u}$. Then adding $u$ to $u'$ we can apply induction again.  If
$-E_2$ is in $u+u'$ then by induction $(u+u')^2=-2$, implying
$u^2=(u')^2=-2$.  If $-E_2$ is not in $u+u'$ then $-E_2$ is in a
vector not altered by adding $u$ to $u'$, implying the result.

Suppose next that all leaves are full.  Then by the proof of
Proposition~\ref{p:full} the graph $\Gamma$ has exactly three leaves
and (again by Proposition~\ref{p:full}) $\Gamma$ is full.  Assume that
$n >4$ and let $u$ be a leaf on the longest leg (of length $>1$),
adjacent to $u'$.  Assume first that on this leg all vectors are of
the form $w=E_{i_w}-\sum _{j\in J_w}E_j$.  By Lemma~\ref{l:noredv} if
$u'$ does not contain $-E_{i_u}$, then on the leg there is a reducible
vector $v=E_{i_v}-\sum _{j\in J_v}E_j$ of valency two. Modifying $v$ to
$E_{i_v}$ we can blow it down, and in the blown--down graph there is
at most one $(-1)$--vertex, to which we can add $v-E_{i_v}$, hence
forming a minimal graph.  Note that the blown down graph $\Gamma '$
has also three leaves, and since before adding $v-E_{i_v}$ it cannot
be full, by Corollary~\ref{c:newfull} the graph $\Gamma '$ must
contain a vector of square $(-1)$; we added $v-E_{i_v}$ to this
vector.  Induction now shows that in this graph $-E_2$ is in a vector
of square $-2$. This vector was either intact by our operation, or it
is one of the two affected by the blow down.  In the first case the
proof is obviously complete, in the second case we see that the two
vectors (after the blow--down) cannot have the shape $E_a-E_2$ and
$E_b+v-E_{i_v}$ since these two vectors give zero pairing with each
other. This shows that $v-E_{i_v}=-E_2$, verifying the statement.  

If $u'$ contains $-E_{i_u}$, then the other $-E_{i_u}$ must be in one
of the leaves of the graph, cf. the proof of Proposition~\ref{p:full}.
Consider the graph $\Gamma '$ we get by adding $u$ to $u'$ and
deleting $-E_{i_u}$ from the other leaf. Since $\Gamma$ was full, the
graph $\Gamma '$ is also full, and so has full leaves.  Then $-E_2$ is
in a vector of square $-2$ in $\Gamma '$, which cannot be the leaf we
modified, since $-E_1$ is there, hence it would have square $-3$ in
$\Gamma '$. Therefore $-E_2$ is either in a vector not changed by our
reduction or in $u+u'$, both cases implying the result.

Finally we have to consider the case when on the longest leg we chose
there is a vector $w$ of the form $-2E_3-\sum _{j\in J_w}E_j$.  In
this case $-E_1$ must be in $w$, since otherwise the sum of vectors
between $w$ and the central vector $v_1$ (incuding both) would contain
$E_1$ and $-2E_j$ for some $j$. (There is no $-E_1$ between $w$ and
$v_1$ since $-E_1$ must be in leaf.) Therefore $w$ is a leaf. If $w'$
is the unique vector adjacent to $w$ then $-E_{i_{w'}}$ is in $w$,
hence the argument given in Lemma~\ref{l:noredv} shows the existence
of a reducible vector on the leg of $w$. From that fact the proof
given above applies. Notice that if there is no vector between $w$ and
$v_1$ then (since we chose the longest leg of $\Gamma$) we get that
$\vert \Gamma \vert =4$, which case has been already discussed.
\end{proof}
As a final preparatory result, we need to examine the case when
$\Gamma $ has a reducible vector of valency three.

\begin{figure}[ht]
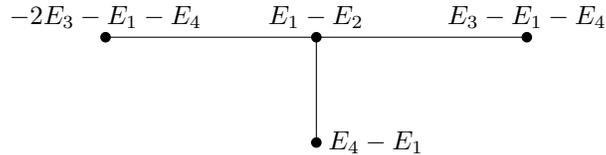

\begin{center}
\setlength{\unitlength}{1mm}
\unitlength=0.7cm
\begin{graph}(13,3)(-2,-2)
\graphnodesize{0.2}

 \roundnode{m1}(0,0)
 \roundnode{m2}(4,0)
 \roundnode{m3}(8,0)
 \roundnode{m4}(4,-2)

\edge{m1}{m2}
\edge{m2}{m3}
\edge{m2}{m4}

  \autonodetext{m1}[n]{{\small $-2E_3-E_1-E_4$}}
  \autonodetext{m2}[n]{{\small $E_1-E_2$}}
  \autonodetext{m3}[n]{{\small $E_3-E_1-E_4$}}
  \autonodetext{m4}[e]{{\small $E_4-E_1$}}
\end{graph}
\end{center}
\caption{\quad The embedding of the first element of $\frc$ into
$4\langle -1\rangle$}
\label{f:cemb}
\end{figure}

\begin{lem}\label{l:red3caseb}
Suppose that $\Gamma$ is of Case {\bf {(B)}}
and it has a reducible vector of valency three. Then $\vert \Gamma \vert =4$
and $\Gamma $ is isomorphic to the graph given by Figure~\ref{f:n=4}(c). 
\end{lem}
\begin{proof}
Obviously $v_1=E_1-\sum _{j\in J_1}E_j$ can be the only reducible
vector of valency three.  Suppose that it is adjacent to
$u_1,u_2,u_3$.  Then by Corollary~\ref{c:alttree} we have that one of
the vectors, say $u_1$ must be of the form $u_1=-2E_3-\sum _{j\in
J_{u_1}}E_j$. (Recall that in Case {\bf {(B)}} the graph contains an
index of type (5).) The shape of $v_1+u_1+u_2+u_3$ implies that
$u_2=E_3-\sum _{j\in J_{u_2}}E_j$.  If $u_3=E_4-\sum _{j\in
J_{u_3}}E_j$ then a straightforward argument (based on the fact that
$Q(u_i,u_j)=0$ for $1\leq i\neq j \leq 3$) shows that $-E_4$ must be
in $u_1$ and $u_2$ (since $u_1,u_2,u_3$ all contain $-E_1$ by the
reducibility assumption).  Let $w=v_1-E_1$; by the reducibility of
$v_1$ we have that $Q(w,u)=0$ for all $u$ distinct from $v_1$.
Suppose that $u\in \Gamma$ contains $-E_2$.  If $u$ is among $\{ v_1,
u_1, u_2, u_3 \}$ then delete $-E_2$ from it, otherwise modify $u$ to
$u+E_2+w$. Finally contract the subgraph $\{ v_1, u_1, u_2, u_3 \}$ to
the single vertex $E_1+u_1+u_2+u_3+E_3+E_4$. In this way we get a new
minimal graph $\Gamma '\in \frs$ on $n-3$ vertices, embedded into
$(n-3)\langle -1 \rangle$, since the indices 2, 3, and 4 have been
eliminated. Notice that $\Gamma '$ contains an index of type (3) (the
index 1 which was of type (9) in $\Gamma $), hence by
Theorem~\ref{t:notype3} (cf. Remark~\ref{r:ifandonlyif}) we get that
$\vert \Gamma ' \vert =n-3=1$.  The classification of graphs on 4
vertices then concludes the proof.
\end{proof}
Now we are ready to prove the classification result for this case:

\begin{proof}[Proof of Theorem~\ref{t:caseb}]
Let $v \in \Gamma $ be the vector containing $-E_2$. By
Proposition~\ref{p:red2} we know that $v+E_2=E_{i_v}$ is of square
$-1$ and therefore it is reducible.  If $v$ is of valency 1 or 2 then
we can blow $v+E_2=E_{i_v}$ down, in the resulting graph $\Gamma '$
there is no $E_{i_v}$, and in addition there is at most one vector of
square $-1$. By adding $-E_2$ to this potential $(-1)$--vector
(resulting $\Gamma ''$), by induction we have $\Gamma ''\in \frc $;
moreover the vector we added $-E_2$ to has square $-2$, hence before
adding $-E_2$ it was of square $-1$, implying that $\Gamma \in
\frc$. If the valency of $v$ is three then the application of
Lemma~\ref{l:red3caseb} concludes the proof.
\end{proof}
As in Proposition~\ref{p:g-s}, we also need that elements of $\frc$ actually
satisfy the conditions defining $\frs$.
\begin{prop}\label{p:c-s}
If $\Gamma \in \frc $ then $\Gamma \in \frs$.
\end{prop}
\begin{proof} 
  The vectors $-2E_3-E_1-E_4, E_1, E_3-E_4-E_1, E_4-E_1$ provide an
  appropriate embedding of the graph of Figure~\ref{f:alap}(c) into $3\langle
  -1 \rangle$.  After $n$ blow--ups the resulting graph on $(n+4)$ vertices
  embeds into $(n+3)\langle -1 \rangle$, and by adding $-E_2$ to the unique
  $(-1)$--vertex, the proof is complete. Notice, for example, that the graph
  of Figure~\ref{f:n=4}(c) embeds into $4 \langle -1 \rangle$ as
  $-2E_3-E_1-E_4, E_1-E_2, E_3-E_4-E_1, E_4-E_1$, cf. Figure~\ref{f:cemb}.
\end{proof}

\section{Classification of graphs of Case {\bf {(C)}}}
\label{s:casec}
Suppose now that $\Gamma$ embeds into $n \langle -1 \rangle$ with indices
listed under Case {\bf {(C)}}. Notice first that in this case for every index
$i$ there is a unique vector $v_i$ such that
\[
v_i=E_i-\sum _{j\in J_i}E_j.
\]
In the following we will distinguish two cases, according to whether 
the leaves of $\Gamma $ are full or not.

\subsection{Graphs with full leaves}

We start with the case when all leaves of the plumbing tree are full.
As before, $i=1$ denotes the unique index of type (9); there are
exactly three leaves, each containing $-E_1$, and if $v=E_{i_v}-\sum
_{j\in J_v}E_j$ is a leaf then one $-E_{i_v}$ is contained by another
leaf, cf. the proof of Proposition~\ref{p:full}.  The unique vector
$v_1=E_1-\sum _{j\in J_1}E_j$ is of valency three and all other
vectors have $\deg v \leq 2$.

\subsubsection*{Graphs with reducible vectors}

First we will  examine the case  when the leaves are full, and
the graph contains reducible vectors.

\begin{thm}
Suppose that $\Gamma $ is a graph with indices of Case {\bf {(C)}},
all leaves are full and there is a vector which is reducible.
Then $\Gamma \in \fra \cup \frb$.
\end{thm}
\begin{proof}
We will proceed by induction.  Consider the reducible vector
$v=E_{i_v}-\sum _{j \in J_v}E_j$. Since the leaves are full, a leaf
cannot be reducible.  If $v$ is of valency two, then change it to
$E_{i_v}$, blow it down and consider the resulting graph $\Gamma
'$. Since $\Gamma '$ with three leaves cannot be full (we erased $v$
and two further $-E_{i_v}$), by Corollary~\ref{c:newfull} it is not
minimal. Add $v-E_{i_v}=-\sum _{j\in J_v}E_j$ to the $(-1)$--vector
and get $\Gamma ''$.  This is a minimal graph, the leaves are still
full, by construction it has a reducible vector (the $(-1)$--vector to
which we added $-\sum _{j\in J_v}E_j$), hence induction applies,
showing that $\Gamma ''\in \fra \cup \frb$. Now $\Gamma $ is gotten
from $\Gamma ''$ by considering the reducible vector $w=E_{i_w}-\sum
_{j\in J_w}E_j\in \Gamma ''$, adding $\sum _{j\in J_w}E_j$ to it,
blowing it up and then subtracting $\sum _{j\in J_w}E_j$ from the
resulting new vertex with decoration $(-1)$.

The blowing down procedure will stop only when $\deg v =3$ for the
    reducible vector $v$. Since the leaves are full, in this case we
    have $n=4$, and we get that the resulting graph is one of
    Figures~\ref{f:n=4}(a) or (b). Since in that case the square of
    the correction term $-\sum _{j\in J_w}E_j$ we added to the
    $(-1)$--vertex is $-2$ or $-3$, it must have been the case
    throughout the whole process, since the two (or three) basis
    vectors cannot be grouped in two groups such that they give zero
    pairing with every other vector. This last observation completes
    the proof.
\end{proof}

\subsubsection*{Graphs with no reducible vectors}

The next case to consider is when the graph $\Gamma$ has full leaves
and it contains no reducible vector.  It turns out that this is the
longest case to discuss, and we will divide our classification into
futher subcases. To do that, recall that since the graph has full
leaves, it is full, and has exactly three leaves $w_1,w_2$ and $w_3$.

\begin{figure}[ht]
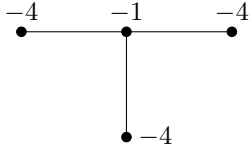

\begin{center}
\setlength{\unitlength}{1mm}
\unitlength=0.7cm
\begin{graph}(10,3)(-3,-2)
\graphnodesize{0.2}

 \roundnode{m1}(0,0)
 \roundnode{m2}(2,0)
 \roundnode{m3}(4,0)
 \roundnode{m4}(2,-2)

\edge{m1}{m2}
\edge{m2}{m3}
\edge{m2}{m4}

  \autonodetext{m1}[n]{{\small $-4$}}
  \autonodetext{m2}[n]{{\small $-1$}}
  \autonodetext{m3}[n]{{\small $-4$}}
  \autonodetext{m4}[e]{{\small $-4$}}
\end{graph}
\end{center}
\caption{\quad (Nonminimal) graph with three free leaves}
\label{f:4-es}
\end{figure}

\begin{defn}
{\rm Suppose that the graph $\Gamma$ of Case {\bf {(C)}} is full and
admits no reducible vectors.
\begin{itemize}
\item The leaves of $\Gamma$ are \emph{of type ${\bfz }_3$} if (after
possibly renaming the leaves) $-E_{i_{w_1}}$ is in $w_2$,
$-E_{i_{w_2}}$ is in $w_3$ and $-E_{i_{w_3}}$ is in $w_1$.
 
\item The leaves of the graph $\Gamma$ as above are \emph{of type
${\bfz }_2$} if (again, after possibly renaming the leaves)
$-E_{i_{w_1}}$ is in $w_2$, another $-E_{i_{w_1}}$ is in $w_3$ and
then it follows that $-E_{i_{w_3}}$ is in $w_2$ and $-E_{i_{w_2}}$ is
in $w_3$.

\item The leaf $w=E_{i_w}-\sum _{j\in J_w}E_j$ of a graph $\Gamma$
with full leaves and no reducible vectors is called \emph{free} if
the single vector $w'$ with $Q(w,w')=1$ does not contain $-E_{i_w}$.
\end{itemize}
}
\end{defn}  
It is easy to see that if a leaf $w$ is free, then (since $\Gamma$
contains no reducible vector) it is adjacent to the central vertex
$v_1$ (of valency three), that is $Q(w,v_1)=1$. It is worth mentioning
that $Q(w,v_1)=1$ might hold without $w$ being free. In the graph
$\Gamma$ there are three leaves, and if all are free then the graph is
defined on four vertices and is nonminimal (admitting the shape given
by Figure~\ref{f:4-es}, with vectors $E_2-E_1-E_3-E_4, E_1,
E_3-E_1-E_2-E_4$ and $E_4-E_1-E_2-E_3$).
In the following therefore we restrict our attention to the cases when
there are zero, one or two free leaves in $\Gamma$.
In fact, by definition if the leaves $\{ w_1, w_2, w_3\}$ are of type 
$\bfz _2$ then $w_1$ is automatically free, hence in this case we have only 
two cases to examine.  Notice first that 
\begin{lem}
Any graph $\Gamma \in \frs$ of Case {\bf {(C)}} with full leaves is
either of type ${\bfz }_3$ or of type ${\bfz }_2$.
\end{lem}
\begin{proof}
Let $w_1=E_{i_{w_1}}-\sum _{j\in J_{w_1}}E_j$ be a leaf. We already know that
one $-E_{i_{w_1}}$ is in another leaf, say in $w_2$. If the third leaf $w_3$
also contains $-E_{i_{w_1}}$ then in order to have
$Q(w_2,w_3)=0$ (since both vectors contain $-E_1$ and $-E_{i_{w_1}}$) we have 
that 
\[
w_2=E_{i_{w_2}}-E_1-E_{i_{w_1}}-E_{i_{w_3}}-\ldots \mbox{ and }
w_3=E_{i_{w_3}}-E_1-E_{i_{w_1}}-E_{i_{w_2}}-\ldots ,
\]
showing that the leaves $\{ w_1,w_2,w_3\}$ are of type $\bfz _2$.

If the third leaf $w_3$ does not contain $-E_{i_{w_1}}$ then (in order
to achieve $Q(w_1,w_3)=0$) we have that $-E_{i_{w_3}}$ is in $w_1$.
To get $Q(w_2,w_3)=0$, we either have that $-E_{i_{w_3}}$ is in $w_2$,
which provides a set of leaves $\{ w_1, w_2, w_3\}$ of type $\bfz _2$
(with the indices permuted), or $-E_{i_{w_2}}$ is in $w_3$, showing that
the leaves are of type $\bfz _3$. Since we have examined all possibilities,
the proof is complete.
\end{proof}
The following two statements now classify graphs with full leaves and
no reducible vectors: The graphs of
Figures~\ref{f:wahltypeembed}--\ref{f:lambdaspec2beagy} (where the
expressions of the vertices are also indicated) might be helpful in
following the proofs.

\begin{prop}\label{p:z3k}
Suppose that $\Gamma \in \frs$ of Case {\bf{(C)}}
has full leaves, no reducible vectors
and its leaves are of type ${\bfz }_3$.
\begin{itemize}
\item
If $\Gamma$ has no free leaf then $\Gamma \in \frw$.
\item
If it has one free leaf then $\Gamma$ is isomorphic to $\Delta
_{p,q,r}$ with $p\geq 1$ and $q,r\geq 0$ (and so $\Gamma \in \frn$).
\item 
If $\Gamma$ has two free leaves then it is isomprhic to $\Lambda
_{p,q,r}$ with $p\geq 1$ and $q,r\geq 0$ (and so $\Gamma \in \farm$).
\end{itemize}
\end{prop}
\begin{proof}
  Suppose first that $\Gamma$ admits no free leaf.  Let us fix a leaf
  $w_1$ and assume that $-E_{i_{w_1}}$ is in $w_2$.  If $w_1'$ is the
  unique vector adjacent to $w_1$, then according to
  Lemma~\ref{l:noredv} the other $-E_{i_{w_1}}$ is in $w_1'$. The same
  principle also shows that if $w_1''$ is the next vector on this leg
  then it contains $-E_{i_{w_1'}}$; and so on throughout the leg. At
  the same time, since $Q(w_2,w_1')=0$, the fact that $-E_{i_{w_1}}$
  is in $w_1'$ implies that $-E_{i_{w_1'}}$ must be in $w_2$ as well:
  recall that since there are no free leaves, $-E_{i_{w_2}}$ cannot be
  in $w_1'$.  Suppose that there are $p$ vectors $u_1, \ldots , u_p$
  on the leg of $w_1$ (not counting $w_1$ and the central vector
  $v_1$). Then the above argument shows that in $w_2$ we have all the
  basis vectors $-E_{i_{u_j}}$, next to $-E_{i_{w_1}}$, $-E_1$ and
  $E_{i_{w_2}}$. This shows that $w_2^2\leq -p-3$. The same argument
  now applies to all the three legs of $\Gamma$. The central element
  $v_1=E_1-\sum _{j\in J_1}E_j$ must be adjacent to the three legs,
  and since all three $E_1$'s are in leaves and the leaves are not
  free, it follows that $v_1^2\leq -4$. A simple count on the possible
  value of $\sum v_i^2$ shows that the estimates sum up to $-3\vert
  \Gamma \vert -1$ (its smallest possible value) therefore the above
  inequalities must be equalities, and the shape of $\Gamma $ is given
  by Figure~\ref{f:wahltype} for some appropriate $p,q,r\geq 0$.  This
  concludes the proof of the first case.

Suppose now that $\Gamma$ admits one free leaf, say $w_1$. This
condition means that one $-E_{i_{w_1}}$ is in a leaf, the other is not
(since the graph is not of $\bfz _2$--type), and it is also not in the
central vector.  Suppose that it is in $u$ on one of the legs.  The
leg is determined in the following way: if $-E_{i_{w_1}}$ is in $w_2$
then the vector $u$ must be on the leg of $w_3$ otherwise the absence
of reducible vectors would imply $Q(w_1,v_1)=0$, a contradiction.  The
vector $u$ divides this leg into two pieces, having $p-1$ vector on
the side towards the leaf, and $q$ towards the central element (not
counting the leaf and the central vector) for some $p\geq 1$ and
$q\geq 0$.  The same estimates as before can be used to estimate the
squares of the leaves, resulting $w_1^2\leq -p-2$, $w_2^2\leq -q-4$
and $w_3^2\leq -r-3$ (where $r$ denotes the number of vertices on the
leg not containing $u$).  The central vector must satisfy only
$v_1^2\leq -3$ since there is one free leaf (which can connect to
$v_1$ by sharing $-E_1$). On the other hand, since $u$ contains
$-E_{i_{w_1}}$ we get that $u^2\leq -3$. Again, the sum of the
estimates already sums up to $-3\vert \Gamma \vert -1$, showing that
all inequalities must be equalities and so $\Gamma$ is of the form
$\Delta _{p,q,r}$. Since the graph is of type $\bfz _3$, the basis
element $-E_{i_{w_1}}$ cannot be in the leaf $w_2$, therefore we need
$p\geq 1$. (The degeneration $p=0$ will lead to a graph of type $\bfz
_2$, cf.  Proposition~\ref{p:z2k}.)

Finally, the presence of two free leaves $w_1,w_2$ means that two legs
are degenerated to a single leaf, hence the second copies of
$-E_{i_{w_1}}$ and $-E_{i_{w_2}}$ are on the same (nontrivial) leg. An
argument similar to the previous one shows that if $-E_{i_{w_2}}$ is
in the leaf of the long leg, then the basis vector $-E_{i_{w_1}}$ is
then one closer to the leaf from the two ``free'' basis elements
$-E_{i_{w_1}}$ and $-E_{i_{w_2}}$. The exact same analyis as above now
shows that if these two basis vectors are in different vertices, then
$\Gamma$ is of the form $\Lambda _{p,q,r}$ with $p,r\geq 1$; and
finally $\Gamma$ is isomorphic to $\Lambda _{p,q,0}$ if both
$-E_{i_{w_1}}$ and $-E_{i_{w_2}}$ are in the same vector. Notice that
since $\Gamma$ is of type $\bfz _3$, this vector cannot be a leaf. The
degeneration when both these basis elements are in a leaf is exactly
$\Lambda _{0,q,0}$, which case --- according to our point of view ---
falls in the category of Case {\bf {(B)}} rather than Case {\bf
{(C)}}. Now this vector $u$ containing both $-E_{i_{w_1}}$ and
$-E_{i_{w_2}}$ will satisfy $u^2\leq -4$, hence the usual argument of
summing the squares provides the result.
\end{proof}
\begin{prop}\label{p:z2k}
Suppose that $\Gamma \in \frs$ of Case {\bf {(C)}} has full leaves, no
reducible vectors and its leaves are of type ${\bfz }_2$.
\begin{itemize}
\item
If it has one free leaf then $\Gamma$ is isomorphic to $\Delta
_{0,q,r}$ with $q,r\geq 0$, in particular $\Gamma \in \frn$.
\item 
if $\Gamma$ has two free leaves then it is isomprhic to $\Lambda
_{0,p,r}$ ($q\geq 0$ and $r\geq 1$), so $\Gamma \in \farm$.
\end{itemize}
\end{prop}
\begin{proof}
Notice that (as we already remarked) the fact that $\Gamma $ 
is of type $\bfz _2$ implies that there is a free leaf $w_1$, and 
both $-E_{i_{w_1}}$ are in the two other leaves. The same idea and the
estimate on the sum of squares shows that such a graph is isomorphic
to $\Delta _{0,q,r}$ for some $q,r\geq 0$.

If there are two free leaves $w_1,w_2$ (and $\Gamma $ is still of type $\bfz
_2$) then for $w_1$ both $-E_{i_{w_1}}$ are in leaves (as is dictated by the
fact that $\Gamma $ is of type $\bfz _2$), but the second $-E_{i_{w_2}}$ is on
the single nontrivial leg. It gives rise to a vector $u$ with $u^2\leq -3$,
and so the usual estimate on the sum of squares identifies $\Gamma$ with a
graph $\Lambda _{0,q,r}$ for $r\geq 1$ and $q\geq 0$. Again, by degenerating
$-E_{i_{w_2}}$ to the same leaf (i.e., considering $r=0$) we get a graph
$\Lambda _{0,q,0}$ of Case {\bf {(B)}}.
\end{proof}

\subsection{The case of nonfull leaves}
The final case to be considered is when the graph $\Gamma \in \frs$ has
indices of Case {\bf {(C)}} and it also admits nonfull leaves.

\begin{lem}\label{l:nemfull}
If $\Gamma $ has a nonfull leaf then it has a reducible vector.
\end{lem}
\begin{proof}
  As usual, we use induction on $\vert \Gamma \vert $. Suppose that
  $v$ is a nonfull leaf adjacent to $v'$. If $v'$ contains $-E_{i_v}$
  then $v$ is reducible, and we are done.  If $v'$ does not contain
  $-E_{i_v}$, then delete $v$ and $-E_{i_v}$ from the graph. (The term
  $-E_{i_v}$ was contained by $w$, and since $Q(v,w)=0$, it is easy to
  see that $(w+E_{i_v})^2\leq -2$, hence the resulting graph is
  minimal on less vertices.) Now if the resulting graph $\Gamma '\in
  \frs$ admits a nonfull leaf, then induction shows that $\Gamma '$
  admits a reducible vector. It is obviously reducible in $\Gamma$ as
  well, unless it is the (potential) new leaf $v'$ of $\Gamma '$, in
  which case (since $-E_{i_v}$ was not in $v'$) we get that $v'$ was
  reducible in $\Gamma$.  In case in $\Gamma '$ all leaves are full,
  then we argue as follows: there are exactly three leaves in $\Gamma
  '$, and by Proposition~\ref{p:full} $\Gamma '$ is full, showing that
  $\sum _{v_i\in \Gamma '}v_i^2=-3\vert \Gamma '\vert -1$. Since
  $\Gamma$ is nonfull, we have that $\sum _{v_i \in \Gamma }v_i^2>
  -3\vert \Gamma \vert -1=-3(\vert \Gamma '\vert +1)-1$, and since we
  have dropped $v$ and $-E_{i_v}$ from a vertex of $\Gamma$, it means
  that $v^2>-2$, contradicting the minimality of $\Gamma$.
\end{proof}
Before the final argument, we need to study the case when the
reducible vector is of valency three, providing the analogue of
Lemma~\ref{l:red3caseb}.

\begin{lem}\label{l:red3casec}
If $\Gamma \in \frs$ is of Case {\bf {(C)}} and it admits a reducible
vector of valency three then $\vert \Gamma \vert =4$ and $\Gamma$ is 
isomorphic to one of the graphs of Figure~\ref{f:n=4}(a) or (b).
\end{lem}
\begin{proof}
The proof will proceed by induction on $\vert \Gamma \vert $.  First,
a simple case--by--case check shows that such graph does not exist if
$\vert \Gamma \vert =5$. Suppose now that $\Gamma \in \frs$ admits a
reducible vector $v_1$ of valency three. Let $v\in \Gamma $ be a leaf
with $Q(v,v_1)=0$. If $\vert \Gamma \vert >5$ then such leaf exists,
and since $v_1$ is reducible, it is nonfull. Following the usual line
of reasoning, either delete $v$ and the unique further $-E_{i_v}$ (if
$-E_{i_v}$ is not in the unique vector $v'$ adjacent to $v$) or add
$v$ to $v'$. In this way we get $\Gamma '\in \frs$ on $\vert \Gamma
\vert -1$ vertices, still admitting a reducible vector of valency
three, hence induction yields that $\Gamma$ must be equal to 4. The
rest of the statement is an easy exercise.
\end{proof}

\begin{thm}
If $\Gamma $ is of Case {\bf {(C)}} and has a nonfull leaf then
$\Gamma \in \fra\cup \frb$.
\end{thm}
\begin{proof}
  By Lemma~\ref{l:nemfull} the assumption implies that there is a
  reducible vector $v\in \Gamma$. If $\deg v =3$ then by
  Lemma~\ref{l:red3casec} the theorem is proved.  If $v$ is of valency
  1 or 2, then we can modify $v$ to $E_{i_v}$, blow down $E_{i_v}$ and
  add $v-E_{i_v}$ to the potential $(-1)$--vector of the resulting
  graph $\Gamma '$. We want to argue that $\Gamma '$ also admits a
  reducible vector. If $\Gamma '$ has a nonfull leaf then this fact
  follows from Lemma~\ref{l:nemfull}. In case all leaves of $\Gamma '$
  are full, then $\Gamma '$ has three leaves and $\Gamma '$ is full,
  hence before adding $v-E_{i_v}$ it must have contained a
  $(-1)$--vector, which is automatically reducible.  By induction on
  $\vert \Gamma \vert $ we get that $\Gamma '\in \fra \cup \frb$,
  therefore in oder to prove the theorem we only need to understand
  how $\Gamma$ can be rebuilt from the blown down graph $\Gamma '$.

 Suppose that $v=E_{i_v}-\sum _{j\in J_v}E_j$ is the reducible vector
  in $\Gamma$, and $w$ is the reducible vector after the blow--down.
  We add $-\sum _{j\in J_v}E_j$ to $w$, and repeat the above
  procedure.  This procedure will terminate only when the reducible
  vector has valency three, in which case $n=4$, hence the sum of the
  tails (i.e., $v-E_{i_v}$) of the reducibles add up to one, two or
  three basis vectors. If it is one, then it is under Case {\bf
  {(B)}}, since it means that that basis vector does not appear
  anywhere else, hence is of type (2). If it is two or three, then
  these come from the first reducible vector, since a single $-E_i$
  cannot be disjoint from all the vectors (since all indices appear in
  some vector with positive multiplicity).  This shows that there are
  two cases, when $(-\sum _{j\in J_v}E_j)^2$ is $-2$ or $-3$,
  corresponding to the two cases growing from the two basic examples
  of Figure~\ref{f:alap} which are of Case {\bf {(C)}}.
\end{proof}

In order to complete the equivalence of $\frs$ with the set $\frg \cup
\frw \cup \frn \cup \farm \cup \fra \cup \frb \cup \frc$ we also need
\begin{prop}\label{p:sok-s}
If $\Gamma \in \frw \cup \frn \cup \farm \cup \fra \cup \frb$ then
$\Gamma \in \frs$.
\end{prop}
\begin{proof}
  The proof of $\fra \subset \frs$ follows the same idea as the proof of
  Proposition~\ref{p:c-s}: the vectors $E_2-E_1-E_3, E_1, E_3-E_1-E_4,
  E_4-E_1-E_2$ embed Figure~\ref{f:alap}(a) into $4\langle -1 \rangle$ such
  that the vector $-E_2-E_3-E_4$ gives zero pairing with all these vectors
  above. After $n$ blow--ups, the addition of $-E_2-E_3-E_4$ to the unique
  $(-1)$--vertex provides the embedding. The proof of the inclusion $\frb
  \subset \frs$ proceeds in the same way: notice that $E_1-E_2-E_3-E_4, E_1,
  E_3-E_1-E_2-E_4, E_4-E_1$ gives an embedding of Figure~\ref{f:alap}(b) into
  $4\langle -1 \rangle$ and $-E_2-E_3$ gives zero pairing with all the vectors
  listed above. Finally, for proving the inclusions $\frw \subset \frs$,
  $\farm \subset \frs$ and $\frn \subset \frs$ we give the explicite
  embeddings by Figures~\ref{f:wahltypeembed} through \ref{f:lambdaspec2beagy}
  of the Appendix.
\end{proof}

With this final argument the proof of Theorem~\ref{t:main} is now complete.

\section{Constructions of rational homology disks through Kirby calculus}
\label{s:kirby}
\begin{figure}[ht]
\begin{center}
\includegraphics[width=10cm]{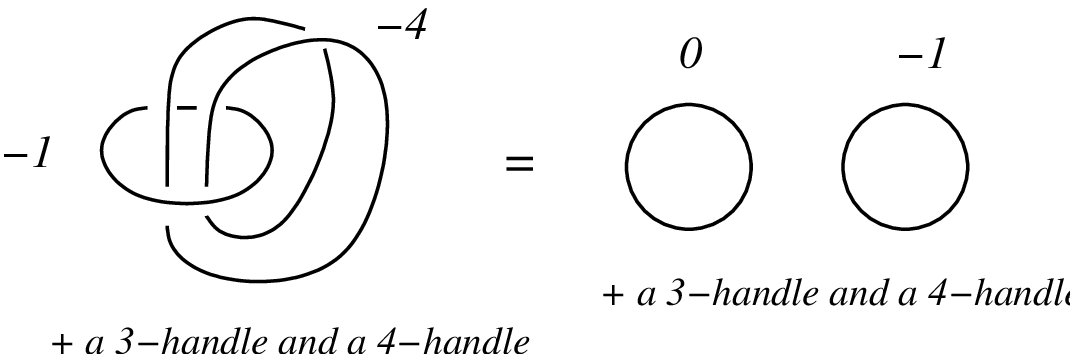}
\end{center}
\caption{The embedding of $\Gamma = (-4)$ into $\cpkk$}
\label{f:minnegy}
\end{figure}

\begin{figure}[ht]
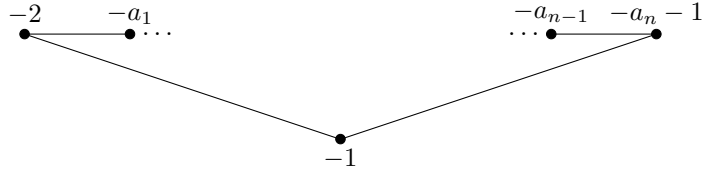

\begin{center}
\setlength{\unitlength}{1mm}
\unitlength=0.7cm
\begin{graph}(-6,5)(9,-2)
\graphnodesize{0.2}

 \roundnode{m1}(0,0)
 \roundnode{m2}(2,0)
 \roundnode{m3}(10,0)
 \roundnode{m4}(12,0)
 \roundnode{m5}(6,-2)

\edge{m1}{m2}
\edge{m3}{m4}
\edge{m4}{m5}
\edge{m1}{m5}

  \autonodetext{m1}[n]{{\small $-2$}}
  \autonodetext{m2}[n]{{\small $-a_1$}}
  \autonodetext{m3}[n]{{\small $-a_{n-1}$}}
  \autonodetext{m4}[n]{{\small $-a_{n}-1$}}
  \autonodetext{m5}[s]{{\small $-1$}}
  \autonodetext{m2}[e]{{\small $\ldots$}}
  \autonodetext{m3}[w]{{\small $\ldots$}}

\end{graph}
\end{center}
\caption{Attachment of the $(-1)$--framed 2-handle to $\mg$}
\label{f:ratbl}
\end{figure}

\begin{figure}[ht]
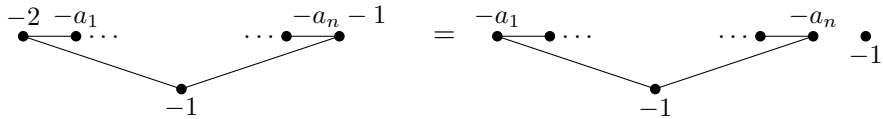

\begin{center}
\setlength{\unitlength}{1mm}
\unitlength=0.7cm
\begin{graph}(-4,5)(9,-2)
\graphnodesize{0.2}

 \roundnode{m1}(-1,0)
 \roundnode{m2}(0,0)
 \roundnode{m3}(4,0)
 \roundnode{m4}(5,0)
 \roundnode{m5}(2,-1)

\edge{m1}{m2}
\edge{m3}{m4}
\edge{m4}{m5}
\edge{m1}{m5}

  \autonodetext{m1}[n]{{\small $-2$}}
  \autonodetext{m2}[n]{{\small $-a_1$}}
  \autonodetext{m4}[n]{{\small $-a_{n}-1$}}
  \autonodetext{m5}[s]{{\small $-1$}}
  \autonodetext{m2}[e]{{\small $\ldots$}}
  \autonodetext{m3}[w]{{\small $\ldots$}}

\freetext(7,0){$=$}

 \roundnode{n1}(8,0)
  \roundnode{n2}(9,0)
 \roundnode{n3}(13,0)
\roundnode{n4}(14,0)
 \roundnode{n5}(11,-1)

\edge{n1}{n2}
\edge{n3}{n4}
\edge{n4}{n5}
\edge{n1}{n5}

  \autonodetext{n1}[n]{{\small $-a_1$}}
  \autonodetext{n4}[n]{{\small $-a_{n}$}}
  \autonodetext{n5}[s]{{\small $-1$}}
  \autonodetext{n2}[e]{{\small $\ldots$}}
  \autonodetext{n3}[w]{{\small $\ldots$}}

 \roundnode{n6}(15,0)
  \autonodetext{n6}[s]{{\small $-1$}}

\end{graph}
\end{center}
\caption{The inductive step}
\label{f:blowdown}
\end{figure}

In this section we prove Theorem~\ref{t:kirby} through embedding the
plumbing 4--manifold $\mg$ corresponding to certain $\Gamma \in \frs$
into $\#_{\vert \Gamma \vert} \cpkk$. Then the complement of $\mg$
(with the reversed orientation) will be an appropriate rational
homology disk $B _{\Gamma}$.  The embedding of $\mg$ into $\#_{\vert
\Gamma \vert} \cpkk$ will be given by Kirby calculus.  We start by
showing that a graph $\Gamma\in\frg $ gives rise to a 4--manifold
which embeds into the connected sum of $\vert \Gamma \vert$ copies of
$\cpkk$. (This statement has been already proved in many places, eg.
\cite{CH}, \cite{FS1}.) The proof uses induction: For the starting case
$\Gamma =(-4)$ the diagram of Figure~\ref{f:minnegy} shows that by
attaching a 2-handle, a 3-handle and a 4-handle to $\mg$ (which in
this case can be represented by the $(-4)$--framed unknot) we get
$\cpkk$, hence the claim follows.
\begin{rem} {\rm This picture is simply the manifestation of the fact that
    twice the generator of $H_2(\cpkk ; \bfz )$ can be represented by a sphere
    of self-intersection $(-4)$. Another way to see this diagram is to regard
    the $(-4)$-sphere as the representative of twice the generator (with
    orientation reversed) and the $(-1)$-sphere as the representative of the
    generator; then by blowing up one of the two (positive) intersections
    repeatedly we can build any graph in $\frg$ and in the meantime also
    verify the embedding we are looking for.}
\end{rem}
Now suppose that $\Gamma=(-2, -a_1, \ldots , -a_n -1)$ is constructed
from $\Gamma '=(-a_1, \ldots , -a_n)$ by the inductive step we used to
construct $\frg$ (cf. Subsection~\ref{ss:contfrac}). We claim that by
attaching a 2-handle (symbolized by the $(-1)$--vertex) to $\mg$ as
shown by the plumbing graph of Figure~\ref{f:ratbl} we get a 4-manifold
which embeds into $\# _{\vert \Gamma \vert } \cpkk$. To see this
embedding we proceed by induction again; suppose that it is true for
$\Gamma '$, and use the diffeomorphism we get by blowing down the
$(-1)$--sphere of Figure~\ref{f:ratbl}, cf. Figure~\ref{f:blowdown}.

\begin{figure}[ht]
\begin{center}
\includegraphics[width=10cm]{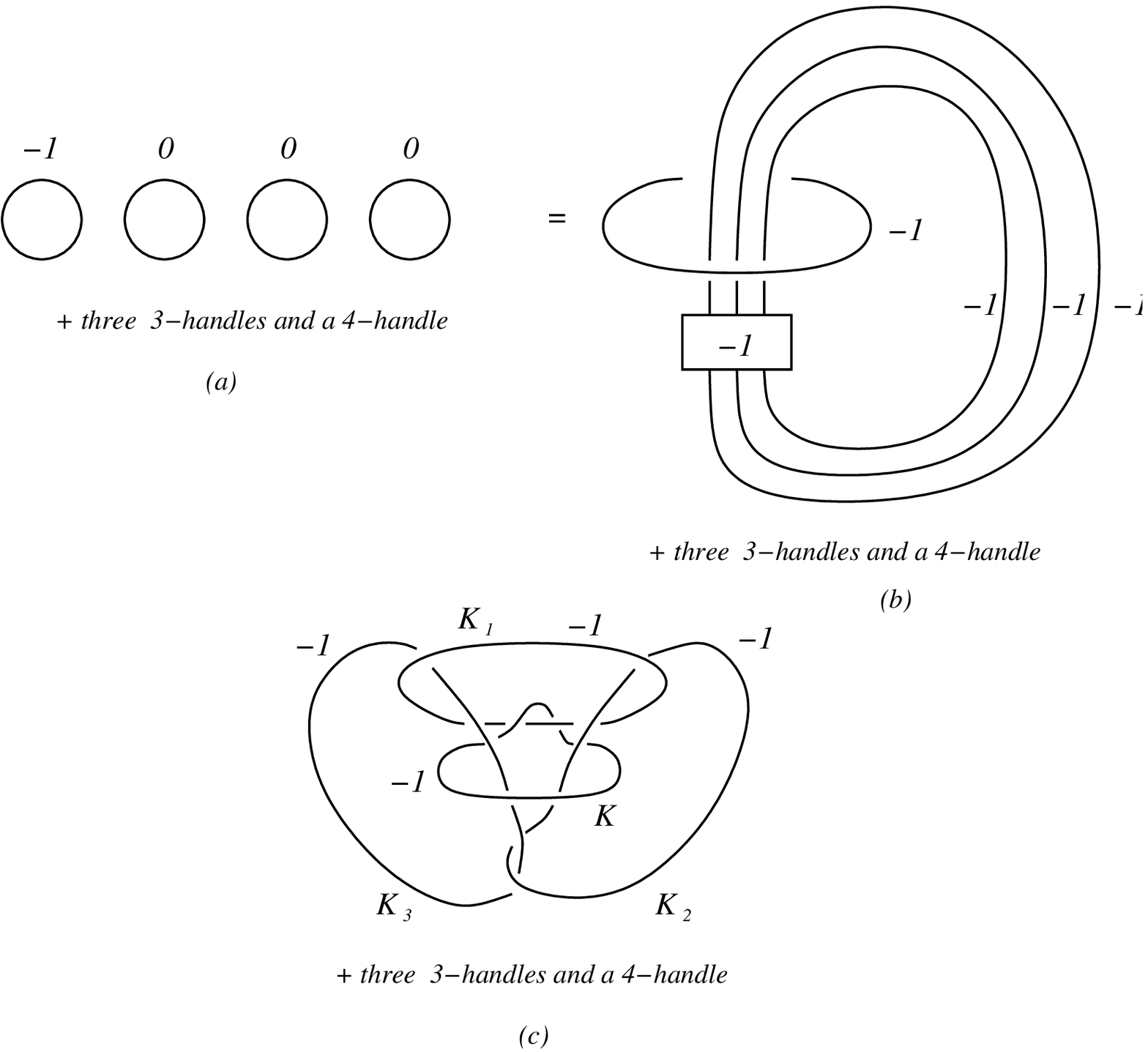}
\end{center}
\caption{Surgery diagrams for $\cpkk$}
\label{f:cpkk}
\end{figure}

\begin{figure}[ht]
\begin{center}
\includegraphics[width=10cm]{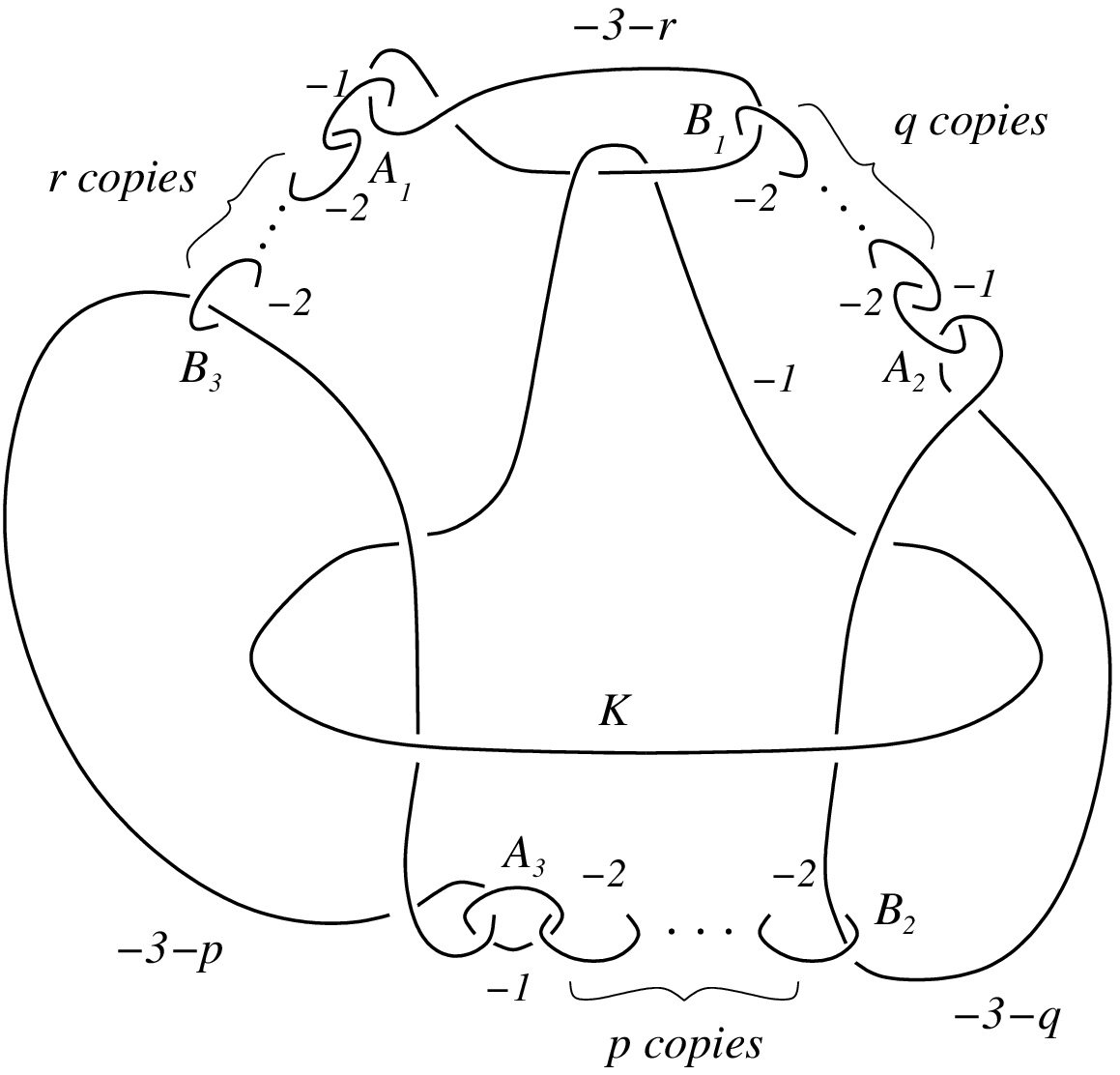}
\end{center}
\caption{The diagram of Figure~\ref{f:cpkk}(c) after $(p+q+r+3)$ blow--ups}
\label{f:embed}
\end{figure}

In order to show the necessary embedding for elements in $\farm\cup
\frn \cup \frw$, notice first that Figure~\ref{f:cpkk} provides three
Kirby diagrams for $\cpkk$: in (a) the 0--framed circles cancel the
3--handles, while (b) can be reduced to (a) by sliding one of the
$(-1)$--circles off the others; whereas (c) differs from (b) by a
simple isotopy.
Now repeatedly blowing up the linkings of the three circles $K_1, K_2$
and $K_3$ ($p+1,q+1$ and $r+1$ times, respectively) we get a
configuration of spheres in $(p+q+r+4)\cpkk$,
cf. Figure~\ref{f:embed}.
Sliding the 2--handles corresponding to the circles $A_i$ over the
handle corresponding to the central circle $K$ we get $M_{\Gamma
_{p,q,r}}\subset (p+q+r+4)\cpkk$, showing that $\partial M_{\Gamma
_{p,q,r}}$ bounds a rational homology disk.  Sliding $A_1$ over $B_1$
and $A_2, A_3$ over $K$ we get the corresponding result for the
4--manifolds given by the graphs $\Delta _{p,q,r}$. Finally, sliding
$A_1$ and $A_2$ over $B_1,B_2$ respectively, and $A_3$ over $K$ we get
the desired results for $\Lambda _{p,q,r}$. Notice that by not
applying blow--ups at some linkings of $K_1,K_2,K_3$ we get the
degeneration of $\Delta _{p,q,r}$ for $p=0$ and of $\Lambda _{p,q,r}$
for $p=0$ and/or $r=0$.  The simple details are left to the
reader. Notice that this approach provides a unified treatment of all
graphs in $\frw \cup \frn \cup \farm$.

\section{Smoothings of normal surface singularities 
with vanishing Milnor number}
\label{s:smoothing}

Next we will show two methods for constructing smoothings of normal
surface singularities which are rational homology disks.  Notice that
even for those cases treated in the previous section we will get
stronger results: in this section the rational homology disks are
smoothings of singularities, hence are equipped with natural Stein
structures; in fact, they will be affine varieties, with computable
fundamental group (and sometimes explicitly given universal covering
spaces).  In addition, the algebro-geometric method to be presented
below applies to cases which are not covered in
Section~\ref{s:kirby}. The constructions in this section will verify
Theorem~\ref{t:smoothing}.

\subsection{Smoothings of negative weights}\label{ss:negwt}
The first method is similar in spirit to the strategy we used in
Section~\ref{s:kirby}.  This method, however, has the limitation of
being applicable only for star--shaped graphs.  Suppose therefore that
$\Gamma$ is a star--shaped graph given by the diagram
$$\xymatrix@R=4pt@C=24pt@M=0pt@W=0pt@H=0pt{\\
\lefttag{\bullet}{n_2/p_2}{8pt}\dashto[ddrr]&
&\hbox to 0pt{\hss\lower 4pt\hbox{.}.\,\raise3pt\hbox{.}\hss}
&\hbox to 0pt{\hss\raise15pt
\hbox{.}\,\,\raise15.7pt\hbox{.}\,\,\raise15pt\hbox{.}\hss}
&\hbox to 0pt{\hss\raise 3pt\hbox{.}\,.\lower4pt\hbox{.}\hss}
&&\righttag{\bullet}{n_{t-1}/p_{t-1}}{8pt}\dashto[ddll]\\
\\
&&\bullet\lineto[dr]&&\bullet\lineto[dl]\\
\lefttag{\bullet}{n_1/p_1}{8pt}\dashto[rr]&&
\bullet\lineto[r]&\overtag{\bullet}{-b}{8pt}\lineto[r]&\bullet
\dashto[rr]&&\righttag{\bullet}{n_{t}/p_{t}}{8pt}\\&~\\&~}
$$ 
(We assume that $t\geq 3$, thus excluding the case of cyclic
quotient singularities.) The strings of $\Gamma$ are described uniquely by the continued
fractions shown, starting from the node.
Define its \emph{dual} $\Gamma '$ by reversing the sign
of all the decorations on $\Gamma$ and applying `handle-calculus' along the
legs to turn the positive numbers into negative again (except for the central
vertex, which might have positive decoration). More explicitely, consider 
$\Gamma '$ given by
$$\xymatrix@R=4pt@C=24pt@M=0pt@W=0pt@H=0pt{\\
\lefttag{\bullet}{n_2/q_2}{8pt}\dashto[ddrr]&
&\hbox to 0pt{\hss\lower 4pt\hbox{.}.\,\raise3pt\hbox{.}\hss}
&\hbox to 0pt{\hss\raise15pt
\hbox{.}\,\,\raise15.7pt\hbox{.}\,\,\raise15pt\hbox{.}\hss}
&\hbox to 0pt{\hss\raise 3pt\hbox{.}\,.\lower4pt\hbox{.}\hss}
&&\righttag{\bullet}{n_{t-1}/q_{t-1}}{8pt}\dashto[ddll]\\
\\
&&\bullet\lineto[dr]&&\bullet\lineto[dl]\\
\lefttag{\bullet}{n_1/q_1}{8pt}\dashto[rr]&&
\bullet\lineto[r]&\overtag{\bullet}{b-t}{8pt}\lineto[r]&\bullet
\dashto[rr]&&\righttag{\bullet}{n_{t}/q_{t}}{8pt}\\&~\\&~}
$$ where now $q_{i}=n_{i}-p_{i}.$  It is not hard to see that $\partial
M_{\Gamma}$ and $\partial M_{\Gamma '}$ are orientation--reversing
diffeomorphic 3--manifolds; in fact, as noted in the next paragraph, putting on appropriate
complex structures, one sees that $M_{\Gamma }\cup M_{\Gamma '}$
is analytically isomorphic to a blow--up of some
$\mathbb P^1$-bundle over $\mathbb P^1$. It follows that for a
negative--definite graph $\Gamma$ its dual graph is no longer
negative--definite, but its intersection form is of type $(1,r)$.
While in Section~\ref{s:kirby} we embedded the manifold $\mg$ into $\#
_{\vert \Gamma \vert }\cpkk$, this time we will find an embedding of
$M_{\Gamma '}$ corresponding to the dual graph $\Gamma '$ into $\cpk \#
(\vert \Gamma ' \vert -1)\cpkk$, the $(\vert \Gamma '\vert -1)$--fold
blow--up of $\cpk$.  As will be discussed, a holomorphic embedding of
$M_{\Gamma '}$ with the above numerical condition will provide the
desired smoothing.  

Suppose $X=Spec(A)$ is a weighted homogeneous normal surface
singularity; thus, $A$ is positively graded, or $X$ has a good
$\C^{*}$-action, or $A$ is the quotient of a positively graded
polynomial ring by weighted homogeneous equations (see e.g. \cite{O-W}
for details as to what follows). Assume also that the link is a
rational homology sphere and suppose that the exceptional configuration of the minimal good resolution
$\tilde{X}\to X$ is given by the graph $\Gamma$.  A natural
$\C^{*}$-compactification of $X$ is given by $\bar{X}=Proj(A[u])$,
where $u$ has weight 1; what has been added is a curve with $t$ cyclic
quotient singularities.  Resolving just these singularities of $\bar{X}$ yields the projective variety $\bar{X'}$, where one has
added to $X$ a curve configuration determined by the dual graph $\Gamma '$.  (Further resolving on $\bar{X'}$ the original singularity thus yields two configurations of rational curves, corresponding to $\Gamma$ and to $\Gamma'$, connected by rational $-1$-curves, and blowing down to a $\mathbb P^1$-bundle over $\mathbb P^1$, as in \cite{OW}.)

A \emph{smoothing of negative weight} of $X$ (and of its compactifications $\bar{X}$ or $\bar{X'}$)
 is a smoothing obtained by adding terms
of lower weight to the defining equations, and gives rise to
deformations of both compactifications which are topologically locally
trivial at $\infty$ (see \cite{pinkdef} for details).  In particular,
such a smoothing of $\bar{X'}$ gives a smooth projective surface $Z$ with a curve
$D$ of type $\Gamma '$, and one has that the Milnor fiber of the
smoothing is diffeomorphic to $Z-D$ (e.g., \cite[Theorem 2.2]{Wahl}).
But under suitable cohomological conditions, there is a converse
\cite[Theorem 6.7]{pinkdef} due to Pinkham: an appropriate pair
$(Z,D)$ gives rise to a smoothing of negative weight.  We give a
special case of Pinkham's result:

\begin{thm}\label{pk}  Let $Z$ be a smooth projective rational surface, and
$D\subset Z$ a union of smooth rational curves whose intersection 
dual graph is 
of type $\Gamma '$.  Assume 
$$\text{rk}\ H_{2}(D; \bfz )=\text{rk}\ H_{2}(Z; \bfz ).$$ If $\Gamma$
is the graph of a rational singularity, then one has a $\mu=0$
smoothing of a rational weighted homogeneous singularity with
resolution dual graph $\Gamma$, and the interior of the Milnor fiber
is diffeomorphic to $Z-D$.
\end{thm}

\begin{proof}
  By \cite[Theorem~6.7]{pinkdef} one must check the vanishing of
  $H^1(Z,\mathcal O (E^{(k)}))$ for all $k \geq 0$, where the $E^{(k)}$ are
  effective divisors supported on what we have called $D$, and $E=E^{(1)}$ is
  the central curve.  The case $k=0$ is just the vanishing of the first cohomology of $Z$.  We follow Pinkham's notation and arguments closely.
  The key point is that the traces $D^{(k)}$ of the divisors $E^{(k)}$ on $E$
  are the familiar divisors used to write down the graded pieces of a weighted
  homogeneous singularity whose resolution dual graph is $\Gamma$; the
  vanishing of their first cohomology on $E$ is equivalent to the rationality
  of the singularity \cite[Theorem~5.7]{pinkham}.  Thus, by the exact
  sequences
$$0\rightarrow \mathcal O (E^{(k)}-E)\rightarrow \mathcal O (E^{(k)})
\rightarrow \mathcal O_E(D^{(k)})\rightarrow 0,$$
one needs only the vanishing of $H^1(Z,\mathcal O (E^{(k)}-E))$ for all $k
\geq 1$.  Write $F_k=E^{(k)}-E^{(k-1)}-E$, for $k\geq 1$; these are reduced
but reducible effective divisors, supported on the chains of rational curves
emanating from $E$.  Via the exact sequences
 $$0\rightarrow \mathcal O (E^{(k-1)})\rightarrow \mathcal O (E^{(k)}-E)
 \rightarrow \mathcal O_{F_k}(E^{(k)}-E)\rightarrow 0,$$ one can proceed
 inductively from the case $k=1$, reducing the claimed vanishing to showing
 that $H^1(F_k, \mathcal O_{F_k}(E^{(k)}-E))=0$ for all $k\geq 1$.  A proof is
 indicated in \cite[Lemma~6.9]{pinkdef}; or, one can argue directly with the
 resolutions of the cyclic quotient singularities whose resolution graphs
 contain the support of the $F_k$.
\end{proof}

\begin{rem} {\rm Note that the hypothesis $\text{rk}\ H_{2}(D ; \bfz
    )=\text{rk}\ H_{2}(Z ; \bfz )$, without the assumption on the rationality
    of the graph $\Gamma$, gives directly that $Z-D$ is a rational homology
    disk whose boundary is the link associated to $\Gamma $; one does not need
    anything about singularities.}
\end{rem}

We will describe pairs $(Z,D)$ by blowing up appropriately certain curves
$C\subset \cpk$, then taking $D$ a subset of the total transform of $C$.
Infinite families will come from systematically continuing the blow--up
procedure.

\begin{exa}{\bf [Graphs in the family $\farm$]}
{\rm 
Let $C=C_{1}\cup C_{2}\cup L_{1} \cup L_{2}$ be a 
plane curve of degree 6, where 
\begin{enumerate} 
    \item $C_{1}$ and $C_{2}$ are smooth conics 
with a triple tangency at one point 
\item $L_{1}$ is the line joining the two intersection points.
\item $L_{2}$ is the tangent line to $C_{1}$ at the simple 
intersection point with $C_{2}$.
\end{enumerate}
   Blowing-up appropriately 8 times yields a smooth surface 
   $Z$, so that the total transform of $C$ has the dual configuration:
\bigskip
$$
\xymatrix@R=6pt@C=24pt@M=0pt@W=0pt@H=0pt{
\\
&&&\overtag{\bullet}{C_{1}}{8pt}\lineto[r]&\lineto[r]
&\lineto[r]&\overtag{\Circ}{e_{6}}{8pt}\\
&&&\lineto[u]&&&\lineto[u]\\
&&&\lineto[u]&&&\lineto[u]\\
&&&\lineto[u]&&&\lineto[u]\\
&&&\lineto[u]&&&\lineto[u]\\
&\lefttag{\bullet}{e_{1}}{6pt}\lineto[r]&
\undertag{\bullet}{e_{3}}{6pt}\lineto[r]
&\undertag{\Circ}{e_{5}}{10pt}\lineto[r]\lineto[u]
&\overtag{\bullet}{C_{2}}{6pt}\lineto[r]
&\overtag{\bullet}{e_{2}}{6pt}\lineto[r]
&\undertag{\bullet}{e_{4}}{6pt}\lineto[u]\lineto[r]
&\undertag{\bullet}{L_{2}}{6pt}\lineto[r]
&\righttag{\bullet}{e_{7}}{6pt}\\
&\lineto[u]&&&\lineto[u]&\lineto[u]&&&\lineto[u]\\
&\lineto[u]&&&\lineto[u]&\lineto[u]&&&\lineto[u]\\
&\lineto[u]&&&\lineto[u]&\lineto[u]&&&\lineto[u]\\
&\lineto[u]&&&\lineto[u]&\lineto[u]&&&\lineto[u]\\
&\lineto[u]&&&\lineto[u]&\lineto[u]&&&\lineto[u]\\
&\lineto[u]\lineto[r]&\lineto[r]\lineto[r]
&\lineto[r]&\lineto[r]&\righttag{\Circ}{L_{1}}{6pt}\lineto[u]&&&\lineto[u]\\
&&&&&&&&\lineto[u]\\
&&&&\lineto[u]&&&&\lineto[u]\\
&&&&\lineto[u]&&&&\lineto[u]\\
&&&&\lefttag{\Circ}{e_{8}}{6pt}\lineto[u]\lineto[r]&\lineto[r]&\lineto[r]&\lineto[r]
&\lineto[u]}
$$
\bigskip

Here, solid unweighted vertices $\bullet$ denote $(-2)$--curves, while circle
vertices are $(-1)$'s.  Proper transforms of the components of $C$ are labeled
as before.  One recovers the blow--up procedure by blowing down sequentially
first $e_{8}$, then $e_{7}$, etc.  Let $\Gamma '$ be the configuration
obtained by deleting the vertices $e_{6}, L_{1}$, and $e_{8}$; it has rank 9
($={\mbox{rk }} H_2(Z; \bfz)$), and is dual to the basic configuration

$$
\xymatrix@R=6pt@C=24pt@M=0pt@W=0pt@H=0pt{
\\
&&\overtag{\bullet}{-2}{8pt}\\
&&\lineto[u]\\
&\undertag{\bullet}{-3}{8pt}\lineto[r]&
\undertag{\bullet}{-2}{8pt}\lineto[u]\lineto[r]&\undertag{\bullet}{-6}{8pt}}
$$
\bigskip

Thus, this construction gives a rational homology disk bounding the 
corresponding $\Sigma$.  
But one can repeatedly blow--up further, between $e_{2}$ and 
  $L_{1}$, and then between the transform of $e_{2}$ and its new 
  neighbour.  The same idea works between $e_{4}$ and $e_{6}$ (and 
  again between the transform of $e_{4}$ and the new exceptional 
  curves), as well as between $C_{2}$ and $e_{8}$.  Doing this 
  respectively $p, q,$ and $r$ times, the new exceptional diagram is
$$
\xymatrix@R=6pt@C=24pt@M=0pt@W=0pt@H=0pt{
\\
&&&\overtag{\bullet}{C_{1}}{8pt}\lineto[r]&\overtag{\bullet}{e_{6}}{6pt}\lineto[r]
&\overtag{\bullet}{g_{1}}{6pt}\dashto[r]&\dashto[r]&\overtag{\bullet}{g_{q-1}}{6pt}\lineto[r]&\overtag{\Circ}{g_{q}}{8pt}\\
&&&\lineto[u]&&&&&\lineto[u]\\
&&&\lineto[u]&&&&&\lineto[u]\\
&&&\lineto[u]&&&&&\lineto[u]\\
&&&\lineto[u]&&&&&\lineto[u]\\
&\lefttag{\bullet}{e_{1}}{6pt}\lineto[r]&
\undertag{\bullet}{e_{3}}{6pt}\lineto[r]
&\undertag{\Circ}{e_{5}}{10pt}\lineto[r]\lineto[u]
&\overtag{\bullet}{-(r+2)}{8pt}\lineto[r]_(.3){C_2}&\lineto[r]_(.7){e_2}
&\overtag{\bullet}{-(p+2)}{8pt}\lineto[r]&\lineto[r]^(.5){-(q+2)}
&\undertag{\bullet}{e_{4}}{6pt}\lineto[u]\lineto[r]
&\undertag{\bullet}{L_{2}}{6pt}\lineto[r]
&\righttag{\bullet}{e_{7}}{6pt}\\
&\lineto[u]&&&\lineto[u]&&\lineto[u]&&&&\lineto[u]\\
&\lineto[u]&&&\lineto[u]&&\lineto[u]&&&&\lineto[u]\\
&\lineto[u]&&&\lineto[u]&&\lineto[u]&&&&\lineto[u]\\
&\lineto[u]&&&\lineto[u]&&\lineto[u]&&&&\lineto[u]\\
&\lineto[u]&&&\lineto[u]&&\lineto[u]&&&&\lineto[u]\\
&\lineto[u]&&&\lineto[u]&&\lineto[u]&&&&\lineto[u]\\
&\lineto[u]&&&\lineto[u]&&\lineto[u]&&&&\lineto[u]\\
&\lefttag{\bullet}{L_{1}}{6pt}\lineto[u]\lineto[r]
&\undertag{\bullet}{f_{1}}{6pt}\dashto[r]
&\dashto[r]&\dashto[r]&\undertag{\bullet}{f_{p-1}}{2pt}\lineto[r]
&\undertag{\Circ}{f_{p}}{6pt}\lineto[u]&&&&\lineto[u]\\
&&&&&&&&&&\lineto[u]\\
&&&&\lineto[u]&&&&&&\lineto[u]\\
&&&&\lineto[u]&&&&&&\lineto[u]\\
&&&&\lineto[u]&&&&&&\lineto[u]\\
&&&&\lineto[u]&&&&&&\lineto[u]\\
&&&&\undertag{\Circ}{h_{r}}{9pt}\lineto[u]\lineto[r]
&\undertag{\bullet}{h_{r-1}}{4pt}\dashto[r]&\dashto[r]
&\dashto[r]&\dashto[r]&\undertag{\bullet}{h_{1}}{4pt}\lineto[r]
&\undertag{\bullet}{e_{8}}{6pt}\lineto[u]}
$$ 

\bigskip

  To reverse the process,  blow down sequentially the 
  $f_{i},\ g_{i},\ h_{i}$ starting from the largest subscript (and 
  define $f_{0}=L_{1}, g_{0}=e_{6}, h_{0}=e_{8}$).  The graph $\Gamma '$ 
   obtained by deleting the $-1$ curves $f_{p},\  g_{q}$, and $h_{r}$ is
   $$
\xymatrix@R=6pt@C=24pt@M=0pt@W=0pt@H=0pt{
\\&&&\bullet\dotto[d]\\
&&&\dotto[u]\\
&&&\dotto[u]&\vbox to 0 pt{\vss\hbox to 0 pt{\hss$\left.\vbox to 20 pt{}\right\}$\hss}\vss}\righttag{}{(q+1)}{11pt}\\
&&&\dotto[u]\\
&(p+2)&&\bullet\dotto[u]&&&&&(r+2)\\
&{\hbox to 0pt{\hss$\overbrace{\hbox to 90pt{}}$\hss}}&&\lineto[u]&&&&&{\hbox to 0pt{\hss$\overbrace{\hbox to 90pt{}}$\hss}}\\
\bullet\dashto[r]&\dashto[r]&\bullet\lineto[r]&\undertag{\bullet}{-1}{6pt}\lineto[r]\lineto[u]&\undertag{\bullet}{-(r+2)}{6pt}\lineto[r]&\undertag{\bullet}{-(p+2)}{6pt}\lineto[r]&\undertag{\bullet}{-(q+2)}{6pt}\lineto[r]&\bullet\dashto[r]&\dashto[r]&\bullet}
$$
\bigskip
   
It is easy to calculate that the dual graph $\Gamma$ is of type
$\Lambda_{q,r,p}$.  (We leave it to the reader to check the special cases when
some of the $p,q,r$ are 0.)  This gives the desired boundaries for all graphs
of type $\mathcal M$.}
\end{exa}

\begin{exa}{\bf [Graphs in the families $\frw $ and $\frn$]}
{\rm A similar construction, starting with $C$ the union of 4 lines
$L_1, L_2,L_3,L_4$ in general position in $\cpk$, yields rational
homology ball smoothings for the triply--infinite families of type
$\mathcal W$ and $\mathcal N$.  For the (3,3,3) type $\mathcal W$,
blow--up 6 times, giving
$$
\xymatrix@R=6pt@C=24pt@M=0pt@W=0pt@H=0pt{
\\
&&\overtag{\Circ}{e_4}{8pt}\lineto[r]
&\overtag{\bullet}{e_{1}}{8pt}\\
&&\lineto[u]&\lineto[u]\\
&&\lineto[u]&\lineto[u]\\
&&\lineto[u]&\lineto[u]\lefttag{\bullet}{L_3}{4pt}\lineto[r]&\lineto[r]&\righttag{\Circ}{e_5}{8pt}\\
&&\lineto[u]&\lineto[u]&&\lineto[u]\\
&&\lineto[u]&\lineto[u]&&\lineto[u]\\
&\lefttag{\bullet}{e_{2}}{6pt}\lineto[r]
&\undertag{\bullet}{L_2}{3pt}\lineto[r]\lineto[u]
&\undertag{\bullet}{L_{1}}{3pt}\lineto[r]^(.2){+1}\lineto[u]
&\overtag{\bullet}{L_4}{9pt}\lineto[r]
&\righttag{\bullet}{e_3}{8pt}\lineto[u]\\
&\lineto[u]&&&\lineto[u]\\
&\lineto[u]&&&\lineto[u]\\
&\lineto[u]&&&\lineto[u]\\
&\lineto[u]\lineto[r]&\lineto[r]&\lineto[r]&\undertag{\Circ}{e_6}{8pt}\lineto[u]}
$$

\bigskip

Blowing up $p$ additional times between $L_4$ and $e_6$ and successor
$(-1)$--curves, $q$ times between $L_3$ and $e_5$, and $r$ times between
$L_2$ and $e_4$, and then removing the three $(-1)$--curves, gives a
$\Gamma '$ which is easily checked to be dual to $\Gamma_{p,q,r}$.

For the (2,4,4) type $\mathcal N$, blow up 7 times, yielding 
$$
\xymatrix@R=6pt@C=24pt@M=0pt@W=0pt@H=0pt{
\\
&&\overtag{\Circ}{e_4}{8pt}\lineto[r]
&\lineto[r]&\overtag{\bullet}{e_{1}}{8pt}\\
&&\lineto[u]&&\lineto[u]\\
&&\lineto[u]&&\lineto[u]\\
&&\lineto[u]&&\lineto[u]\\
&&\lineto[u]&&\lineto[u]\\
&\lefttag{\bullet}{e_{3}}{6pt}\lineto[r]
&\undertag{\bullet}{L_1}{6pt}\lineto[r]\lineto[u]
&\overtag{\bullet}{L_{2}}{8pt}\lineto[r]
&\undertag{\bullet}{L_3}{6pt}\lineto[r]^(.1){0}\lineto[u]
&\overtag{\bullet}{L_{4}}{8pt}\lineto[r]
&\overtag{\bullet}{e_{2}}{6pt}\lineto[r]
&\righttag{\bullet}{e_{6}}{8pt}\\
&\lineto[u]&&\lineto[u]&&\lineto[u]&&\lineto[u]\\
&\lineto[u]&&\lineto[u]&&\lineto[u]&&\lineto[u]\\
&\lineto[u]&&\lineto[u]&&\lineto[u]&&\lineto[u]\\
&\lineto[u]&&\lineto[u]&&\lineto[u]&&\lineto[u]\\
&\lineto[u]&&\lineto[u]&&\lineto[u]&&\lineto[u]\\
&\lineto[u]\lineto[r]&\lineto[r]
&\lineto[r]&\lineto[r]&\righttag{\Circ}{e_{5}}{6pt}\lineto[u]&&\lineto[u]\\
&&&&&&&\lineto[u]\\
&&&\lineto[u]&&&&\lineto[u]\\
&&&\lineto[u]&&&&\lineto[u]\\
&&&\lefttag{\Circ}{e_{7}}{6pt}\lineto[u]\lineto[r]&\lineto[r]&\lineto[r]&\lineto[r]
&\lineto[u]}
$$
\bigskip

(The self-intersection of $L_3$ is 0.)  Blowing up $p$ additional
times between $L_1$ and $e_4$ and successor $(-1)$--curves, $q$ times
between $L_2$ and $e_7$, and $r$ times between $L_4$ and $e_5$, and
then removing the three $(-1)$--curves, gives $\Gamma '$ which is easily
checked to be dual to $\Delta_{p,q,r}$.}
\end{exa}

\begin{rem}  
{\rm How would one come up with such a configuration, obtainable by blowing
up curves in $\cpk$, given a star--shaped candidate $\Gamma$?
One first forms $\Gamma '$ as above.  The lattice generated by the
vertices is supposed to have a finite overlattice which is the second
homology of a surface $Z$; this must be the unimodular odd lattice
with signature $(1,n-1)$, where $n$ is the number of curves in
$\Gamma '$.  In particular, generators of the overlattice are certain
rational combinations of the vertices in $\Gamma '$, and the class $K$
of a canonical divisor is also such a combination.  To find suitable
rational $(-1)$--curves on the putative $Z$, one looks for elements $e$ of
the overlattice for which $e\cdot e=-1, K\cdot e=-1$, and which
intersect exactly two curves in $\Gamma '$ (or intersects one curve twice).  One might be able to use
some of these to produce a configuration which blows down.  This
approach has succeeded in a number of cases, including all the ones in
the above examples.  We present two more, of type $\mathcal C$, though
there are quite a few others.}
\end{rem}

\begin{exa}{\bf [A family in $\frc$]}
{\rm
Next we give an infinite family of graphs in
$\frc$ which arise from normal surface singularities which admit rational
homology disk smoothings. Let $C=C_1\cup C_2$ be a plane curve of degree 5,
where $C_1$ is an irreducible nodal cubic, and $C_2$ is a smooth quadric
intersecting $C_1$ in only one point (hence, with intersection multiplicity
6).  An affine version of such a curve is $$(y^2-x(1-x)^2)(y^2-x(1-2x))=0.$$
Blowing up appropriately 8 times yields a smooth surface for which the total
transform of $C$ has dual configuration
$$
\xymatrix@R=6pt@C=24pt@M=0pt@W=0pt@H=0pt{
\\
\lineto[r]\lineto[d]&\overtag{\Circ}{e_3}{6pt}\lineto[d]&\overtag{\bullet}{C_2}{8pt}\lineto[d]\\
\lineto[d]&\lineto[d]&\lineto[d]\\
\lineto[d]&\lineto[d]&\lineto[d]\\
\undertag{\bullet}{e_1}{8pt}\lineto[r]&\undertag{\bullet}{C_1}{6pt}\lineto[r]&\undertag{\Circ}{e_8}{11pt}\lineto[r]&\undertag{\bullet}{e_7}{8pt}\lineto[r]&\undertag{\bullet}{e_6}{8pt}\lineto[r]&\undertag{\bullet}{e_5}{8pt}\lineto[r]&\undertag{\bullet}{e_4}{8pt}\lineto[r]&\undertag{\bullet}{e_2}{8pt}\\
&&&&&\\
&&&&&\\
&&&&&\\
}
$$

Blowing up $p$ additional times between $C_1$ and $e_3$ and its successor
$(-1)$'s, and then removing the last $(-1)$--curve, gives $\Gamma '$ 
$$
\xymatrix@R=6pt@C=24pt@M=0pt@W=0pt@H=0pt{
\\
&(p+1)&&&\bullet&&&&&\\
&{\hbox to 0pt{\hss$\overbrace{\hbox to 90pt{}}$\hss}}&&&\lineto[u]&&&&&\\
\bullet\dashto[r]&\dashto[r]&\bullet\lineto[r]&\undertag{\bullet}{-(p+2)}{6pt}\lineto[r]&\undertag{\bullet}{-1}{6pt}\lineto[r]\lineto[u]&\bullet\lineto[r]&\bullet\lineto[r]&\bullet\lineto[r]&\bullet\lineto[r]&\bullet\\
&&&&&\\
&&&&&\\
&&&&&\\
}
$$
The dual graph $\Gamma$ in this case is of type $\mathcal C$, and is
$$
\xymatrix@R=6pt@C=24pt@M=0pt@W=0pt@H=0pt{
\\
&&&&&p&&\righttag{\bullet}{-2}{6pt}&&&&&\\
&&&&&{\hbox to 0pt{\hss$\overbrace{\hbox to 60pt{}}$\hss}}&&\lineto[u]&&&&&\\
&&&\undertag{\bullet}{-(p+3)}{6pt}\lineto[r]&\undertag{\bullet}{-2}{6pt}\dashto[r]&\dashto[r]&\undertag{\bullet}{-2}{6pt}\lineto[r]&\undertag{\bullet}{-2}{6pt}\lineto[r]\lineto[u]&\undertag{\bullet}{-6}{6pt}\\
&&&&&\\
&&&&&\\
&&&&&\\
}
$$
Thus, this graph corresponds to a singularity with a $\mu=0$ smoothing.
}
\end{exa}

\begin{exa}\label{e:newfourv}
{\bf [A family in $\frc$ with a 4--valent node]} {\rm We finish
    this subsection by constructing an infinite family of graphs of type
    $\frc$ with one node of valency 4, which arise from a rational homology
    disk smoothing of negative weight.  We thank E. Shustin for assistance with
    this example.

We consider four conics in the plane,
    all contained in the real unit circle:
\[
F= \qquad \qquad \ \ \ \ \ \ \ \ \ \ \ \ \ \ \ \ \ \ \  \{ x^2+y^2-1=0\}
\]
\[
G= \quad \{ x^2+y^2-1+(1/2)(x+1)(y+1)=0 \}
\]
\[
H= \qquad \ \ \  \ \  \{ x^2+y^2-1-(3/4)(x^2-1)=0 \}
\]
\[
J= \quad \{ x^2+y^2-1-(1/2)(x-1)(y+1)=0\} .
\]
Relevant intersection points are $(1,0)$, where $F,H, J$ are pairwise tangent;
$(-1,0)$, where $F, G,H$ are pairwise tangent; $(0,1/2)$, through which
$G,H,J$ pass transversally; and $(0,-1)$, where $G$ and $J$ have a triple
tangency, and $F$ has a simple tangency with each of them.  Let $L$ be the
common tangent line $y=-1$.  There are two further important points: $H$ and
$J$ also meet transversally at $(-3/5,-2/5)$, which lies on the straight line
$M$ through $(-1,0)$ and $(0,-1)$; and finally $H$ and $G$ meet transversally
at $(3/5,-2/5)$, which lies on the straight line $N$ through $(1,0)$ and
$(0,-1)$.

We now consider the (complex) plane curve which is the union of the 4 conics $F,G,H,J$
and three lines $L,M,N$, and blow up minimally so that the reduced total
transform has only normal crossings.  We finally blow up either of the two
intersection points of the line $L$ and the conic $H$ (which are in the
imaginary complex plane).  We reach the diagram of Figure~\ref{f:compl}, where
the proper transforms of the seven curves in the plane are given the same
name.
\begin{figure}[ht]
\begin{center}
\includegraphics[width=10cm]{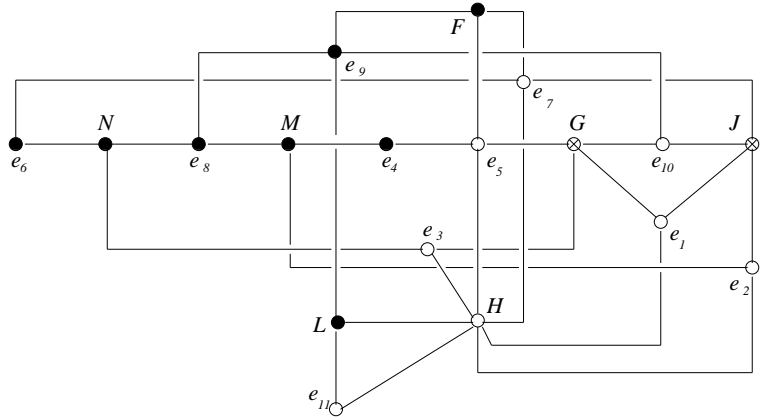}
\end{center}
\caption{The configuration of curves after 11 blow--ups; solid circles
  represent $(-2)$--curves, the empty ones $(-1)$--curves, the
ones with a cross correspond to $(-3)$--curves, while $H\cdot H=-4$}
\label{f:compl}
\end{figure}

Again, blowing down the curves $e_i$ starting from the largest subscript will
show how to reverse the procedure.  Now, removing the vertices
$e_1,e_2,e_3,e_7$ and $e_9$ from the diagram, and then blowing down the
$(-1)$--curve $e_{10}$, yields the simpler diagram of
Figure~\ref{f:fourvalent}.

\begin{figure}[ht]
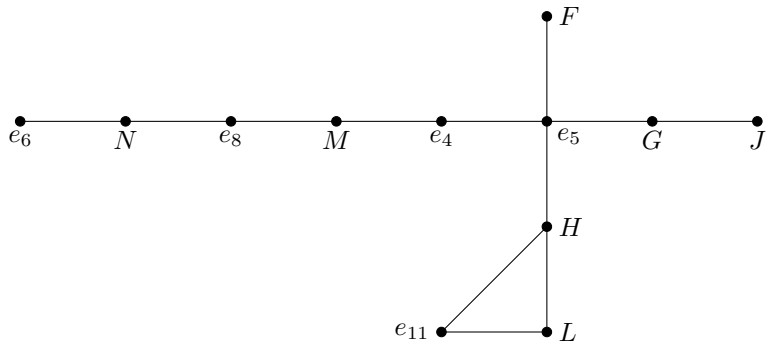

\begin{center}
\setlength{\unitlength}{1mm}
\unitlength=0.7cm
\begin{graph}(-6,7)(9,-5)
\graphnodesize{0.2}

 \roundnode{m1}(-2,0)
 \roundnode{m2}(0,0)
 \roundnode{m3}(2,0)
 \roundnode{m4}(4,0)
 \roundnode{m5}(6,0)
 \roundnode{m6}(8,0)
 \roundnode{m7}(10,0)
 \roundnode{m8}(12,0)
 \roundnode{m9}(8,2)
 \roundnode{m10}(8,-2)
 \roundnode{m11}(8,-4)
 \roundnode{m12}(6,-4)

\edge{m1}{m2}
\edge{m2}{m3}
\edge{m3}{m4}
\edge{m4}{m5}
\edge{m5}{m6}
\edge{m7}{m6}
\edge{m7}{m8}
\edge{m9}{m6}
\edge{m10}{m6}
\edge{m10}{m11}
\edge{m12}{m11}
\edge{m10}{m12}

  \autonodetext{m1}[s]{{\small $e_6$}}
  \autonodetext{m2}[s]{{\small $N$}}
  \autonodetext{m3}[s]{{\small $e_8$}}
  \autonodetext{m4}[s]{{\small $M$}}
  \autonodetext{m5}[s]{{\small $e_4$}}
  \autonodetext{m6}[se]{{\small $e_5$}}
  \autonodetext{m7}[s]{{\small $G$}}
  \autonodetext{m8}[s]{{\small $J$}}
  \autonodetext{m9}[e]{{\small $F$}}
  \autonodetext{m10}[e]{{\small $H$}}
  \autonodetext{m11}[e]{{\small $L$}}
  \autonodetext{m12}[w]{{\small $e_{11}$}}

\end{graph}
\end{center}
\caption{The configuration of Figure~\ref{f:compl} after blow--ups ($e_5$ and 
$e_{11}$ are $(-1)$--curves while $H$ is a $(-4)$--curve)}
\label{f:fourvalent}
\end{figure}

It is clear that removing $e_{11}$ gives a curve configuration whose
components span the second rational homology group of the blown--up plane
(recall that $e_{10}$ has been blown down, so after deleting $e_{11}$, the
remaining eleven curves are in $\cpk \# 10\cpkk$). Furthermore, the resulting
$\Gamma'$ is dual to the graph of a rational singularity, hence one
can apply Theorem~\ref{pk}.  Blowing up further between $H$ and $e_{11}$, and
then removing the exceptional curve yield rational homology disk smoothings of
negative weight of singularities with resolution graph

\bigskip

$$
\xymatrix@R=6pt@C=24pt@M=0pt@W=0pt@H=0pt{
\\&&&&p+1\\
&&\overtag{\bt}{-3}{8pt}
&&{\hbox to 0pt{\hss$\overbrace{\hbox to 60pt{}}$\hss}}&\\
&&\lineto[u]\\
&\undertag{\bt}{-6}{6pt}\lineto[r]&
\bt\lineto[u]\lineto[r]_(.2){-3}&\undertag{\bt}{-2}{6pt}\dashto[r]&\dashto[r]
&\undertag{\bt}{-2}{6pt}\lineto[r]
&\undertag{\bt}{-(p+2)}{6pt}\\
&&\lineto[u]\\
&&\lineto[u]\\
&&\undertag{\bt}{-2}{6pt}\lineto[u]&&&}
$$

\bigskip

(In the case $p=0$ one has to blow down $e_{11}$ and remove the curve $L$.)
}
\end{exa}  

\subsection{Examples of $\mu=0$ smoothings via quotients}

Finally we show a second algebro--geometric way for constructing smoothings
with $\mu =0$. Besides discussing the (already familiar) cases we have dealt
with in the previous subsection, there will be some families of singularities
(with resolution graphs in $\fra\cup \frb$) which can be treated in this way.
We will show further examples of graphs involving nodes with valency 4
giving rise to singularities with rational homology disk smoothing.  Let us
start by describing the general strategy first.

Suppose $(\mathcal Y,o)$ is the germ of an isolated 3-dimensional
normal singularity, $G$ a finite group of automorphisms acting freely
off $o$, and $f\in \mathcal O_{\mathcal Y ,o}$ a $G$--invariant
function whose zero locus $(Y,o)$ has an isolated normal singularity.  (This implies that $\mathcal Y$ is Cohen--Macaulay.)
Then $f$ defines a smoothing of $(Y,o)$, with Milnor fiber $M_{Y}$;
but it also defines a smoothing of $X= Y/G$, whose Milnor fiber
$M_{X}$ is the quotient of the free action of $G$ on $M_{Y}$.  In
particular $$\chi (M_{Y})=|G|\cdot \chi (M_{X}).$$ Consequently, if
$|G|=\chi (M_{Y})=1+\mu _Y,$ then $\chi (M_X)=1$ and so we have constructed a
smoothing of $X$ with Milnor number $\mu _X= 0$.

In the following, we will denote by $(1/m)[a_{1},\cdots,a_{n}]$ the
transformation of $\mathbb C ^{n}$ which is multiplication of the
$j$th coordinate by exp$(2\pi ia_{j}/m)$ (where $a_{j}\in \mathbb Z$).
All examples below are weighted homogeneous, simplifying calculations.
It is easy to check that a given group action is free on $\mathcal Y
-\{o\}$.  Describing the resolution graph of $X$ is achieved by
lifting the action of $G$ to the weighted blow--up of $Y$, then
dividing and resolving the quotient, cf.  \cite[Section~8]{nw1}.  One
needs to compute the Milnor number $\mu _Y$ for $Y$ in order to show
that it satisfies $1+\mu _Y=|G|$.  But, in our examples $Y$ is a
complete intersection defined by one or two weighted homogeneous
equations, so one can easily apply \cite[Korollar~3.10(a)~or~(c)]{gh}.

\begin{exa}{\bf [Graphs in $\frg$]}
{\rm
Let us start with the case of the standard rational blow--down along $L(p^2,
pq-1)$ (cf. also \cite[(5.9.1)]{Wahl}). Suppose $0<q<p,\ (q,p)=1$.  Let
$(\mathcal Y,o)=(\mathbb C ^{3},o)$, $G$ the cyclic group generated by
$(1/p)[1,q,-1]$, and $f=xz-y^{p}$.  Then one gets $\mu =0$ smoothings of $X$,
the cyclic quotient singularity of type $p^{2}/(pq-1)$, whose resolution graph
is a chain as in Definition 1.3.  The Milnor fiber is a free quotient of the
(simply connected) Milnor fiber of the $A_{p-1}$-singularity $xz-y^{p}=0$.
}
\end{exa}

The three basic graphs of Figure~\ref{f:n=4} (Section~\ref{s:caseb}),
corresponding to the spherical triples $(3,3,3)$, $(2,4,4)$,
$(2,3,6)$, are resolution graphs of certain \emph{log--canonical}
surface singularities; each is the quotient of a \emph{simple
elliptic} singularity (i.e., cone over an elliptic curve) by a cyclic
group acting freely off the singular point.  As mentioned in the
Introduction, each gives rise to a triply--infinite family of
singularities having a smoothing with Milnor number 0, and these
examples fill out the classes $\frw, \frn$, and $\farm$.

\begin{exa}{\bf [Graphs in $\frw$]}
  {\rm We apply the above principle to graphs of type $\frw$; the following is
    exactly \cite[Example~(5.9.2)]{Wahl}. Suppose $p,q,r \geq 0$, and let
    $N=(p+2)(q+2)(r+2)+1.$ Let $(\mathcal Y,o)=(\mathbb C ^{3},o)$, $G$ the
    cyclic group generated by $(1/N)[1,(q+2)(r+2),-(r+2)]$, and
    $f=xy^{p+2}+yz^{q+2}+zx^{r+2}$.  Then one gets $\mu =0$ smoothings of
    rational singularities whose resolution graphs we claim are the
    $\Gamma_{p,q,r}$ as in Figure~\ref{f:wahltype}. The Milnor fiber has
    fundamental group $G$.
    
    We indicate how to find the quotient graph $\Gamma$.  Give the variables
    $x,y,z$ weights $a,b,c$, respectively, where $a=(p+2)(q+1)+1,
    b=(q+2)(r+1)+1,\ c=(r+2)(p+1)+1$; this makes $f$ weighted homogeneous, of
    weight $N$.  Since $G$ preserves weights, the quotient is still weighted
    homogeneous; and it is rational, with discriminant group of order $N^2$
    (by Proposition \ref{rat}).  We show that the graph $\Gamma$ of the
    quotient has three chains of rational curves, corresponding to the
    continued fractions (from the outside to the center) $a/(q+1),\ b/(r+1),\
    c/(p+1)$.  This will give $\Gamma_{p,q,r}$, except that the central
    self--intersection of $-4$ is deduced at the end from the calculation of
    the discriminant order.
    
    Use the weights $a,b,c$ to do a weighted blow--up of $\mathbb C^3$, as in
    \cite[Section~4]{nw1}. This space is covered by three open sets.  The
    first is the quotient of a copy of $\mathbb C^3$, with coordinates
    $u,v,w$, modulo the cyclic action of $S=(1/a)[-1,b,c]$, and the
    coordinates are related by $x=u^a,y=u^bv,z=u^cw.$ The proper transform of
    $f=0$ is checked to be the smooth surface given by $v^{p+2}+vw^{q+2}+w=0.$
    A lifting of the the generator of $G$ to this $\mathbb C^3$ is
    $T=(1/aN)[1,q+1,-1]$; one divides the smooth surface by the group
    generated by $T$ and $S$.  But $T^N$ is computed to be the inverse of $S$,
    so one only needs to divide by $T$.  Further, $T^a=(1/N)[1,0,0]$ is a
    pseudo--reflection; so one divides out first by this transformation,
    replacing the variable $u$ by $u'=u^N$.  Thus, in the $u',v,w$ space, one
    has the smooth surface $v^{p+2}+vw^{q+2}+w=0,$ and a cyclic group action
    generated by $(1/a)[1,q+1,-1]$.  The only fixed point occurs at
    $u'=v=w=0$; at such a point, $u'$ and $v$ are local analytic coordinates,
    and one divides out by the group $(1/a)[1,q+1]$.  The image of $u'=0$ is
    the central rational curve, so the continued fraction expansion of
    $a/(q+1)$ goes from the outer curve into the central curve.  This is the
    cyclic quotient singularity claimed above.  }
\end{exa}

\begin{exa}{\bf [Graphs in $\frn$]} {\rm Suppose now that $p,q,r \geq 0$, and
    let us define $N=(p+1)(q+3)(r+2)+q+2$.  Let $(\mathcal Y,o)\subset(\mathbb
    C ^{4},o)$ be the hypersurface singularity given by the equation
    $(x^{q+2}+z^{r+2}-yw=0)$, $G$ the cyclic group generated by
    $(1/N)[1,(q+3), -(p+1)(q+3), -1]$, and $f=xw-y^{p+1}z.$ Then one gets
    $\mu=0$ smoothings of rational singularities with resolution graph
    $\Delta_{p,q,r}$ as in Figures~\ref{f:masik} and \ref{f:specmasik}.  The
    Milnor fiber has fundamental group $G$.}
\end{exa}

\begin{exa}\label{ex:236} {\rm Embed the plane cubic curve
    $A^{3}+B^{3}-AC^{2}=0$ into ${\mathbb C} {\mathbb P}^{5}$ via the Veronese
    embedding $$[A,B,C]\mapsto [x_{0},x_{1},x_{2},y_{0},y_{1},y_{2}]\equiv
    [A^{2},B^{2},C^{2},BC,AC,AB],$$ and let $Y\subset \C^{6}$ be the affine
    cone.  Let $\mathcal Y$ be the total space of the one-parameter smoothing
    (with parameter $T$) of $Y$, given $Y$ is defined by 9 equations, as in
    \cite[(9.6)]{pink}, and a 1-parameter smoothing (with parameter $T$) is
    given by the vanishing of the nine expressions
    $$x_{0}x_{1}-y_{2}^{2}+Ty_{1}$$ $$x_{1}x_{2}-y_{0}^{2}+Ty_{1}$$
    $$x_{0}x_{2}-y_{1}^{2}+Ty_{0}$$ $$x_{0}y_{0}-y_{1}y_{2}+Tx_{1}$$
    $$x_{1}y_{1}-y_{0}y_{2}+Tx_{0}$$ $$x_{2}y_{2}-y_{0}y_{1}+T^{2}$$
    $$x_{0}^{2}+x_{1}y_{2}-y_{1}^{2}$$ $$x_{1}^{2}+x_{0}y_{2}-y_{0}y_{1}$$
    $$x_{1}y_{0}+x_{0}y_{1}-x_{2}y_{1}-Ty_{2}.$$ 
(This requires some checking,
    but easily follows as in \cite[(9.6)]{pink}.)  The cyclic group $G$ of
    order 6, acting on the $X_{i},Y_{j},T$ via $(1/6)[1,3,1,2,4,5,0]$, acts
    freely off the origin of $\mathcal Y$.  $T$ is $G$--invariant, and $Y$ mod
    $G$ is the by--now--familiar singularity with resolution graph
$$
\xymatrix@R=6pt@C=24pt@M=0pt@W=0pt@H=0pt{
\\
&&\overtag{\bt}{-3}{8pt}\\
&&\lineto[u]\\
&\undertag{\bt}{-2}{8pt}\lineto[r]&
\undertag{\bt}{-2}{8pt}\lineto[u]\lineto[r]&\undertag{\bt}{-6}{8pt}}
$$

\bigskip

\bigskip

The Milnor fiber of $Y$ is simply connected, with Euler characteristic 6; so,
one has a $\mu=0$ smoothing of $X$, whose Milnor fiber has fundamental group
$G$.}
\end{exa}

\begin{exa}\label{e:intex} {\bf [A family of graphs in $\fra$]} 
{\rm Let $p\geq 0$, and for convenience write $m=3p^{2}+9p+7, \
r=3p^{2}+6p+2=m-(3p+5),\ n=9(p+2)$ and  $a=3p+5$.  Consider the
\emph{metacyclic group} $G\subset GL(3,\C)$, acting freely on
$\C^{3}-\{0\}$, generated by $S=(1/m)[1,r,r^{2}]$ and $T$ defined by
$T(x,y,z)=(y,z,\omega x)$, where $\omega$=exp$(2\pi i/(a+1))$.  One
easily finds that $S^{m}=T^{n}=I,\ TST^{-1}=S^{r}$, and $T^{3}$ is
multiplication by $\omega.$ Further, $G/[G,G]$ is cyclic of order $n$,
generated by the image of $T$, and $|G|=mn.$ Let $\mathcal Y\
=\C^{3}$, and $f=x^{a}y+y^{a}z+\omega xz^{a}.$ Then $f$ is
$G$--invariant, and defines a hypersurface singularity $Y$ whose Milnor
fiber has Euler characteristic $1+a^{3}=mn=|G|$.  Thus, one gets a
$\mu=0$ smoothing of a rational singularity whose resolution graph has
valency 4:
\bigskip
$$
\xymatrix@R=6pt@C=24pt@M=0pt@W=0pt@H=0pt{
\\&&&&p\\
&&\overtag{\bt}{-3}{8pt}
&&{\hbox to 0pt{\hss$\overbrace{\hbox to 60pt{}}$\hss}}&\\
&&\lineto[u]\\
&\undertag{\bt}{-3}{6pt}\lineto[r]
&\bt\lineto[u]\lineto[r]_(.2){-3}&\undertag{\bt}{-2}{6pt}\dashto[r]
&\dashto[r]
&\undertag{\bt}{-2}{6pt}\lineto[r]
&\undertag{\bt}{-4}{6pt}\lineto[r]
&\undertag{\bt}{-(p+2)}{6pt}\\
&&\lineto[u]\\
&&\lineto[u]\\
&&\undertag{\bt}{-3}{6pt}\lineto[u]}
$$
\bigskip

To find the resolution graph, let $\bar{G}$ be the image of $G$ in
$PGL(3,\C)$ (modding out by $T^{3}$), and locate the orbits in the
projective plane curve defined by $f$ on which $\bar{G}$ acts with
non--trivial isotropy.  Then, resolve the hypersurface singularity via
one blow--up, divide out by $T^{3}$, and find and describe the fixed
points of $\bar{G}$ on this smooth surface.  One finds three orbits on
which every point has isotropy of order 3; and 3 points where the
isotropy is generated by $S$, and the image of $T$ puts them into one
$\bar{G}$--orbit.  Dividing out then gives all the information on the
resolution, except for the weight of the central curve; this can be
computed by noting that the discriminant group has order $n^{2}$.

Since the Milnor fiber of $Y$ is simply connected, the fundamental
group of the $\mu=0$ Milnor fiber is the non--abelian group $G$.
(This perhaps explains why this example, originally discovered in
1983, was surprising --- a Seifert manifold with 4 branches which bounds
a rational homology ball.)  Note that this example is of type
$\mathcal A$.}
\end{exa}
\begin{rem}
{\rm There are very few ``interesting'' subgroups of $GL(3,\C)$ which
act freely off the origin, see \cite{wolf}.}
\end{rem}
\begin{exa} {\bf [A family in $\frb$]}
\label{e:valency4}
{\rm Let $p\geq 2$, and for convenience write $m=2p^{2}-2p-1,\
r=2(p-1)^{2}, \ n=8p, \ a=2p-1$.  Consider the group $G\subset
GL(4,\C)$, generated by $S=(1/m)[1,r,-1,-r]$ and $T$ defined by
$T(x,y,z,w)=(\eta^{4} w,x,y,z),$ where $\eta$ = exp$(2\pi i/4n)$.  One
easily finds that $S^{m}=T^{n}=I,\ T^{-1}ST=S^{r}$, and $T^{4}$ is
multiplication by $\eta^{4}$.  Again, $G/[G,G]$ is cyclic of order
$n$, generated by the image of $T$, and $|G|=mn$.  Let $\mathcal Y\
=\{xz+\eta^{2}yw=0\}\subset \C^{4}$, on which $G$ can be checked to
act freely off the origin.  Then
$f=x^{a}y+y^{a}z+z^{a}w+\eta^{4a}w^{a}x$ is $G$--invariant, and defines
a complete intersection singularity $Y$ whose Milnor fiber has Euler
characteristic $mn=|G|$.  Thus, one gets a $\mu =0$ smoothing of a
rational surface singularity.  The same approach as above (first
finding fixed points of some $\bar{G}$ on the projective curve) yields
that the resolution graph of the singularity is
\bigskip
$$
\xymatrix@R=6pt@C=24pt@M=0pt@W=0pt@H=0pt{
\\&&&&p\\
&&\overtag{\bt}{-4}{8pt}
&&{\hbox to 0pt{\hss$\overbrace{\hbox to 60pt{}}$\hss}}&\\
&&\lineto[u]\\
&\undertag{\bt}{-4}{6pt}\lineto[r]&
\bt\lineto[u]\lineto[r]_(.2){-3}&\undertag{\bt}{-2}{6pt}\dashto[r]&\dashto[r]
&\undertag{\bt}{-2}{6pt}\lineto[r]
&\undertag{\bt}{-3}{6pt}\lineto[r]
&\undertag{\bt}{-(p+2)}{6pt}\\
&&\lineto[u]\\
&&\lineto[u]\\
&&\undertag{\bt}{-2}{6pt}\lineto[u]&&&}
$$
}
\end{exa}

\section{Appendix: Embeddings of graphs in $\frw \cup \farm\cup \frn$
into diagonal lattices}
\label{s:append}

In this final section we explicitely describe embeddings of the graphs
in $\frw \cup \farm\cup \frn$ into diagonal lattices of equal rank.
Formally, these embeddings verify one direction of the equivalence in 
Theorem~\ref{t:main} --- although this is the 'less interesting'
direction of the theorem.

\begin{figure}[hm]
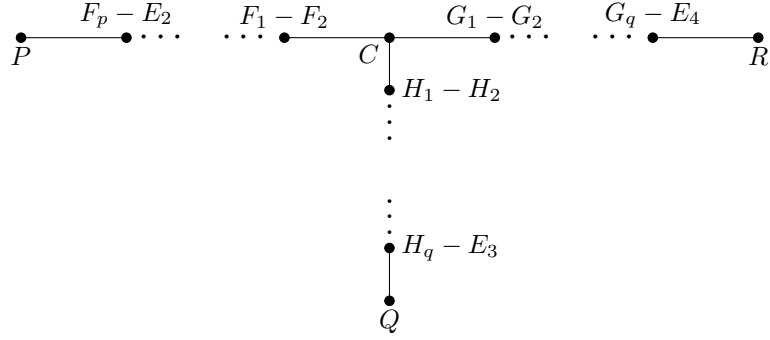

\begin{center}
\setlength{\unitlength}{1mm}
\unitlength=0.7cm
\begin{graph}(12,6)(-1,-5.5)
\graphnodesize{0.2}

 \roundnode{m1}(-2,0)
 \roundnode{m2}(0,0)
 \roundnode{m3}(3,0)
 \roundnode{m4}(5,0)
 \roundnode{m5}(7,0)
 \roundnode{m6}(10,0)
 \roundnode{m7}(12,0)
 \roundnode{m8}(5,-1)  
 \roundnode{m9}(5,-4)
 \roundnode{m10}(5,-5)

\edge{m1}{m2}
\edge{m3}{m4}
\edge{m4}{m5}
\edge{m6}{m7}
\edge{m4}{m8}
\edge{m9}{m10}

  \autonodetext{m1}[s]{{\small $P$}}
  \autonodetext{m2}[n]{{\small $F_p-E_2$}}
  \autonodetext{m3}[n]{{\small $F_1-F_2$}}
  \autonodetext{m4}[sw]{{\small $C$}}
  \autonodetext{m5}[n]{{\small $G_1-G_2$}}
  \autonodetext{m6}[n]{{\small $G_q-E_4$}}
  \autonodetext{m7}[s]{{\small $R$}}
  \autonodetext{m8}[e]{{\small $H_1-H_2$}}
  \autonodetext{m9}[e]{{\small $H_q-E_3$}}
  \autonodetext{m10}[s]{{\small $Q$}}
  \autonodetext{m2}[e]{{\Large $\cdots$}}
  \autonodetext{m3}[w]{{\Large $\cdots$}}
  \autonodetext{m5}[e]{{\Large $\cdots$}}
  \autonodetext{m6}[w]{{\Large $\cdots$}}

\freetext(5,-1.3){\Large $.$}
\freetext(5,-1.6){\Large $.$}
\freetext(5,-1.9){\Large $.$}
\freetext(5,-3.1){\Large $.$}
\freetext(5,-3.4){\Large $.$}
\freetext(5,-3.7){\Large $.$}

\end{graph}
\end{center}
\caption{\quad The embedding of the graph $\Gamma _{p,q,r}\in \frw$ into the
  diagonal lattice; here $C=E_1-F_1-G_1-H_1$,
  $P=E_2-E_1-E_3-H_p-\ldots - H_1, Q=E_3-E_1-E_4-G_q-\ldots - G_1,
  R=E_4-E_1-E_2-F_r-\ldots - F_1$ (where the basis of the diagonal lattice is
  given by $\{ E_1, \ldots , E_4, F_1, \ldots , F_r, G_1, \ldots , G_q, H_1,
  \ldots , H_p\}$)}
\label{f:wahltypeembed}
\end{figure}

\begin{figure}[hb]
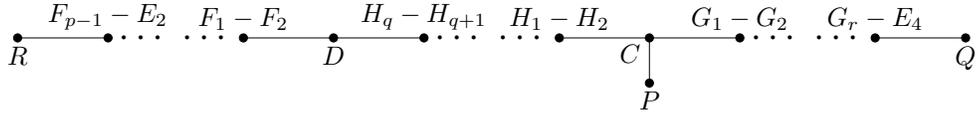

\begin{center}
\setlength{\unitlength}{1mm}
\unitlength=0.6cm
\begin{graph}(6,2)(-1.5,-1)
\graphnodesize{0.2}

 \roundnode{m12}(-9,0)
 \roundnode{m11}(-7,0)
 \roundnode{m10}(-4,0)
 \roundnode{m1}(-2,0)
 \roundnode{m2}(0,0)
 \roundnode{m3}(3,0)
 \roundnode{m4}(5,0)
 \roundnode{m5}(7,0)
 \roundnode{m6}(10,0)
 \roundnode{m7}(12,0)
 \roundnode{m8}(5,-1)  

\edge{m1}{m2}
\edge{m3}{m4}
\edge{m4}{m5}
\edge{m6}{m7}
\edge{m4}{m8}
\edge{m12}{m11}
\edge{m10}{m1}

  \autonodetext{m1}[s]{{\small $D$}}
  \autonodetext{m2}[n]{{\small $H_q-H_{q+1}$}}
  \autonodetext{m3}[n]{{\small $H_1-H_2$}}
  \autonodetext{m4}[sw]{{\small $C$}}
  \autonodetext{m5}[n]{{\small $G_1-G_2$}}
  \autonodetext{m6}[n]{{\small $G_r-E_4$}}
  \autonodetext{m7}[s]{{\small $Q$}}
  \autonodetext{m8}[s]{{\small $P$}}
 \autonodetext{m10}[n]{{\small $F_1-F_2$}}
 \autonodetext{m11}[n]{{\small $F_{p-1}-E_2$}}
 \autonodetext{m12}[s]{{\small $R$}} 
 \autonodetext{m2}[e]{{\Large $\cdots$}}
  \autonodetext{m3}[w]{{\Large $\cdots$}}
  \autonodetext{m5}[e]{{\Large $\cdots$}}
  \autonodetext{m6}[w]{{\Large $\cdots$}}
 \autonodetext{m10}[w]{{\Large $\cdots$}}
 \autonodetext{m11}[e]{{\Large $\cdots$}}

\end{graph}
\end{center}
\caption{\quad The embedding of the graph $\Delta _{p,q,r}$ in for
$p\geq 1$ and $q,r\geq 0$; here $C=E_1-G_1-H_1$, $D=H_{q+1}-E_3-F_1$,
$P=E_3-E_1-E_2-F_{p-1}-\ldots - F_1$, $Q=E_4-E_1-E_3-H_{q+1}-\ldots -
H_1$ and $R=E_2-E_1-E_4-G_r-\ldots - G_1$.  The basis in the diagonal
lattice is now $\{ E_1, \ldots , E_4, F_1, \ldots , F_{p-1}, G_1,
\ldots , G_r, H_1, \ldots , H_{q+1}\}$}
\label{f:masikbeagy}
\end{figure}

\begin{figure}[ht]
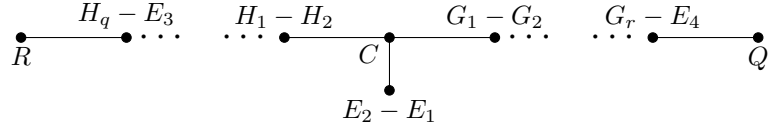

\begin{center}
\setlength{\unitlength}{1mm}
\unitlength=0.7cm
\begin{graph}(10,4)(0,-2)
\graphnodesize{0.2}

 \roundnode{m1}(-2,0)
 \roundnode{m2}(0,0)
 \roundnode{m3}(3,0)
 \roundnode{m4}(5,0)
 \roundnode{m5}(7,0)
 \roundnode{m6}(10,0)
 \roundnode{m7}(12,0)
 \roundnode{m8}(5,-1)  

\edge{m1}{m2}
\edge{m3}{m4}
\edge{m4}{m5}
\edge{m6}{m7}
\edge{m4}{m8}

  \autonodetext{m1}[s]{{\small $R$}}
  \autonodetext{m2}[n]{{\small $H_q-E_3$}}
  \autonodetext{m3}[n]{{\small $H_1-H_2$}}
  \autonodetext{m4}[sw]{{\small $C$}}
  \autonodetext{m5}[n]{{\small $G_1-G_2$}}
  \autonodetext{m6}[n]{{\small $G_r-E_4$}}
  \autonodetext{m7}[s]{{\small $Q$}}
  \autonodetext{m8}[s]{{\small $E_2-E_1$}}
  \autonodetext{m2}[e]{{\Large $\cdots$}}
  \autonodetext{m3}[w]{{\Large $\cdots$}}
  \autonodetext{m5}[e]{{\Large $\cdots$}}
  \autonodetext{m6}[w]{{\Large $\cdots$}}

\end{graph}
\end{center}
\caption{\quad The embedding of the graph $\Delta _{0,q,r}$;
$C=E_1-G_1-H_1$, $Q=E_4-E_1-E_2-E_4-H_q-\ldots -H_1$,
$R=E_3-E_1-E_2-E_4-G_r-\ldots -G_1$ (the basis in the diagonal lattice
is now $\{ E_1, \ldots , E_4, G_1, \ldots , G_r, H_1, \ldots ,
H_q\}$)}
\label{f:specmasikbeagy}
\end{figure}

\begin{figure}[ht]
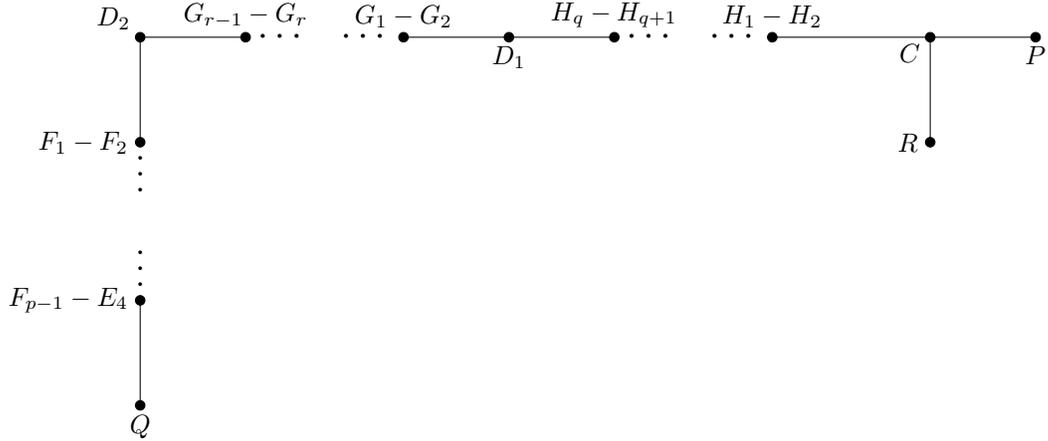

\begin{center}
\setlength{\unitlength}{1mm}
\unitlength=0.7cm
\begin{graph}(14,10)(-9,-10)
\graphnodesize{0.2}

 \roundnode{m15}(-10,-7)
 \roundnode{m14}(-10,-5)
 \roundnode{m13}(-10,-2)
 \roundnode{m12}(-10,0)
 \roundnode{m11}(-8,0)
 \roundnode{m10}(-5,0)
 \roundnode{m1}(-3,0)
 \roundnode{m2}(-1,0)
 \roundnode{m3}(2,0)
 \roundnode{m4}(5,0)
 \roundnode{m5}(7,0)
  \roundnode{m8}(5,-2)  

\edge{m1}{m2}
\edge{m3}{m4}
\edge{m4}{m5}
\edge{m4}{m8}
\edge{m12}{m11}
\edge{m10}{m1}
\edge{m15}{m14}
\edge{m13}{m12}

  \autonodetext{m1}[s]{{\small $D_1$}}
  \autonodetext{m2}[n]{{\small $H_q-H_{q+1}$}}
  \autonodetext{m3}[n]{{\small $H_1-H_2$}}
  \autonodetext{m4}[sw]{{\small $C$}}
  \autonodetext{m5}[s]{{\small $P$}}
  \autonodetext{m8}[w]{{\small $R$}}
 \autonodetext{m10}[n]{{\small $G_1-G_2$}}
 \autonodetext{m11}[n]{{\small $G_{r-1}-G_r$}}
 \autonodetext{m12}[nw]{{\small $D_2$}} 
  \autonodetext{m13}[w]{{\small $F_1-F_2$}}
 \autonodetext{m14}[w]{{\small $F_{p-1}-E_4$}}  
 \autonodetext{m15}[s]{{\small $Q$}} 
\autonodetext{m2}[e]{{\Large $\cdots$}}
\autonodetext{m3}[w]{{\Large $\cdots$}}
 \autonodetext{m10}[w]{{\Large $\cdots$}}
 \autonodetext{m11}[e]{{\Large $\cdots$}}
\freetext(-10,-2.3){\Large $.$}
\freetext(-10,-2.6){\Large $.$}
\freetext(-10,-2.9){\Large $.$}
\freetext(-10,-4.1){\Large $.$}
\freetext(-10,-4.4){\Large $.$}
\freetext(-10,-4.7){\Large $.$}

\end{graph}
\end{center}
\caption{\quad The embedding of the graph $\Lambda _{p,q,r}$ for
$p,r\geq 1$ and $q\geq 0$; here $C=E_1-H_1$, $D_1=H_{q+1}-G_1-E_3$,
$D_2=G_r-F_1-E_2$, $P=E_2-E_1-E_4-F_{p-1}-\ldots - F_1$,
$Q=E_4-E_1-E_3-H_{q+1}-\ldots - H_1$, $R=E_3-E_1-E_2-G_r-\ldots -
G_1$}
\label{f:lambdabeagy}
\end{figure}

\begin{figure}[ht]
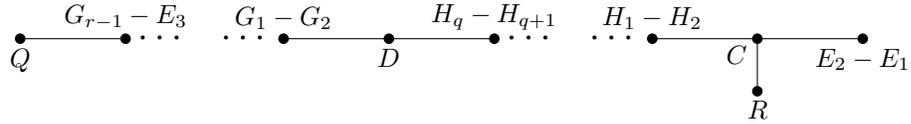

\begin{center}
\setlength{\unitlength}{1mm}
\unitlength=0.7cm
\begin{graph}(10,6)(-6,-5)
\graphnodesize{0.2}

 \roundnode{m12}(-9,0)
 \roundnode{m11}(-7,0)
 \roundnode{m10}(-4,0)
 \roundnode{m1}(-2,0)
 \roundnode{m2}(0,0)
 \roundnode{m3}(3,0)
 \roundnode{m4}(5,0)
 \roundnode{m5}(7,0)
  \roundnode{m8}(5,-1)  

\edge{m1}{m2}
\edge{m3}{m4}
\edge{m4}{m5}
\edge{m4}{m8}
\edge{m12}{m11}
\edge{m10}{m1}

  \autonodetext{m1}[s]{{\small $D$}}
  \autonodetext{m2}[n]{{\small $H_q-H_{q+1}$}}
  \autonodetext{m3}[n]{{\small $H_1-H_2$}}
  \autonodetext{m4}[sw]{{\small $C$}}
  \autonodetext{m5}[s]{{\small $E_2-E_1$}}
  \autonodetext{m8}[s]{{\small $R$}}
 \autonodetext{m10}[n]{{\small $G_1-G_2$}}
 \autonodetext{m11}[n]{{\small $G_{r-1}-E_3$}}
 \autonodetext{m12}[s]{{\small $Q$}} 
\autonodetext{m2}[e]{{\Large $\cdots$}}
\autonodetext{m3}[w]{{\Large $\cdots$}}
 \autonodetext{m10}[w]{{\Large $\cdots$}}
 \autonodetext{m11}[e]{{\Large $\cdots$}}

\end{graph}
\end{center}
\caption{\quad The embedding of the graph $\Lambda _{0,q,r}$ for
$r\geq 1$ and $q\geq 0$; $C=E_1-H_1$, $D=H_{q+1}-E_3-G_1$,
$Q=E_3-E_1-E_4-E_2-H_{q+1}-\ldots -H_1$,
$R=E_4-E_1-E_3-E_2-G_{r-1}-\ldots - G_1$}
\label{f:lambdaspec1beagy}
\end{figure}

\begin{figure}[ht]
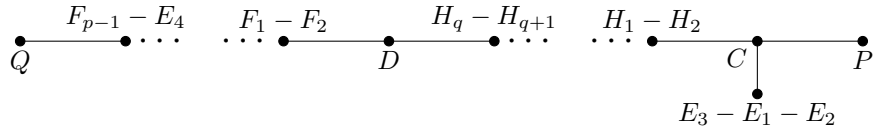

\begin{center}
\setlength{\unitlength}{1mm}
\unitlength=0.7cm
\begin{graph}(10,6)(-6,-5)
\graphnodesize{0.2}

 \roundnode{m12}(-9,0)
 \roundnode{m11}(-7,0)
 \roundnode{m10}(-4,0)
 \roundnode{m1}(-2,0)
 \roundnode{m2}(0,0)
 \roundnode{m3}(3,0)
 \roundnode{m4}(5,0)
 \roundnode{m5}(7,0)
  \roundnode{m8}(5,-1)  

\edge{m1}{m2}
\edge{m3}{m4}
\edge{m4}{m5}
\edge{m4}{m8}
\edge{m12}{m11}
\edge{m10}{m1}

  \autonodetext{m1}[s]{{\small $D$}}
  \autonodetext{m2}[n]{{\small $H_q-H_{q+1}$}}
  \autonodetext{m3}[n]{{\small $H_1-H_2$}}
  \autonodetext{m4}[sw]{{\small $C$}}
  \autonodetext{m5}[s]{{\small $P$}}
  \autonodetext{m8}[s]{{\small $E_3-E_1-E_2$}}
 \autonodetext{m10}[n]{{\small $F_1-F_2$}}
 \autonodetext{m11}[n]{{\small $F_{p-1}-E_4$}}
 \autonodetext{m12}[s]{{\small $Q$}} 
\autonodetext{m2}[e]{{\Large $\cdots$}}
\autonodetext{m3}[w]{{\Large $\cdots$}}
 \autonodetext{m10}[w]{{\Large $\cdots$}}
 \autonodetext{m11}[e]{{\Large $\cdots$}}
 
\end{graph}
\end{center}
\caption{\quad The embedding of the graph $\Lambda _{p,q,0}$ for
$p\geq 1$ and $q\geq 0$; $C=E_1-H_1$, $D=H_{q+1}-E_2-E_3-F_1$,
$P=E_2-E_1-E_4-F_{p-1}-\ldots - F_1$, 
$Q=E_4-E_1-E_3- H_{q+1}-\ldots - H_1$}
\label{f:lambdaspec2beagy}
\end{figure}

\FloatBarrier


\begin{thebibliography}{AAA}
\bibitem{CH} A. Casson and J. Harer, 
{\it Some homology lens spaces
which bound rational homology balls}, Pacific J. Math. {\bf96} (1981),
23--36.



\bibitem{D1}
S. Donaldson,
{\it Irrationality and the h-cobordism conjecture}, J. Diff. Geom.
{\bf 26} (1987), 141--168.

\bibitem{Donor}
S. Donaldson,
{\it The orientation of Yang--Mills moduli spaces and 4--manifold topology},
J. Differential Geometry {\bf26} (1987), 397--428.


\bibitem{FS2} 
R. Fintushel and R. Stern,
{\it Immersed spheres in 4-manifolds and the immersed Thom Conjecture},
Turkish J. Math. {\bf 19} (1995), 145--157.

\bibitem{FS1} 
R. Fintushel and R. Stern,
{\it Rational blowdowns of smooth 4--manifolds}, J. Diff. Geom. {\bf46}
(1997), 181--235.

\bibitem{FSexo}
R. Fintushel and R. Stern,
{\it Double node neighborhoods and families of simply connected 4--manifolds
with $b^+_2=1$}, J. Amer. Math. Soc.  {\bf19}  (2006), 171--180.

\bibitem{GaS}
D. Gay and A. Stipsicz,
{\it Symplectic rational blow--down along Seifert fibered 3--manifolds},
Int. Math. Res. Notices, 2007 Volume {\bf2007}, article ID rnm084.

\bibitem{GSuj}
D. Gay and A. Stipsicz,
{\it Symplectic surgeries and normal surface singularities},
arXiv:0708.1417


\bibitem{gh} G.--M. Greuel and H. Hamm, {\it Invarianten quasihomogener
    vollst\"andiger Durchschnitte}, Invent. Math. {\bf 49} (1978), 67--86.

\bibitem{hirz} F. Hirzebruch, W. Neumann, and S. Koh, 
{\it Differentiable Manifolds and Quadratic Forms}, 
Lecture Notes in Pure and Applied Math. {\bf 4}, Marcel Dekker, New York, 1971.

\bibitem{kollar}
J. Koll\'ar,
{\it Is there a topological Bogomolov--Miyaoka--Yau inequality?}
arXiv:math.AG/0602562.

\bibitem{laufert} H. Laufer, 
{\it Taut two-dimensional singularities},
 Math. Ann. {\bf 205} (1973), 131--164.

\bibitem{laufermin} H. Laufer, 
{\it On minimally elliptic singularities}, Amer. J. Math.
{\bf99} (1977),  1257--1295.



\bibitem{LPark}
Y. Lee and J. Park, 
{\it Simply connected surfaces of general type with $p_g=0$ and $K^2=2$},
arXiv:math.AG/0609072.


\bibitem{LW}
E. Looijenga and J. Wahl,
{\it Quadratic functions and smoothing surface singularities}, Topology  
{\bf25}  (1986), 261--291. 



\bibitem{nemethi}
A. N\'emethi,
{\it On the Ozsv\'ath-Szab\'o invariant of negative definite 
plumbed 3-manifolds},
Geom. Topol.  {\bf9}  (2005), 991--1042.
      



\bibitem{raym} W. Neumann and F. Raymond,
{\it Seifert manifolds, plumbing, $\mu$-invariant and orientation 
reversing maps}, Alg. and Geom. Topology, Proc. Santa Barbara 1977, 
Lecture Notes in Math. {\bf 664} (Springer, Berlin, 1978), 163-196.




\bibitem{nw1} W. Neumann and J. Wahl,
{\it Complete intersection singularities of splice
 type as universal abelian covers}, Geom. Topol. {\bf 9} (2005), 699--755.


\bibitem{OW}
P. Orlik and P. Wagreich, 
{\it Algebraic surfaces with $k^*$--action},
Acta Math. {\bf138} (1977), 43--81.

 \bibitem{O-W} P. Orlik and P. Wagreich, {\it Isolated singularities
of algebraic surfaces with $\mathbb C ^{*}$ action}, Ann. of
Math. {\bf 93} (1971), 205--228.

\bibitem{Oszabs} P. Ozsv\'ath and Z. Szab\'o, 
{\it Absolutely graded Floer homologies and intersections forms
for four--manifolds with boundary},
Adv. Math. {\bf173} (2003), 179--261.



\bibitem{upa}
H. Park, J. Park and D. Shin: {\it A simply connected
surface of general type with $p_g=0$ and $K^3=3$},
arXiv:0708.0273


\bibitem{Pratb} 
J. Park,
{\it Seiberg--Witten invariants of generalized rational blow--downs},
Bull. Austral. Math. Soc. {\bf56} (1997), 363--384.

\bibitem{P} 
J. Park, 
{\it Simply connected symplectic 4--manifolds with $b_2^+=1$ and $c^2_1=2$},
Invent. Math. {\bf159} (2005), 657--667.

\bibitem{PSS}
J. Park, A. Stipsicz and Z. Szab\'o,
{\it Exotic smooth structures on $\cpot$},  
Math. Res. Lett.  {\bf12}  (2005), 701--712. 
 	

\bibitem{pink} H. Pinkham, 
{\it Deformations of algebraic varieties with
$\mathbb G_{m}$ action},  Ast\'{e}risque {\bf 20} (1974), 1-131.

\bibitem{pinkham} H. Pinkham, 
{\it Normal surface singularities with
$\C ^{*}$ action},  Math. Ann. {\bf 227} (1977), 183--193.

\bibitem{pinkdef} 
H. Pinkham, 
{\it Deformations of normal surface
singularities with $\C^{*}$-action},  Math. Ann. {\bf 232} (1978),
65--84.

\bibitem{SPAMS} 
A. Stipsicz, 
{\it On the ${\overline {\mu}}$ invariant of rational surface
singularities} Proc. Amer. Math. Soc., to appear.

\bibitem{SS}
A. Stipsicz and Z. Szab\'o,
{\it An exotic smooth structure on $\cphat$},
Geom. Topol. {\bf9} (2005), 813--832.

\bibitem{Sym0}
M. Symington,
{\it Symplectic rational blowdowns}, J. Diff. Geom. {\bf50} (1998),
505--518.

\bibitem{Sym} 
M. Symington,
{\it Generalized symplectic rational blowdowns}, 
Algebr. Geom. Topol. {\bf1} (2001), 503--518.


\bibitem{wahl1}
J. Wahl,
{\it Elliptic deformations of minimally elliptic singularities}, Math. Ann. {\bf253} (1980), 241-262.

\bibitem{Wahl}
J. Wahl,
{\it Smoothings of normal surface singularities},
Topology {\bf20} (1981), 219--246.

\bibitem{wolf} J. Wolf, \emph{Spaces of constant curvature},
  McGraw-Hill, New York-London-Sydney (1967).


\end{thebibliography}
\end{document}